\documentclass[preprint,review,12pt]{elsarticle}

\usepackage{amssymb}

 \usepackage{lineno}

\usepackage{amsmath,hyperref}
\usepackage{mathrsfs}
\usepackage{mathtools}
\usepackage[margin=1in]{geometry}
\usepackage{multirow}
\usepackage{subfig}
\biboptions{numbers,sort&compress}

\usepackage{subfig}

\numberwithin{equation}{section}

\newcommand{\pStokeslet}{\mathbb{G}}
\newcommand{\fStokeslet}{\mathbb{G}}
\newcommand{\rStokeslet}{\mathbb{G}^{(r)}}
\newcommand{\wStokeslet}{\mathbb{G}^{(w)}}
\newcommand{\mbf}[1]{\mathbf{#1}}
\newcommand{\vf}[1]{\boldsymbol{#1}}
\newcommand{\op}[1]{\boldsymbol{\mathcal{#1}}}
\let\div\undefined
\newcommand{\div}{\boldsymbol{\nabla} \cdot}
\newcommand{\grad}{\boldsymbol{\nabla}}
\newcommand{\lapl}{\boldsymbol{\nabla}^2}
\newcommand*\diff{\,\mathrm{d}}
\newcommand{\ie}{i.e.}
\newcommand{\vslip}{\breve{\vf{v}}}
\newcommand{\bfslip}{\breve{\mbf{v}}}
\newcommand{\mbfhalf}[1]{\mbf{#1}^{1/2}}
\newcommand{\vfhalf}[1]{\vf{#1}^{\frac{1}{2}}}
\newcommand{\ophalf}[1]{\op{#1}^{\frac{1}{2}}}
\newcommand{\stochdrift}[2]{\partial_{#1} \cdot #2}

\newcommand{\Mwave}{\mbf{M}^{(w)}}
\newcommand{\Mreal}{\mbf{M}^{(r)}}
\newcommand{\sqrtM}{\mbf{M}^{1/2}}
\newcommand{\sqrtMwave}{\left(\mbf{M}^{(w)}\right)^{1/2}}

\newcommand{\xir}{\xi^2 r^2}
\newcommand{\kxi}{k^2/4\xi^2}

\newcommand{\kbt}{k_BT}
\newcommand{\rcutoff}{r_c}
\newcommand{\rAlpert}{r_{\text{Alpert}}}
\newcommand{\Nbox}{N_{\text{box}}}
\newcommand{\tolEwald}{\epsilon_{\text{Ewald}}}
\newcommand{\tolIter}{\epsilon_{\text{tol}}}

\usepackage{color}
\PassOptionsToPackage{normalem}{ulem}
\usepackage{ulem}

\newcommand{\delete}[1]{}

\definecolor{mygreen}{rgb}{0,.6,0}

\allowdisplaybreaks

\usepackage[nameinlink,noabbrev]{cleveref}
\crefname{equation}{eq.}{eqs.} 
\Crefname{equation}{Eq.}{Eqs.}

\journal{J. Comput. Phys.}

\begin{document}

\begin{frontmatter}



\title{A Fluctuating Boundary Integral Method for Brownian Suspensions}

\author[cims]{Yuanxun Bao\corref{cor1}}
\ead{billbao@cims.nyu.edu}

\author[yale]{Manas Rachh}
\ead{manas.rachh@yale.edu}

\author[ic]{Eric E.~Keaveny}
\ead{e.keaveny@imperial.ac.uk}

\author[cims,flatiron]{Leslie Greengard}
\ead{greengard@cims.nyu.edu}

\author[cims]{Aleksandar Donev}
\ead{donev@courant.nyu.edu}

\address[cims]{Courant Institute of Mathematical Sciences, New York University, 251 Mercer Street, New York, NY, USA}
\address[yale]{Applied Mathematics Program, Yale University, New Haven, CT 06511}
\address[ic]{Department of Mathematics, Imperial College London, London SW7 2AZ, United Kingdom}
\address[flatiron]{Center for Computational Biology, Flatiron Institute, New York, NY}

\cortext[cor1]{Corresponding author}



\begin{abstract}
We present a fluctuating boundary integral method (FBIM) for overdamped Brownian Dynamics (BD) of two-dimensional periodic suspensions of rigid particles of complex shape immersed in a Stokes fluid. We develop a novel approach for generating Brownian displacements that arise in response to the thermal fluctuations in the fluid. Our approach relies on a first-kind boundary integral formulation of a mobility problem in which a random surface velocity is prescribed on the particle surface, with zero mean and covariance proportional to the Green's function for Stokes flow (Stokeslet). This approach yields an algorithm that scales linearly in the number of particles for both deterministic and stochastic dynamics, handles particles of complex shape, achieves high order of accuracy, and can be generalized to three dimensions and other boundary conditions. We show that Brownian displacements generated by our method obey the discrete fluctuation-dissipation balance relation (DFDB). Based on a recently-developed Positively Split Ewald method [A. M. Fiore, F. Balboa Usabiaga, A. Donev and J. W. Swan, J. Chem. Phys., 146, 124116, 2017], near-field contributions to the Brownian displacements are efficiently approximated by iterative methods in real space, while far-field contributions are rapidly generated by fast Fourier-space methods based on fluctuating hydrodynamics. FBIM provides the key ingredient for time integration of the overdamped Langevin equations for Brownian suspensions of rigid particles. We demonstrate that FBIM obeys DFDB by performing equilibrium BD simulations of suspensions of starfish-shaped bodies using a random finite difference temporal integrator.

\end{abstract}


\end{frontmatter}


\section{Introduction}

Complex fluids containing colloidal particles are ubiquitous in science and industrial applications. Colloidal particles span length scales from {several nanometers}, such as magnetic nano-propellers \cite{MagneticNanoPropeller} and molecular motors \cite{BrownainMotor_Peskin,BioNanoMotors_Simulation}, to a few microns, such as self-phoretic Janus particles \cite{JanusParticle} and motile microorganisms \cite{ActiveSuspensions}. 
In the last decade, increasing attention has been given to the emerging field of {\it active} colloidal suspensions \cite{Hematites_Science, CatalyticNanomotors, FlippingNanorods, Self-MotileColloidalParticles, ActiveNanoSwimmersReview}, in which particles move autonomously or in response to external forces.
Despite the advances in the theory and experimental design of passive and active colloids, developing accurate and  efficient computational methods that are capable of simulating tens or hundreds of thousands of particles, as well as handling particles of complex shape, still remains a formidable challenge. 
Here we develop a novel algorithm for generating the Brownian (stochastic) displacements required to perform overdamped Brownian dynamics of a suspension of rigid particles immersed in a Stokes fluid. Our method is based on boundary integral techniques, scales linearly in the number of particles for both deterministic and stochastic dynamics, handles particles of complex shape, and achieves high order accuracy. 
Because of its close connection to fluctuating hydrodynamics, we refer to our method as the {\it Fluctuating Boundary Integral Method} (FBIM).
We restrict our attention to two-dimensional periodic domains. However, our approach can be extended to three dimensions and confined suspensions.

The two key ingredients that need to be included in a computational method for colloidal suspensions are the long-ranged hydrodynamic interactions (HIs) and the correlated Brownian motion of the particles. 
In the absence of active and Brownian motion, describing the hydrodynamics of Stokesian suspensions requires the accurate solution of mobility problems \cite{MicrohydrodynamicsBook,BoundaryIntegral_Pozrikidis}, \ie, computing the linear and angular velocities of the particles in response to applied (external) forces and torques. This process defines the action of a {\it body mobility matrix}, which converts the applied forces to the resulting particle motions. The mobility matrix encodes the many-body hydrodynamic interactions among the particles, and is therefore dense with long-ranged interactions between \emph{all} bodies. In the presence of Brownian motion, fluctuation-dissipation balance requires the Brownian displacements to have zero mean and covariance proportional to the hydrodynamic mobility matrix. This necessitates an algorithm for generating Gaussian random variables with covariance equal to the body mobility matrix. In this work, we present a new linear-scaling boundary integral method to generate, together, both the deterministic and Brownian displacements.


For passive suspensions of spherical particles {in zero Reynolds number flows (infinite Schmidt number)}, the methods of  Brownian \cite{BrownianDynamics_DNA,BrownianDynamics_DNA2, BrownianDynamics_OrderN,BD_IBM_Graham,RPY_FMM} and Stokesian Dynamics (SD) \cite{BrownianDynamics_OrderNlogN, StokesianDynamics_Brownian,SD_SpectralEwald} have dominated the chemical engineering community. These methods are tailored to sphere suspensions and utilize a multipole hierarchy truncated at either the monopole (BD) or dipole (SD) level in order to capture the far-field behavior of the hydrodynamic interactions. Modern fast algorithms can apply the action of the truncated mobility matrix with linear-scaling  by using the Fast Multipole Method (FMM) for an unbounded domain \cite{RPY_FMM}, and using Ewald-like methods for periodic \cite{SpectralEwald_Stokes,SD_SpectralEwald} and confined domains \cite{BrownianDynamics_OrderN}. The Brownian (stochastic) displacements are typically generated iteratively by a Chebyshev polynomial approximation method \cite{BD_Fixman_sqrtM}, or by the Lanczos algorithm for application of the matrix square root \cite{SquareRootKrylov}. However, since the hydrodynamic interactions among particles decay slowly like the inverse of distance {in three dimensions and diverge logarithmically in two dimensions}, the condition number of the mobility matrix grows as the number of particles increases (keeping the packing fraction fixed, see \cite[Fig. 1]{SpectralRPY}). Therefore, the overall computational scaling for {generating Brownian displacements using iterative methods} is only superlinear in general.

{The fluctuating Lattice Boltzmann (FLB) method has been used for Brownian suspensions for some time \cite{LB_ThermalFluctuations,LB_SoftMatter_Review}. This is an explicit solvent method which includes fluid inertia and thus operates at finite Schmidt number instead of in the overdamped limit we are interested in; furthermore, FLB relies on artificial fluid compressibility to avoid solving Poisson problems for the pressure. While the cost of each time step is linear in the number of particles (more precisely, the number of fluid cells) $N$, in three dimensions $\mathcal{O}\left( N^{2/3} \right)$ time steps are required for vorticity to diffuse throughout the system volume, leading to superlinear $\mathcal{O}\left(N^{5/3}\right)$ overall complexity \cite{SpectralRPY}.}

As an alternative, methods such as the Fluctuating Immersed Boundary method (FIB) \cite{BrownianBlobs} and the fluctuating Force Coupling Method (FCM) \cite{ForceCoupling_Fluctuations,FluctuatingFCM_DC}  
utilize fluctuating hydrodynamics to generate the Brownian increments in linear time by solving the fluctuating \emph{steady} Stokes equation on a grid.
This ensures that the computational cost of Brownian simulation is only marginally larger than the cost of deterministic simulations, in stark contrast to traditional BD approaches. 
The FIB/FCM approach to generating the Brownian displacements is further improved in the recently-developed Positively Split Ewald (PSE) method \cite{SpectralRPY}. In PSE, the Rotne-Prager-Yamakawa (RPY) tensor \cite{RotnePrager} is used to capture the long-ranged hydrodynamic interactions, and its action is computed with spectral accuracy by extending the Spectral Ewald \cite{SpectralEwald_Stokes} method for the RPY tensor. 
The key idea in PSE is to use the Hasimoto splitting \cite{Mobility2D_Hasimoto} to decompose the RPY tensor into near-field (short-ranged) and far-field (long-ranged) contributions, in a way that ensures that {\em both} contributions are independently {\it symmetric positive definite} (SPD).
This makes it possible to apply a Lanczos algorithm \cite{SquareRootKrylov} to generate the near-field   
contribution with only a small ($\mathcal{O}(1)$) number of iterations, while the far-field contribution is computed by fast Fourier-space methods based on fluctuating hydrodynamics using only a few FFTs. Later in \autoref{sec:sqrtM}, we will apply the same SPD splitting idea to the Green's function of steady Stokes flow to achieve linear scaling in the FBIM.

Essentially all commonly-used methods for overdamped Brownian suspension flows are limited to spherical particles (with some extensions to spheroids \cite{ForceCoupling_Ellisoids}), and generalization to include particles of complex shape is generally difficult. Further, these methods employ an {\it uncontrolled} truncation of a multipole expansion hierarchy and therefore become inaccurate when particles get close to one another as in dense suspensions.

For deterministic Stokes problems, the Boundary Integral Method (BIM) \cite{BoundaryIntegral_Pozrikidis} is very well-developed \cite{Stokes3D_Greengard,Stokes2D_Manas,Stokes3D_FMM} and allows one to handle complex particle shapes and achieve {\em controlled accuracy} even for dense suspensions \cite{BoundaryIntegral_Periodic3D,BoundaryIntegral_SpheroidQBX}.
In the boundary integral framework, the steady Stokes equation are reformulated as an integral equation of unknown densities that are defined on the boundary, using a first-kind (single-layer densities) or second-kind (double-layer densities) formulation, or a mixture of both.
Suspended particles of complex geometry can be directly discretized by a surface mesh, and, by a suitable choice of surface quadrature, higher-order (or even spectral) accuracy can be achieved.
The key difficulty is handling the singularity of the Green's functions appearing in the boundary integral formulation.
One possibility is to regularize the singularity, as done in the method of regularized Stokeslets \cite{RigidRegularizedStokeslets} and the recently-developed linear-scaling rigid multiblob method \cite{RigidMultiblobs}, both of which rely on a regularized first-kind formulation.
Regularization, however, comes at a drastic loss of accuracy, and to resolve near-field hydrodynamic interactions accurately one must {make use of} high-order accurate
{\em singular quadratures} for the singular or near-singular kernel in the first- and second-kind integral operators, {such as} quadrature-by-expansion (QBX) \cite{QBX_Epstein,QBX_Klockner}.

Discretizing the boundary integral equation typically leads to a dense linear system, and fast algorithms for performing the dense matrix-vector product are required to achieve linear-scaling, such as the Fast Multipole Method (FMM) \cite{Stokes3D_Greengard, Stokes3D_FMM}, and Spectral Ewald methods \cite{BoundaryIntegral_Periodic3D,BoundaryIntegral_SpheroidQBX}. 
Much of the state-of-the-art BIM work on Stokes mobility problems \cite{BoundaryIntegral_Periodic3D,BoundaryIntegral_SpheroidQBX,Stokes3D_FMM,StokesSecondKind_Shelley} uses (completed) {\em second-kind} (double layer) formulations, and we will not review it in detail here since we will rely on a first-kind formulation.
The current development of BIM for {Stokesian} suspensions is limited to the deterministic case only. 
In this work, we combine singular quadrature (namely, Alpert quadrature in two dimensions) and the idea of SPD splitting from the PSE method to develop a linear-scaling method for generating the Brownian displacements based on a {\em first-kind} boundary integral formulation.
{The key idea is that, instead of adding a stochastic stress to the fluid equations as done in fluctuating hydrodynamics, we prescribe a random surface velocity (distribution) on the particle surface that has zero mean and covariance proportional to the (singular) periodic Green's function of the Stokes equations (the periodic {\it Stokeslet}). Reformulating the resulting stochastic boundary value problem as a first-kind boundary integral equation reveals that the random surface velocity has covariance proportional to the single-layer integral operator, which suggests that one ought to handle the singularity of the covariance using the same machinery used to handle the singularity of the Green's function in the first-kind BIM. This allows us to develop a numerical method that satisfies discrete fluctuation-dissipation balance (DFDB) to within solver and roundoff errors.}

{The outline of this paper is as follows. In the continuum formulation (\autoref{sec:continuum-formulation}), we show that computing the action of the hydrodynamic mobility matrix and its square root can be formulated equivalently as the solution of a Stokes boundary value problem, herein referred to as the {\it Stochastic Stokes Boundary Value Problem} (SSBVP). Reformulating the SSBVP as a first-kind boundary integral equation allows us to develop a numerical method (\autoref{sec:FBIM-method}) that satisfies discrete fluctuation-dissipation balance (DFDB) to within solver and roundoff errors.   We apply the FBIM to several benchmark problems in two dimensions (\autoref{sec:FBIM-results}), and assess the effectiveness of the method by its accuracy and convergence, robustness, and scalability for suspensions of many rigid disks. We also couple the FBIM with stochastic temporal integrators based on the idea of random finite differences \cite{BrownianBlobs,BrownianMultiBlobs,MagneticRollers} to perform BD simulations with starfish-shaped particles, and confirm that DFBD is obtained for sufficiently small time step sizes by comparing the numerical equilibrium distribution to the correct Gibbs-Boltzmann distribution. }

\section{Continuum formulation}
\label{sec:continuum-formulation}
This section presents the continuum formulation of the equations of motion for Brownian suspensions. We consider a suspension of $N$ rigid Brownian particles of complex shape immersed in a viscous incompressible fluid with constant density $\rho$, viscosity $\eta$, and temperature $T$ in a domain $\mathcal{V}$ with periodic boundary conditions. Each particle or body, indexed by $\beta$, is described by the position of a chosen ``tracking point'', denoted by $\vf{q}_{\beta}$, and its rotation relative to a chosen reference configuration, denoted by $\vf{\theta}_{\beta}$. 
{In this paper, we restrict the formulation and implementation to two spatial dimensions only, and $\theta_{\beta} \in \mathbb{R}$ is simply the angle of rotation. {The generalization of the continuum formulation} to three dimensional systems is conceptually straightforward but technically complicated by the fact that orientation in three dimensions cannot be represented by a vector in $\mathbb{R}^3$; instead, one can use normalized quaternions \cite{BrownianMultiBlobs}.} We denote the particle surface by $\Gamma_{\beta}$ which encloses an interior domain denoted by $D_{\beta}$. 
The fluid domain exterior to the particles is defined by $ E = \mathcal{V} \backslash \left \{ \cup_{\beta=1}^{N} \overline{D}_{\beta}  \right\} $, where $\overline{D}_{\beta} = D_{\beta} \cup \Gamma_{\beta}$. 

{A number of prior works \cite{LLNS_FD_Fox,VACF_Langevin,LangevinDynamics_Theory,VACF_FluctHydro}
have arrived at a fluctuating hydrodynamics formulation of the equations of motion
when fluid and particle inertia (and perhaps compressibility \cite{BrownianCompressibility_Zwanzig}) is accounted for\footnote{Note that Hauge and Martin-Lof \cite{VACF_FluctHydro} explain that there is some ambiguity in whether the stochastic traction
is taken to be zero or nonzero on the particle surface; this choice does not, however, affect the resulting equations of motion for the bodies.
We consider the formulation given here to be the more physically meaningful and
follow Hinch \cite{VACF_Langevin}, see in particular Section 3 of the comprehensive work of Roux \cite{LangevinDynamics_Theory}.}.}
The fluid velocity $\vf{v}(\vf{r},t)$ and the pressure $\pi$ 
follows the time-dependent fluctuating Stokes equations for all $\vf{r} \in E$,
\begin{subequations}
\begin{align}
\rho \, \partial_t \vf{v} + \grad \pi &= \eta \lapl \vf{v} + \sqrt{2\eta k_B T} ~ \div \op{Z}, \\
\div \vf{v} &= 0,
\end{align}
\end{subequations}
where $k_B$ is Boltzmann's constant and  $\op{Z}(\vf{r},t)$ is a random Gaussian tensor field whose components are white in space and time with mean zero and covariance
\begin{equation}
 \left \langle \op{Z}_{ij}(\vf{r},t) ~ \op{Z}_{kl}(\vf{r}',t) \right \rangle  = (\delta_{ik} \delta_{jl} + \delta_{il} \delta_{jk}) \delta(\vf{r}-\vf{r}') \delta(t-t'). 
\label{StochStressCov}
\end{equation}
On the surface of $\Gamma_{\beta}$, we assume that \textit{no-slip} boundary condition (BC) holds,
\begin{equation}
\vf{v}(\vf{x},t) = \vf{u}_{\beta}+ \vf{\omega}_{\beta}  \times (\vf{x} -\vf{q}_{\beta}),  \quad \forall  \vf{x} \in \Gamma_{\beta}, 
\label{noslipBC}
\end{equation}
where $\vf{u}_{\beta}$ and $\vf{\omega}_{\beta}$ are the translational and rotational velocities of the particle. Each particle is also subject to applied force $\vf{f}_{\beta}$ and torque $\vf{\tau}_{\beta}$, which are related to the fluid stress tensor $\vf{\sigma}$ and stochastic stress tensor $\vf{\sigma}^{(s)}$ by 
\begin{equation}
\begin{aligned}
m_{\beta} \frac{d \vf{u}_{\beta}}{dt} &= \vf{f}_{\beta} - \int_{\Gamma_{\beta}} ( \vf{\lambda}_{\beta} + \vf{\lambda}_{\beta}^{(s)} ) (\vf{x}) \diff S_{\vf{x}}, \\
\mbf{I}_{\beta} \cdot \frac{d \vf{\omega}_{\beta}}{dt} &= \vf{\tau}_{\beta} - \int_{\Gamma_\beta}  (\vf{x}-\vf{x}_{\beta}) \times ( \vf{\lambda}_{\beta} + \vf{\lambda}_{\beta}^{(s)}  ) (\vf{x}) \diff S_{\vf{x}},
\end{aligned}
\label{FTbalance}
\end{equation}
where $m_{\beta}$ and $\mbf{I}_{\beta}$ are mass and moment of inertia tensor of body $\beta$. Here, $\vf{\lambda}_{\beta}(\vf{x}) = (\vf{\sigma} \cdot \vf{n}_{\beta})(\vf{x})$  and $\vf{\lambda}_{\beta}^{(s)}(\vf{x}) = \left( \vf{\sigma}^{(s)} \cdot \vf{n}_{\beta}  \right)(\vf{x})$ are the normal components of the stress tensors on the outside of the surface of the body, 
\begin{equation}
\begin{aligned}
 \vf{\sigma} &= -\pi \vf{I} + \eta (\grad \vf{v}  + \grad^{\top}{\vf{v}} ), \\
 \vf{\sigma}^{(s)} &=  \sqrt{2\eta \kbt } \op{Z},  
\end{aligned}
\label{fluidstress}
\end{equation}
and $\vf{n}_{\beta}$ is the unit outward normal vector of $\Gamma_{\beta}$ pointing into the fluid domain.

In the overdamped (large Schmidt number) limit where the fluid velocity is eliminated as a fast variable through an adiabatic mode elimination procedure \cite{AveragingHomogenization},  the diffusive motion of the rigid bodies can be described by the stochastic differential equation (SDE) of {\it Brownian Dynamics} (BD)\footnote{We use the differential notation of SDEs common in the physics literature.} \cite{LangevinDynamics_Theory}:
\begin{equation}
 \frac{d \vf{Q}}{dt} = \op{N} \vf{F} + \sqrt{2\kbt} \vfhalf{\op{N}} \op{W} + (k_B T) (\stochdrift{\vf{Q}}{\op{N}}) , \label{BDeqns}
\end{equation}
where 
$\vf{Q}_{\beta} = \left\{ \vf{q}_{\beta}, \vf{\theta}_{\beta} \right\}$ and $\vf{Q} = \left\{ \vf{Q}_{\beta} \right\}_{\beta = 1}^{N}$ is a composite vector that collects the positions and orientations of  particles (in two dimensions, $\vf{Q}_{\beta} \in \mathbb{R}^3$). 
The first term on the right-hand-side of \Cref{BDeqns}  is the deterministic motion of rigid bodies, where  $\op{N}(\vf{Q}) \succeq \vf{0}$ is the symmetric positive semi-definite {\it body mobility} matrix that converts applied forces and torques $\vf{F}(\vf{Q}) =  \left\{  \vf{f}_{\beta}(\vf{Q}), \vf{\tau}_{\beta}(\vf{Q}) \right\}_{\beta=1}^{N}$ to rigid body motion, and can be obtained by solving the standard mobility problem for rigid body motion, 
\begin{equation}
\left\{
\begin{aligned}
-\div \vf{\sigma} &= \grad \pi - \eta \lapl \vf{v} = \vf{0}, \\
\div \vf{v} &= 0, \\
\vf{v}(\vf{x}) &= \vf{u}_{\beta} + \vf{\omega}_{\beta} \times (\vf{x}-\vf{q}_{\beta}), \quad \forall \, \vf{x} \in \Gamma_{\beta}, \\
\int_{\Gamma_{\beta}} \vf{\lambda}_{\beta}(\vf{x}) \diff S_{\vf{x}} &= \vf{f}_{\beta} \quad \text{and} \quad \int_{\Gamma_{\beta}} (\vf{x}-\vf{x}_c) \times \vf{\lambda}_{\beta}(\vf{x})\diff S_{\vf{x}} = \vf{\tau}_{\beta}.
\end{aligned}
\right.
\end{equation}
The random Brownian motion of the particles involves computing the ``square root'' of the body mobility matrix, denoted by $\op{N}^{\frac{1}{2}}(\vf{Q})$  acting on a vector of independent white noise processes $\op{W}(t)$. {In order for fluctuation-dissipation balance to hold, 
the matrix $\op{N}^{\frac{1}{2}}$ must satisfy $\op{N}^{\frac{1}{2}} \left(\op{N}^{\frac{1}{2}} \right)^{\top} = \op{N}$, and does not necessarily have to be square.}
The last term on the right-hand-side of \Cref{BDeqns} is the stochastic drift term due to the Ito interpretation of the SDEs, where the divergence operator  $(\partial_{\vf{x}} \cdot )$ for a matrix-valued function $\vf{A}(\vf{x})$ is defined by $(\partial_{\vf{x}} \cdot \vf{A})_i = \sum_j \partial A_{ij}/\partial x_j$. 

Developing numerical schemes to integrate \Cref{BDeqns} has two main challenges. The first challenge is that, at every time step, one needs to generate the random velocity 
\begin{equation}
 \vf{U} = \left\{ \vf{u}_{\beta}, \vf{\omega}_{\beta}  \right\}_{\beta=1}^{N_{b}} 
 = \vf{\bar{U}} + \vf{\tilde{U}} 
 = \op{N} \vf{F} + \op{N}^{\frac{1}{2}} \vf{W}, \label{bodyvelocity}
\end{equation}
where $\vf{\bar{U}}$ is the deterministic particle velocity due to applied forces and torques, $\vf{\tilde{U}}$ is the random velocity due to the stochastic stress tensor, and $\vf{W}$ is a vector of independent and identically distributed (i.i.d.) Gaussian random variables with mean zero and covariance
\begin{equation}
\langle  W_i \, W_j \rangle =  \frac{2\kbt}{\Delta t} \delta_{ij},
\end{equation}
where $\Delta t$ is the time step size. More precisely, our task is to efficiently and accurately apply the action of $\op{N}$ and $\op{N}^{\frac{1}{2}}$. 
The second challenge is to compute or approximate the stochastic drift term $(\kbt)(\partial_{\vf{Q}} \cdot \op{N})$, which is conventionally handled by developing specialized stochastic temporal integrators \cite{BD_Fixman,MultiscaleIntegrators,BrownianBlobs,BrownianMultiblobSuspensions,FluctuatingFCM_DC}. The main focus of this paper is to tackle the first challenge by developing schemes  that generate the action of 
$\op{N}$ and $\op{N}^{\frac{1}{2}}$ based on a boundary integral formulation.



\subsection{Boundary value problem formulation}\label{sec:SSBVP}
For simplicity, we now consider only a single particle $\Gamma$ described by $\{ \vf{q}, \vf{\theta} \}$ immersed in a periodic domain $\mathcal{V}=[0,L]^2$, and we drop the subscript $\beta$. The generalization to account for many-body interaction is straightforward. In this section, we show that the random velocity given by \Cref{bodyvelocity}
can be obtained  by solving the {\it Stochastic Stokes Boundary Value Problem } (SSBVP):
\begin{equation}
\left \{
\begin{aligned}
-\div \vf{\sigma} &= \grad \pi - \eta \grad^2 \vf{v} = 0, \quad \vf{r} \in \mathcal{V} \backslash \overline{D},  \\ 
\div \vf{v} &= 0,   \\
\vf{v}(\vf{x}) &= \vf{u}+ \vf{\omega} \times (\vf{x}-\vf{q}) - \vf{\breve{v}}(\vf{x}) , \quad \vf{x} \in \Gamma, \\
\int_{\Gamma} \vf{\lambda}(\vf{x}) \diff S_{\vf{x}} &= \vf{f} \quad \text{and} \quad 
\int_{\Gamma} ( \vf{x} - \vf{q}) \times \vf{\lambda}(\vf{x})  \diff S_{\vf{x}} = \vf{\tau},
\end{aligned}
\right.
\label{SSBVP}
\end{equation}
where $\vslip(\vf{x})$ is a random surface velocity prescribed on the particle, that has zero mean and covariance
\begin{equation}
 \langle \vf{\breve{v}}(\vf{x}) \, \vf{\breve{v}}(\vf{y})  \rangle =   \frac{2\kbt}{\Delta t}\pStokeslet(\vf{x}- \vf{y}),  \quad \text{for all} ~ \left( \vf{x} \neq \vf{y} \right) \in \Gamma,
 \label{slipcovariance} 
\end{equation}
Here $\pStokeslet(\vf{r})$ is the Green's function for steady Stokes flow with viscosity $\eta$, and includes the specified boundary conditions (periodic BCs in our case). {For $\vf{x} = \vf{y}$, \Cref{slipcovariance} is not well-defined since $\pStokeslet$ is singular, which implies that $\vslip$ is a distribution and not a function; a more precise definition is given later in \Cref{KLexpansion}.}

By linearity of Stokes flow, the solution of \Cref{SSBVP} is the superposition of
\begin{align}
\vf{v} = \vf{\bar{v}} + \vf{\tilde{v}} ~,~ \vf{\sigma} = \vf{\bar{\sigma}} + \vf{\tilde{\sigma}} ~,~ 
\vf{U} = \vf{\bar{U}} + \vf{\tilde{U}}, 
\end{align}
where $\vf{\bar{U}} = \{ \vf{\bar{u}}, \vf{\bar{\omega}} \}$ and $ \vf{\tilde{U}} =  \{ \vf{\tilde{u}}, \vf{\tilde{\omega}} \}$ with $\{ \vf{\bar{v}}, \vf{\bar{\sigma}}, \vf{\bar{u}}, \vf{\bar{\omega}} \}$  satisfying a Stokes BVP {\it without} random surface velocity,
\begin{equation}
\left\{
\begin{aligned}
  -\div \vf{\bar{\sigma}} &= \grad \bar{\pi} - \eta \grad^2 \vf{\bar{v}} = \vf{0},  \\ 
  \div \vf{\bar{v}} &= 0, \\
  \vf{\bar{v}}(\vf{x}) &= \vf{\bar{u}}+ \vf{\bar{\omega}} \times (\vf{x}-\vf{q}), \quad \vf{x} \in \Gamma \\
  \int_{\Gamma}  \vf{\bar{\lambda}}(\vf{x}) \diff S_{\vf{x}} &=  \vf{f} \quad \text{and} \quad
\int_{\Gamma} (\vf{x}-\vf{q}) \times \vf{\bar{\lambda}}(\vf{x}) \diff S_{\vf{x}} = \vf{\tau},
\end{aligned}
\right.
\label{BVPnoslip}
\end{equation}
and $\{ \vf{\tilde{v}}, \vf{\tilde{\sigma}}, \vf{\tilde{u}}, \vf{\tilde{\omega}} \}$ satisfying a force- and torque-free Stokes BVP {\it with} a random surface velocity $\vslip$ with zero mean and covariance \eqref{slipcovariance},
\begin{equation}
\left\{
\begin{aligned}
  -\div \vf{\tilde{\sigma}} &= \grad \tilde{\pi} - \eta \grad^2 \vf{\tilde{v}} = \vf{0},  \\ 
  \div \vf{\tilde{v}} &= 0, \\
  \vf{\tilde{v}}(\vf{x}) &= \vf{\tilde{u}}+ \vf{\tilde{\omega}} \times (\vf{x} - \vf{q})- \vf{\breve{v}}(\vf{x}) , \quad \vf{x} \in \Gamma, \\
  \int_{\Gamma}  \vf{\tilde{\lambda}}(\vf{x}) \diff S_{\vf{x}} &=  \vf{0} \quad \text{and} \quad
\int_{\Gamma} (\vf{x}-\vf{q}) \times \vf{\tilde{\lambda}}(\vf{x})  \diff S_{\vf{x}} = \vf{0}.
\end{aligned}
\right.
\label{BVPslip}
\end{equation}
The BVP given by \Cref{BVPnoslip} is the standard mobility problem for rigid body motion that solves for the deterministic part of the particle velocity $\vf{\bar{U}} = 
\{ \vf{\bar{u}}, \vf{\bar{\omega}} \} = \op{N} \vf{F}$. The BVP given by \Cref{BVPslip}  generates its stochastic part $\vf{\tilde{U}}= \{\vf{\tilde{u}}, \vf{\tilde{\omega}} \} $. 

To show that the random particle velocity $\vf{\tilde{U}}$ determined by the mobility problem \eqref{BVPslip} indeed obeys the fluctuation-dissipation balance,
\begin{equation}
\langle \vf{\tilde{U}}  \vf{\tilde{U}}^{\top} \rangle =  \frac{2\kbt}{\Delta t} \op{N},
\end{equation}
we will invoke the \textit{Lorentz Reciprocal Theorem} (LRT) \cite[see Eq.~(1.4.5)]{BoundaryIntegral_Pozrikidis},
\begin{equation}
\int_{\Gamma} \vf{{u}} \cdot \vf{\lambda}' \diff S = \int_{\Gamma} \vf{u}' \cdot \vf{{\lambda}} \diff S, \label{LRT1}
\end{equation}
where $\{ \vf{u}, \vf{\lambda} \}$ and $\{ \vf{u}', \vf{\lambda}' \}$ are two arbitrary velocity-traction pairs corresponding to solutions of the homogeneous Stokes equations. 
To apply the LRT, we substitute  $\{ \vf{u}, \vf{\lambda} \}$ by $\{ \vf{\tilde{v}}, \vf{\tilde{\lambda}} \}$ from \Cref{BVPslip}, and $\{ \vf{u}', \vf{\lambda}' \}$ by $\{\vf{v}^{(i)}, \vf{\lambda}^{(i)}\}$,
where $\vf{v}^{(i)}$ and $\vf{\lambda}^{(i)}$ are the velocity and traction of the standard mobility problem \eqref{BVPnoslip} with applied force and torque $\vf{F} = \{ \vf{f}, \vf{\tau} \} = \vf{e}^{(i)}$, where $\vf{e}^{(i)} \in \mathbb{R}^3$ is the $i^{th}$ column of the identity matrix. By using the associated BCs in \Cref{BVPnoslip} and \Cref{BVPslip} to express $\vf{\tilde{v}}$ and $\vf{v}^{(i)}$ on $\Gamma$, and then making use of the force and torque balance conditions for $\vf{\tilde{\lambda}}$ and $\vf{\lambda}^{(i)}$, we can rewrite \Cref{LRT1} as
\begin{equation}
  \vf{\tilde{U}} \cdot \vf{e}^{(i)} = \tilde{U}_i = \int_{\Gamma} \vf{\breve{v}}(\vf{x}) \cdot \vf{\lambda}^{(i)} (\vf{x}) \diff S_{\vf{x}}.
\label{Udotei}
\end{equation}
This allows to express the covariance of $\vf{\tilde{U}}$ as:
\begin{equation}
\begin{aligned}
\langle  \tilde{U}_i  \tilde{U}_j \rangle &=
 \int_{\Gamma} \int_{\Gamma} \vf{\lambda}^{(j)} (\vf{x}) \cdot \langle \vf{\breve{v}}(\vf{x}) \vf{\breve{v}}(\vf{y})  \rangle \cdot \vf{\lambda}^{(i)} (\vf{y}) \diff S_{\vf{y}} \diff S_{\vf{x}} \\ 
 &=  
 \frac{2\kbt}{\Delta t} \int_{\Gamma} \int_{\Gamma} \vf{\lambda}^{(j)}(\vf{x}) \cdot \pStokeslet(\vf{x} - \vf{y}) \cdot \vf{\lambda}^{(i)}(\vf{y}) \diff S_{\vf{y}} \diff S_{\vf{x}}.
\end{aligned}
\label{deltaUcovariance}
\end{equation}
It can be shown (see \ref{appendix:LRT}) that the last integral in \Cref{deltaUcovariance} is equal to the ${(i,j)}^{th}$ element of the body mobility matrix $\op{N}$, giving the desired result,
\begin{equation}
\langle  \tilde{U}_i  \tilde{U}_j \rangle = \frac{2\kbt}{\Delta t} \mathcal{N}_{ij}.
\label{FDB_LRT}
\end{equation}
This shows that solving the SSBVP \eqref{SSBVP} gives the desired deterministic and stochastic particle velocity
\begin{equation}
\vf{U} = \vf{\bar{U}} + \vf{\tilde{U}} = \op{N} \vf{F} + \op{N}^{\frac{1}{2}} \vf{W}.
\label{DetPlusStochVelocities}
\end{equation}

While the LRT has been used in the past to analyze the nonhomogeneous Stokes BVP involving the fluctuating stress \cite{LLNS_FD_Fox}, here we employ it to establish that the homogeneous Stokes BVP given by  \Cref{BVPslip} and \Cref{slipcovariance} yields the correct statistics for the rigid body motion of the immersed particles.  The removal of the fluctuating stress driving the surrounding fluid allows for the eventual application of boundary integral techniques to solve the SSBVP.


\subsection{First-kind integral formulation}\label{sec:FirstKindContinuum}
For rigid particles moving in a Stokes fluid, we observe that the details of what happens inside the particle do not actually matter for its motion and its hydrodynamic interactions with other particles or boundaries. Therefore, it is possible to extend the fluid to the entire domain so that the fluid inside the body moves with a velocity that is continuous across the boundary of the body.  Once we extend the fluid to the interior of bodies, we may write down an alternative formulation of the SSBVP \eqref{SSBVP} as a {\it first-kind} boundary integral equation \cite{BoundaryIntegral_Pozrikidis},
\begin{equation}
 \vf{v}(\vf{x} \in \Gamma) = \vf{u} + \vf{\omega} \times (\vf{x}-\vf{q}) - \breve{\vf{v}}(\vf{x}) =  \int_{\Gamma} \pStokeslet(\vf{x}-\vf{y}) \, \vf{\psi}(\vf{y}) \diff S_{\vf{y}},  
\label{firstkindBIE}
\end{equation}
along with the force and torque balance conditions 
\begin{equation}
\int_{\Gamma} \vf{\psi}(\vf{x}) \diff S_{\vf{x}} = \vf{f} \quad \text{and} \quad 
\int_{\Gamma} ( \vf{x} - \vf{q}) \times \vf{\psi}(\vf{x})  \diff S_{\vf{x}} = \vf{\tau}.
\label{FTbalance3}
\end{equation}
Equations \eqref{firstkindBIE} and \eqref{FTbalance3}  together define a linear system of equations to be solved for the {\it single-layer density} $\vf{\psi}(\vf{x} \in \Gamma)$ and particle velocity $\vf{U} = \left\{ \vf{u}, \vf{\omega} \right\}$. We remark that $\vf{\psi}$ is the jump in the normal component of the stress when going across the body surface from the ``interior'' flow to the ``exterior'' flow. If $\vslip = 0$, then $\vf{\psi} = \vf{\lambda}$ is the traction. 

In operator notation, we can write the system formed by \Cref{firstkindBIE,FTbalance3} as a {\it saddle-point} problem,
\begin{equation}
 \left[
  \begin{array}{cc}
   \op{M} & -\op{K} \\
   -\op{K}^{*} & \vf{0}
  \end{array}
 \right]
 \left[
   \begin{array}{c}
    \vf{\psi} \\
    \vf{U}
   \end{array} 
 \right] = 
 -\left[
  \begin{array}{c}
    \vf{\breve{v}} \\
    \vf{F}
   \end{array} 
 \right],
\label{saddlesys}
\end{equation}
where $\op{M}$ denotes the {\it single-layer} integral operator defined by
\begin{equation}
 (\op{M} \vf{{\psi}})(\vf{x} \in \Gamma) = \int_{\Gamma} \pStokeslet(\vf{x} - \vf{y}) \, \vf{{\psi}}(\vf{y}) \diff S_{\vf{y}}, \label{singlelayerM}
\end{equation}
and $\op{K}$ is a geometric operator that relates particle velocity to surface velocities, 
\begin{equation}
(\op{K} \vf{U})(\vf{x} \in \Gamma) = \vf{u} + \vf{\omega} \times (\vf{x}-\vf{q}), \label{operatorK}
\end{equation}
and its adjoint $\op{K}^{\star}$ is an integral operator that converts the single-layer density $\vf{\psi}$ to a  force  and torque, 
\begin{equation}
\op{K}^{\star} \vf{\psi} =
\left(
\int_{\Gamma} \vf{\psi}(\vf{x}) \diff S_{\vf{x}} ~ , ~
\int_{\Gamma}  (\vf{x}-\vf{q}) \times \vf{\psi}(\vf{x}) \diff S_{\vf{x}}
\right) 
= (\vf{f}, \vf{\tau}). 
\label{Kstar}
\end{equation}

The covariance of the random surface velocity $\vslip$ can be written as 
\begin{equation}
\langle \vslip \vslip\rangle = \frac{2\kbt}{\Delta t} \op{M}, \label{vslipcovM}
\end{equation}
by which we mean that
$  \op{M} \vf{\psi}' = \langle (\vslip, \vf{\psi}') \, \vf{\vslip} \rangle $, for all $\vf{\psi}'$ in $L^2$-space, and $(\cdot, \cdot)$ denotes the $L^2$- inner product defined by
$(\vf{f}, \vf{g}) = \int_{\Gamma} \vf{f}(\vf{x}) \cdot \vf{g}(\vf{x}) \diff S_{\vf{x}}$. {If the random surface velocity $\vslip$ were a function and could therefore be evaluated pointwise,  \Cref{vslipcovM} would simply be a formal rewriting of \Cref{slipcovariance}. However, we reminder the reader again that \Cref{slipcovariance} is also formal (in the same way that \Cref{StochStressCov} is) and  $\vslip$ is a distribution and therefore {\em cannot} be evaluated pointwise. We can define $\vslip$ more precisely as follows.}
Since $\op{M}$ is a compact, self-adjoint and positive-semidefinite operator in the $L_2$ sense, it has countably infinitely many eigenvalues $\lambda_i \geq 0$ and orthonormal eigenfunctions $\vf{w}_i$, so we may write $\vslip$ in the  Karhunen-Lo\`{e}ve expansion
\begin{equation}
\vslip \overset{\mathrm{d.}}{=}  \sum_{i=1}^{\infty} \sqrt{\lambda_i} \, W_i \, \vf{w}_i,
\label{KLexpansion}
\end{equation}
where as before $W_i$ are independent Gaussian random variables with mean zero and variance $2\kbt / \Delta t$. 
In formal operator notation, we will write \Cref{KLexpansion}  as $\vf{\vslip} = \op{M}^{\frac{1}{2}}\vf{W}$, where $\left(\ophalf{M}\right) \left(\ophalf{M} \right)^{\top} = \op{M}$. {We have found \Cref{vslipcovM} (equivalently, \Cref{KLexpansion}), rather than the deceptively simple \Cref{slipcovariance}, to be a suitable starting point for a finite-dimensional discretization of $\vslip$, as we explain shortly.}

In this formulation, we require that $\vslip(\vf{x})$ is consistent with a divergence-free velocity field in the extended domain, \ie, 
\begin{equation} 
\int_{\Gamma} \vslip(\vf{x}) \cdot \vf{n}(\vf{x}) \diff S_{\vf{x}} = 0, \label{vslipnull}
\end{equation}
which is required  for \Cref{saddlesys} to be solvable since the single-layer operator  $\op{M}$ has a nontrivial null space consisting of single-layer densities that are normal to the boundary,
\begin{equation}
(\op{M} \vf{n})(\vf{x} \in \Gamma) =  \int_{\Gamma} \pStokeslet(\vf{x}-\vf{y})\, \vf{n}(\vf{y}) \diff S_{\vf{y}} = \vf{0}.
\end{equation}
From \Cref{KLexpansion}, we note that $\vslip = \op{M}^{\frac{1}{2}} \vf{W}$ is perpendicular to the null space of $\op{M}$, and hence, the solvability condition \Cref{vslipnull} is fulfilled. 
Formally\footnote{For the operator $\op{M}^{-1}$ to be well-defined, we need to resort to in its precise definition  either the space of band-limited functions or the fractional Sobolov space, but here, we will simply use this formal computation to inform our discretization and show later in \autoref{sec:FBIM-discrete} that our finite-dimensional discretization converges numerically.},
taking the Schur complement of \Cref{saddlesys} to eliminate $\vf{{\psi}}$ and solving for the body motion $\vf{U}$, we obtain
\begin{equation}
\begin{aligned}
 \vf{U} 
 = \vf{\bar{U}} + \vf{\tilde{U}} &= (\op{K}^{\star} \op{M}^{-1} \op{K})^{-1}\vf{F} + (\op{N} \op{K}^{\star} \op{M}^{-1})\vf{\breve{v}} \\
 &= (\op{K}^{\star} \op{M}^{-1} \op{K} )^{-1}\vf{F} + ( \op{N} \op{K}^{\star} \op{M}^{-1} \op{M}^{\frac{1}{2}} ) \vf{W} \\
 &= \op{N} \vf{F} + \op{N}^{\frac{1}{2}} \vf{W},
\end{aligned}
\end{equation}
which allows us to formally define $\op{N}$ and $\op{N}^{\frac{1}{2}}$ explicitly as,
 \begin{equation}
\begin{aligned}
 \op{N} &= (\op{K}^{\star} \op{M}^{-1} \op{K})^{-1}, \\
 \op{N}^{\frac{1}{2}} &= \op{N} \op{K}^{\star} \op{M}^{-1} \op{M}^{\frac{1}{2}}.
\end{aligned}
\label{contNandNhalf}
\end{equation}
Note that the random surface velocity $\vslip$ is in the range of $\op{M}$ by the construction of \Cref{KLexpansion}, so that $\op{M}^{-1}\vslip $ is well-defined.
We can therefore formally show that  fluctuation-dissipation balance holds in the  continuum operator sense, 
\begin{equation}
\begin{aligned}
 \op{N}^{\frac{1}{2}} \left( \op{N}^{\frac{1}{2}} \right)^{\top}
 &= \op{N} \op{K}^{\star} \op{M}^{-1} \ophalf{M} \left( \ophalf{M} \right)^{\top}\op{M}^{-1} \op{K}  \op{N} \\
 &= \op{N} ( \op{K}^{\star} \op{M}^{-1} \op{K} ) \op{N}  \\
 &=  \op{N} (\op{N})^{-1} \op{N}  =  \op{N}.
\end{aligned}
\label{Nhalfcov}
\end{equation}
This shows that the desired random velocity $\vf{U}$ in \Cref{bodyvelocity} can be generated by solving  \Cref{saddlesys}, which is the first-kind boundary integral formulation of the SSBVP \eqref{SSBVP}.
{While at first sight it may appear that we have simply formally rederived \Cref{FDB_LRT} by appealing to a first-kind BVP formulation, we will demonstrate next that the formal continuum formulation presented here has a well-defined finite-dimensional truncation that is very suitable for numerical computations.}

\section{Fluctuating boundary integral method}
\label{sec:FBIM-method}
This section presents the fluctuating boundary integral method (FBIM) for suspensions of Brownian rigid particles in two dimensions. {Our discrete formulation closely follows the (formal) continuum first-kind formulation presented in \autoref{sec:FirstKindContinuum}. We have found this to be much more effective than following the more general BVP formulation presented in \autoref{sec:SSBVP}. Specifically, one approach to discretizing the SSBVP \Cref{SSBVP} is to first generate a (smooth) random surface velocity {\em function} $\vslip$ using a {\em regularized} \footnote{The most direct way to regularize the singular Green's function is to represent it in Fourier space and then simply truncate the finite-dimensional sum to a finite number of Fourier modes.}  variant of \Cref{slipcovariance}, and then to use a standard spectrally-accurate second-kind boundary integral formulation to solve the resulting boundary-value problem. Our preliminary investigations of such an approach have revealed that regularization (truncation) of the covariance in \Cref{slipcovariance} leads to a drastic loss of accuracy in numerical fluctuation-dissipation balance. Instead, by relating the covariance of the random surface velocity to the first-kind operator as in \Cref{vslipcovM}, the regularization of the {\em distribution} $\vslip$ becomes directly connected to the singular quadrature used to discretize $\op{M}$. As we demonstrate here, this leads to a first-kind formulation that satisfies discrete fluctuation-dissipation (DFDB) to within solver tolerances, while preserving the underlying accuracy of the singular quadrature.}

We begin by presenting a discrete formulation of the mobility problem \eqref{saddlesys} by first discretizing the continuum operators $\op{M}$, $\op{K}$ and $\op{K}^{\star}$, to obtain a discrete saddle-point linear system, whose solution strictly obeys DFDB without any approximation. To efficiently solve the saddle-point linear system with Krylov iterative methods, we  present the two key components of FBIM: a fast routine for computing matrix-vector product of the  single-layer matrix $\mbf{M}$ (\ie, the discretized operator $\op{M}$), and a fast method for generating the random surface velocity $\bfslip = \sqrtM \mbf{W}$. To address the inherent ill-conditioning of the linear system arising from the first-kind integral formulation, we will also discuss block-diagonal preconditioning for the iterative solvers. 

\subsection{Discrete formulation of the mobility problem}
\label{sec:FBIM-discrete}
We present the discrete formulation of the mobility problem \eqref{saddlesys} by first discretizing the continuum operators $\op{K}$,  $\op{K}^{\star}$ and $\op{M}$. In the discrete formulation, for generality, we describe our method for many-body suspensions. Let us assume that $\Gamma_{\beta}$ is parametrized by $\vf{\gamma}_{\beta}: [0,2\pi] \rightarrow \mathbb{R}^2 $. We introduce a collection of $N_p$ equispaced points $s_j = j \Delta s$, $j \in \{ 1, \dots, N_p \}$, where $\Delta s = 2\pi/N_p$, and $\Gamma_{\beta}$ is discretized by the collection of nodes $\mbf{x}_{\beta} = \left\{ \mbf{x}_{\beta,j} \right\}_{j=1}^{N_p}$, where $\mbf{x}_{\beta,j} = \vf{\gamma}_{\beta}(s_j)$. We also denote the composite vector of all nodes or points by $\mbf{X} = \left\{ \mbf{x}_{\beta}\right\}_{\beta=1}^{N}$.

The discrete operator $\mbf{K}$ is a geometric matrix defined by
\begin{equation}
(\op{K}\vf{U})\left(\mbf{x}_{\beta,j}\right) = (\mbf{K}_{\beta} \vf{U}_{\beta})_j = \vf{u}_{\beta} +  \vf{\omega}_{\beta} \times(\mbf{x}_{\beta,j} - \vf{q}_{\beta}), 
\end{equation}
where $\mbf{K}_{\beta}$ is the sub-block of $\mbf{K}$ that maps the particle velocity $\vf{U}_{\beta} = \left\{ \vf{u}_{\beta}, \vf{\omega}_{\beta} \right\}$ to the velocity at the node $\mbf{x}_{\beta,j}$ on $\Gamma_{\beta}$.

The adjoint operator $\op{K}^{\star}$ defined by \Cref{Kstar} can be discretized by the periodic trapezoidal rule for each body, 
\begin{subequations}
\begin{align}
\int_{\Gamma_{\beta}} \vf{\psi} (\vf{x}) \diff S_{\vf{x}} &~\approx~ \sum_{j=1}^{N_p} \vf{\mu}_{\beta,j},   \\
\int_{\Gamma_{\beta}} (\vf{x} - \vf{q}_{\beta}) \times \vf{\psi} (\vf{x}) \diff S_{\vf{x}} &~\approx~\sum_{j=1}^{N_p} (\mbf{x}_{\beta,j} - \vf{q}_{\beta}) \times \vf{\mu}_{\beta,j}, 
\end{align}
\end{subequations}
where $\vf{\mu}_{\beta,j} = \vf{\psi}\left(\mbf{x}_{\beta,j} \right)\Delta s$ denotes an unknown discrete boundary force  at the node $\mbf{x}_{\beta,j}$ on $\Gamma_{\beta}$. We can therefore write
\begin{equation}
(\op{K}^{\star} \vf{\psi})(\mbf{X}) \approx \mbf{K}^{\top} \vf{\mu},
\end{equation}
where $\vf{\mu} = \left\{ \vf{\mu}_{\beta} \right\}_{\beta=1}^{N}$ is a composite vector that collects all the boundary forces. 

To approximate the single-layer operator $\op{M}$, we need to employ a singular quadrature that can handle the singularity of $\pStokeslet$. The development of  quadrature rules for singular or near-singular integrals  is an active research area of its own \cite{Alpert_Quadrature,SingularIntegral_Beale,SingularIntegral_Gimbutas,QBX_Epstein,QBX_Klockner,SingularIntegral_Kolm}. In two dimensions, for simplicity, we will use the Alpert quadrature which is based on a modification of  the trapezoidal rule with auxiliary nodes whose weights are configured to achieve the desired order of accuracy \cite{Alpert_Quadrature}. 

In matrix notation, the approximation of the single-layer integral operator evaluated on the vector of quadrature nodes $\mbf{X}$  can be compactly written as  
\begin{equation}
(\op{M} \vf{\psi})(\mbf{X})  \approx \mbf{M} \vf{\mu},
\label{SingleLayerMatrixDef}
\end{equation}
where $\mbf{M}$ is the discretized matrix operator of $\op{M}$, hereinafter referred to as the single-layer matrix.
The details of constructing $\mbf{M}$ will be discussed in \autoref{sec:fastMatVec}.
{We have chosen to keep track of surface forces in the discrete representation (rather than surface tractions as in the continuous case) as this yields both a symmetric saddle-point system in the usual sense, as well as the desirable property that $\mbf{M}_{ij} \rightarrow \pStokeslet(\mbf{x}_i - \mbf{x}_j)$ as $|\mbf{x}_i-\mbf{x}_j| \rightarrow \infty$. The matrix $\mbf{M}$ has a physical interpretation of a mobility matrix relating surface forces to surface velocities, which in turn suggests that $\mbf{M}$ should be symmetric (self-adjoint) in the standard $L_2$ sense.}


Following the continuum formulation (see \Cref{vslipcovM}), we require the discrete random surface velocity\footnote{The discrete $\bfslip$ can be thought of as being a suitably scaled finite-volume representation of the distribution $\vslip$.} to satisfy
\begin{equation}
\langle  \bfslip \bfslip ^{\top} \rangle = \frac{2\kbt}{\Delta t}\mbf{M}. \label{discreteCovariance} 
\end{equation}
More precisely, we need to generate a vector of Gaussian random variables $\bfslip$ whose covariance is given by \Cref{discreteCovariance}, and we denote it by $\bfslip = \sqrtM \mbf{W}$, where  $\mbf{W}$ is a finite-dimensional vector of Gaussian random variables with mean zero and covariance
\begin{equation}
\langle W_i W_j \rangle = \frac{2\kbt}{\Delta t} \delta_{ij}, 
\end{equation}
and $\sqrtM$ satisfies
\begin{equation}
\left( \sqrtM  \right) \left( \sqrtM \right)^{\top} = \mbf{M}.
\end{equation}
The discrete formulation of \Cref{saddlesys} is a saddle-point linear system for the boundary forces $\vf{\mu}$ and the rigid body motion $\vf{U}$, 
\begin{equation}
\left[
\begin{array}{cc}
\mbf{M} & -\mbf{K} \\
-\mbf{K}^{\top} & \mbf{0} 
\end{array}
\right]
\left[
\begin{array}{c}
\vf{\mu}  \\
\vf{U}
\end{array}
\right]
= 
- \left[
\begin{array}{c}
\sqrtM \mbf{W} \\
\vf{F}
\end{array}
\right]. 
\label{discreteSaddlePoint}
\end{equation}
By taking the Schur complement of $\mbf{M}$ and eliminating $\vf{\mu}$, we obtain
\begin{equation}
\begin{aligned}
\vf{U} &= \vf{N} \vf{F} + \vfhalf{N} \mbf{W} \\ 
&=  (\mbf{K}^{\top} \mbf{M}^{\dagger} \mbf{K} )^{-1}\vf{F} + ( \vf{N} \mbf{K}^{\top} \mbf{M}^{\dagger} \sqrtM ) \mbf{W}, 
\end{aligned}
\label{discreteParticleVelocity}
\end{equation}
from which we can define
\begin{equation}
\begin{aligned}
\vf{N} &= (\mbf{K}^{\top} \mbf{M}^{\dagger} \mbf{K} )^{-1}, \\
\vfhalf{N} &= \vf{N} \mbf{K}^{\top} \mbf{M}^{\dagger} \sqrtM,
\end{aligned}
\end{equation}
where $\mbf{M}^{\dagger}$ is the pseudo-inverse of $\mbf{M}$, and $\vf{N}$ is an approximation of  $\op{N}$ up to the order of accuracy of Alpert quadrature.  Under the definition of $\vf{N}$ and $\vfhalf{N}$, we can verify that our discrete formulation satisfies DFDB without any approximation,
\begin{equation}
\begin{aligned}
 \vf{N}^{\frac{1}{2}} \left( \vf{N}^{\frac{1}{2}} \right)^{\top}
 &= \vf{N} \mbf{K}^{\top} \mbf{M}^{\dagger} \mbfhalf{M}\left( \mbfhalf{M} \right)^{\top} \mbf{M}^{\dagger} \mbf{K}  \vf{N}  ,\\
 &= \vf{N} ( \mbf{K}^{\top} \mbf{M}^{\dagger} \mbf{K} ) \vf{N} , \\
 &=  \vf{N} (\vf{N})^{-1} \vf{N} =  \vf{N}.
\end{aligned}
\label{discreteFDB}
\end{equation}
Note that these relations are well-defined finite-dimensional versions of the formal continuum equations \eqref{contNandNhalf} and \eqref{Nhalfcov}. 

To generate the random velocity given by \Cref{discreteParticleVelocity}  efficiently, we solve the saddle-point linear system \eqref{discreteSaddlePoint} by GMRES. In the remaining sections,  we present efficient numerical that compute the matrix-vector product $\mbf{M} \vf{\mu}$ and generate the random surface velocity $\bfslip = \mbfhalf{M}\mbf{W}$. 

\subsection{Fast matrix-vector multiplication for the single-layer matrix} 
\label{sec:fastMatVec}


In this section we develop a fast method to efficiently perform the matrix-vector product $\mbf{M} \vf{\mu}$. The key idea for achieving linear-scaling is to use Ewald splitting  to decompose the periodic Stokeslet as 
\begin{equation}
\pStokeslet = \wStokeslet_{\xi} +  \rStokeslet_{\xi} = H * \pStokeslet + ( \pStokeslet - H*\pStokeslet) , 
\end{equation}
where ``$*$'' denotes convolution, and $\wStokeslet_{\xi}$ is a far-field (long-ranged) smooth kernel that decays exponentially in Fourier space, and $\rStokeslet_{\xi}$ is a near-field (short-ranged) singular kernel that decays exponentially in real space.
The choice of splitting function $H(r ; \xi)$ by Hasimoto \cite{Mobility2D_Hasimoto} is defined in Fourier space as
\begin{equation}
\widehat{H}(k ; \xi) = \left(1+\frac{k^2}{4\xi^2} \right) e^{-k^2/4\xi^2}, \label{HasimotoFunction}
\end{equation}
where $\xi$ is the splitting parameter that controls the rate of exponential decay.

Using the Fourier representation of the periodic Green's function of Stokes flow we can compute the far-field kernel in wave space as 
\begin{equation}
\wStokeslet_{\xi}(\vf{r}) = \frac{1}{\eta V}  \sum_{\mbf{k} \neq 0}  \frac{\widehat{H}(k ; \xi)}{k^2} (\mbf{I} - \hat{\mbf{k}} \hat{\mbf{k}}^{\top} ) \, e^{-i\mbf{k} \cdot \vf{r} }, \label{wave-space-kernel}
\end{equation}
where  $V = |\mathcal{V}| = L^2$, $\mbf{k} \in \{ 2\pi \kappa_i/L:  \kappa_i \in \mathbb{Z}, \, i=1,2\}$, and $\hat{\mbf{k}} = \mbf{k}/k$ for $k = |\mbf{k}|$. 
The near-field kernel is analytically computed in real space using the inverse Fourier transform\footnote{We gratefully thank Anna-Karin Tornberg for sharing with us notes on the Hasimoto splitting of the Stokeslet in two dimensions.},
\begin{equation}
\rStokeslet_{\xi}(\vf{r}) =  \frac{1}{4\pi \eta} \left[ \frac{1}{2} E_1(\xir) \mbf{I} + \left(\frac{\vf{r} \otimes \vf{r}}{r^2} - \mbf{I} \right) e^{-\xir} \right], \label{real-space-kernel}
\end{equation}
where $E_1(z)$ is the exponential integral defined by
\begin{equation}
E_1(z) = \int_{1}^{\infty} \frac{e^{-zt}}{t}  dt = \int_{z}^{\infty} \frac{e^{-t}}{t} dt. \label{expint}
\end{equation}
We observe from \Cref{real-space-kernel} that $\rStokeslet_{\xi}$ also has the logarithmic singularity of $\pStokeslet$, since in the limit $z \rightarrow 0$,
\begin{equation}
E_1 (z) = -\gamma - \log z + O(z), 
\end{equation}
where $\gamma$ is the Euler-Mascheroni constant. An important remark on the Hasimoto splitting is that it ensures both $\rStokeslet_{\xi}$ and $\wStokeslet_{\xi}$ are SPD, because $0\leq \widehat{H}(k ; \xi) \leq 1$ for all $\mbf{k}$ and $\xi$ \cite{SpectralRPY}. 

The splitting of $\pStokeslet$ naturally induces the splitting of the single-layer integral operator $\op{M} = \op{M}^{(r)} + \op{M}^{(w)}$, and subsequently, the corresponding splitting of the single-layer matrix 
\begin{equation}
\mbf{M} = \Mreal + \Mwave,
\end{equation}
where the elements of $\mbf{M}^{(w)}$ are obtained by applying the regular trapezoidal rule to $\op{M}^{(w)}$, which gives
\begin{equation}
\left(\mbf{M}^{(w)} \right)_{mn} = \wStokeslet_{\xi} (\mbf{x}_m - \mbf{x}_n). \label{Mwave-elements}
\end{equation}
The elements of $\mbf{M}^{(r)}$ are obtained by applying Alpert's hybrid Gauss-trapezoidal quadrature to  $\op{M}^{(r)}$, and can be futher decomposed as
\begin{equation}
\mbf{M}^{(r)} =  \mbf{M}^{(t)} + \mbf{M}^{(a)}.
\label{Mreal-decompose}
\end{equation}
where $\mbf{M}^{(t)}$ is the trapezoidal rule, \ie, $\left (\mbf{M}^{(t)} \right)_{mn} = \pStokeslet^{(r)}_{\xi}(\mbf{x}_m-\mbf{x}_n)$ for $m \neq n$, and we define it to be zero for $m=n$. The singular quadrature correction $\mbf{M}^{(a)}$ is a banded matrix that contains Alpert weights for the logarithmic singularity of $\pStokeslet_{\xi}^{(r)}$. The Alpert  weights are defined on a set of auxiliary nodes that do not coincide with the trapezoidal nodes, so a local Lagrangian interpolation from the auxiliary nodes to the trapezoidal nodes is needed to obtain the elements of $\mbf{M}^{(a)}$.
Note that $\Mreal$, $\Mwave$ and $\mbf{M}^{(a)}$ depend on $\xi$, but for conciseness of notation, the subscript $\xi$ is omitted.

A key ingredient of FBIM is a fast method to compute  the matrix-vector product
\begin{equation}
\mbf{M}\vf{\mu} = \left( \mbf{M}^{(a)} + \mbf{M}^{(t)} + \mbf{M}^{(w)} \right) \vf{\mu}.
\end{equation}
We recall that the Alpert quadrature assigns only local correction weights to the trapezoidal nodes. 
Since  $\mbf{M}^{(a)}$ is block-diagonal and banded, matrix-vector products involving $\mbf{M}^{(a)}$ can be computed in $\mathcal{O}(N)$ operations using vector rotations and sparse matrix-vector multiplications,
\begin{equation}
\left( \mbf{M}^{(a)} \vf{\mu} \right)_{\beta}= \mbf{R}_{\beta} \, \mbf{M}^{(a)}_{\text{ref}} \, \mbf{R}_{\beta}^{\top} \, \vf{\mu}_{\beta}, \quad \beta  = 1, \dots, N,
\end{equation}
where $\mbf{M}^{(a)}_{\text{ref}}$ is a precomputed Alpert matrix for some reference configuration, and $\mbf{R}_{\beta}$ is the rotation matrix from the chosen reference configuration to the configuration of $\Gamma_{\beta}$. 

To accelerate matrix-vector products involving $\Mreal$ and $\Mwave$ which are not sparse, we rely on the 
Spectral Ewald method \cite{SpectralEwald_Stokes}. 
For the near-field contribution, due to the short-ranged nature of $\rStokeslet_{\xi}$, the action of $\mbf{M}^{(t)}$ can be computed by using the cell list algorithm, commonly-used in Molecular Dynamics \cite{Allen_Tildesley_book}. First, we partition the computational domain into $N_\text{box} \times N_{\text{box}}$ cells and sort the points into these cells, which takes $\mathcal{O}(N)$ work. 
The splitting parameter $\xi$ is chosen such that the real-space sum converges to within a prescribed tolerance $\epsilon$, at a cutoff radius $\rcutoff = L/N_{\mathrm{box}}$. For each target point, the real-space sum is reduced to a local interaction with source points in its own cell and in all adjacent cells (with periodicity), \ie, nine cells in two dimensions. If the density of points in each cell is approximately held fixed as the system size grows, the complexity of the direct summation in the near field is also $\mathcal{O}(N)$.


In the SE method, the far-field contribution given by \Cref{wave-space-kernel,Mwave-elements} can be factorized as
\begin{equation}
\Mwave = \mbf{D}^{\star} \mbf{B} \mbf{D},
\label{Mwave-factorization}
\end{equation}
where the block-diagonal matrix $\mbf{B}$ is defined in the Fourier space as
\begin{equation}
\mbf{B}(\mbf{k}, \xi) = \frac{\widehat{H}(k; \xi)}{k^2}  (\mbf{I} - \hat{\mbf{k}}\hat{\mbf{k}}^{\top}),
\end{equation}
which essentially maps the Fourier representation of forces to velocities. 
The operator $\mbf{D}$ is the non-uniform Discrete Fourier Transform (NUDFT) that converts point forces $\vf{\mu} = \{ \vf{\mu}_n \}$ on a collection of non-uniform source points $\{\mbf{x}_n\}$ to Fourier space. 
The operator $\mbf{D}$ and its adjoint $\mbf{D}^{\star}$ can be efficiently applied using the non-uniform Fast Fourier Transform (NUFFT) (see \cite{NUFFT} and the references therein).  In operator notation, we can express the NUFFT as
\begin{equation}
\mbf{D} = \op{C} \op{F} \op{S} \quad \text{and} \quad \mbf{D}^{\star} = \op{S}^{\star} \op{F}^{\star} \op{C}^{\star},
\end{equation}
where $\op{S}$ and $\op{S}^{\star}$ are a pair of spreading and interpolation operators using Gaussian kernels, and $\op{F}$ and $\op{F}^{\star}$  are the forward and inverse FFT operators on a uniform grid, and $\op{C}$ and $\op{C}^{\star}$ are ``deconvolution'' operators.

The main idea of NUFFT is to first smear (spread) the point forces to a uniform grid using a Gaussian kernel, then make use of the FFT on the uniform grid, and finally apply the deconvolution (see \cite{NUFFT}).
In the SE method, the Gaussian kernel with the Fourier representation $e^{-\eta k^2/ 8\xi^2}$ is used for spreading and interpolation, where $\eta$ is a free parameter that sets the shape of Gaussian \cite{SpectralEwald_Stokes,SpectralEwald_Electrostatics}.
We observe that this Gaussian with parameter $\eta$ is different from the one in the Hasimoto function for Ewald splitting. The main purpose for introducing an extra parameter $\eta$ is that it allows the SE method to independently control the width of spreading and interpolation, so that the accuracy of evaluating the wave-space sum with NUFFT is decoupled from the accuracy of the Ewald sum. 

Next we address the choice of parameters in the SE method. The parameters that need to be set by the user include the number of partition cells $N_{\mathrm{box}}$ (or $\rcutoff$), the splitting parameter $\xi$, the size of the Fourier grid $M$ (even), the number of points $P$ (odd) for spreading and interpolation in the NUFFT, and the Gaussian shape parameter $m$ (which is related to $\eta$, see \cite[Eq.~(22)]{SpectralEwald_Stokes}). First, we set the number of partition cells $N_{\mathrm{box}}$, primarily based on balancing the computational work between the real- and wave-space sums (see \autoref{sec:SuspensionOfDisks} for details), and set $\rcutoff = L/N_{\mathrm{box}}$. For a user-specified error tolerance level $\epsilon$, the splitting parameter $\xi$ and $M$ are determined by the truncation error estimates in the real and wave spaces \cite{SpectralEwald_Stokes}, respectively,
\begin{subequations}
\begin{align}
 C_r e^{-\xi^2 \rcutoff^2} & \leq \epsilon, \\
 C_w e^{-k_{\max}^2 / 4\xi^2}  & \leq \epsilon,
\end{align}
\end{subequations}
where $k_{\max}$ corresponds to the largest mode of a grid of size $M\times M$. The constants in the error estimates are estimated empirically as in \cite{SpectralEwald_Stokes}, and we have found that $C_r \approx 100$ and $C_w \approx 1$ are suitable for the Stokeslet in two dimensions. 
The remaining parameters $P$ and $m$ (or $\eta$) together determine the accuracy of evaluating the wave-space sum using NUFFT. From the error estimation of Lindbo and Tornberg \cite{SpectralEwald_Stokes,SpectralEwald_Electrostatics}, 
the optimal choice is $m \sim \sqrt{\pi P}$, so that the approximation error arising from the NUFFT (including quadrature error and Gaussian truncation error) is approximately $e^{-\pi P /2} \lesssim \epsilon$, from which $P$ is determined.

The overall complexity of the wave-space sum is $\mathcal{O}(N)$ for spreading and interpolation and $\mathcal{O}(M^2 \log M^2)$ for the FFTs. For dense suspensions in two dimensions, we have found in our implementation that the computational work is dominated by spreading and interpolation instead of the FFTs.

\subsection{Fast sampling of the random surface velocity}
\label{sec:sqrtM}
The second key component of FBIM is a fast routine for sampling the random surface velocity $\bfslip =\mbfhalf{M} \mbf{W}$. The action of $\mbf{M}^{1/2}$ can be computed using iterative methods such as the Chebyshev polynomial approximation  \cite{BD_Fixman_sqrtM} or the Lanczos algorithm for matrix square root \cite{SquareRootKrylov}. For the Stokeslet in two dimensions, applying iterative methods directly for $\mbf{M}$ are expected to scale poorly for dense suspensions because of its logarithmic growth in the far field. The condition number of $\mbf{M}$, as well as the number of iterations required to converge for a given tolerance level, would grow with the system size. This would make the overall complexity for applying the action of $\mbfhalf{M}$ superlinear, even though the action of $\mbf{M}$ can be applied with $\mathcal{O}(N)$ work.

To improve the superlinear complexity for generating the random surface velocity with iterative methods, {we apply the {\it Positively Split Ewald} (PSE) method (developed for the RPY tensor) to the two-dimensional periodic Stokeslet \cite{SpectralRPY}.} The main idea of PSE is to split the action of $\mbf{M}^{1/2}$ into a near-field part $\left(\mbf{M}^{(r)} \right)^{1/2}$  and a far-field part $\sqrtMwave$, and generate the random surface velocity as
\begin{equation}
\mbf{M}^{1/2} \mbf{W} \overset{\mathrm{d.}}{=} \left(\mbf{M}^{(r)} \right)^{1/2} \mbf{W}^{(r)} + \sqrtMwave \mbf{W}^{(w)}, 
\label{sqrtM-split}
\end{equation}
so that the near-field contribution can be rapidly generated by the Lanczos algorithm \cite{SquareRootKrylov} in the real space, and the far-field contribution can be efficiently handled in the wave space by NUFFT.
The right-hand-side of \Cref{sqrtM-split} defines one way of computing the action of $\mbfhalf{M}$ provided that $\mbf{W}^{(r)}$ and $\mbf{W}^{(w)}$ are two independent Gaussian random vectors.

For a sufficiently small cut-off radius or a sufficiently large $\xi$, the real-space kernel $\rStokeslet_{\xi}$ decays exponentially, so that the near-field interaction is localized and the condition number of $\mbf{M}^{(r)}$ does not grow with the number of bodies, while the packing fraction is held fixed. However, due to the singular nature of $\rStokeslet_{\xi}$, the condition number of $\mbf{M}^{(r)}$ may grow if the number of points per body increases. 
We have found that the Lanczos algorithm \cite{SquareRootPreconditioning} with block-diagonal preconditioning (see \autoref{sec:block-diagonal-precond})  can significantly reduce the number of iterations. In the case when two bodies nearly touch, so that the problem itself becomes ill-conditioned, the block-diagonal approximation becomes worse and the number of iterations for all iterative solvers increases. Nevertheless, we have found that the block-diagonal preconditioner is still effective, and it takes a reasonable number of iterations for the Lanczos algorithm to converge, even for dense suspensions (see \autoref{fig:IterativeSolverConv}).

Note that the validity of \Cref{sqrtM-split} relies on the property that both $\mbf{M}^{(r)}$ and $\mbf{M}^{(w)}$ need to be SPD.
We observe that, even though $\rStokeslet_{\xi}$ is a SPD kernel, its singular contribution given by the Alpert correction matrix $\mbf{M}^{(a)}$ is not symmetric for a general-shaped particle. We have confirmed numerically that using only the symmetric part of $\mbf{M}^{(a)}$ does not affect the accuracy of the Alpert quadrature. In practice, we observe that the Lanczos algorithm is rather insensitive to spurious negative eigenvalues of small magnitude that may exist for $(\mbf{M}^{(r)}+(\mbf{M}^{(r)})^{\top})/2$.

The far-field matrix $\mbf{M}^{(w)}$ is SPD by construction because of  \Cref{wave-space-kernel}, and we can rewrite \Cref{Mwave-factorization} as 
\begin{equation}
\Mwave = \left( \mbf{D}^{\star} \mbf{B}^{1/2} \right) \left( \mbf{D}^{\star} \mbf{B}^{1/2} \right)^{\star},
\end{equation}
with $\mbf{B}^{1/2}$ defined in the wave space as
\begin{equation}
\mbf{B}^{1/2}(\mbf{k},\xi) = \frac{1}{k} \widehat{H}^{1/2}(k;\xi)(\mbf{I} - \hat{\mbf{k}}\hat{\mbf{k}}^{\top}),
\end{equation}
The far-field contribution of the random surface velocity can  be generated by
\begin{equation}
\sqrtMwave \mbf{W}^{(w)} = \mbf{D}^{\star} \mbfhalf{B}\mbf{W}^{(w)} = \op{S}^{\star} \op{F}^{\star} \op{C}^{\star} \mbfhalf{B} \mbf{W}^{(w)},
\label{sqrtMwave}
\end{equation}
where $\mbf{W}^{(w)}(\vf{\kappa})$ is a complex-valued random vector in the wave space, and $\kappa_i \in \left\{ -\frac{M}{2}, \dots, \frac{M}{2}-1 \right\}$.
The sequence of operations in \Cref{sqrtMwave} can be interpreted as follows. We first generate random numbers in the wave space, and project onto the divergence-free subspace by the projection $\mbf{I} - \hat{\mbf{k}}\hat{\mbf{k}}^{\top}$, then scale by $\widehat{H}^{1/2}(k;\xi)/k$, and apply deconvolution and (inverse) FFTs to obtain  velocities in the real space on the grid, and finally, perform interpolation using a Gaussian kernel to obtain the random surface velocities on the particles. This is equivalent to how random velocities are generated in methods like FIB \cite{BrownianBlobs} and fluctuating FCM \cite{ForceCoupling_Fluctuations,FluctuatingFCM_DC}. We remark that the cost of applying the action of $\sqrtMwave$ in \Cref{sqrtMwave} is even cheaper than the action of $\Mwave$, since it only requires half the work. 
 
We note that certain complex-conjugate symmetry of $\mbf{W}^{(w)}(\vf{\kappa})$ must be maintained to ensure its Fourier transform gives a real-valued Gaussian random vector with the correct covariance in the real space. Specifically, we require that, the zeroth mode $\vf{\kappa}=(0,0)$ is set to be zero, and the Nyquist modes $\vf{\kappa} \in \left\{ (-\frac{M}{2},0), (0,-\frac{M}{2}), (-\frac{M}{2},-\frac{M}{2}) \right\}$  are real-valued and generated from  $\mathscr{N}(\mbf{0}, \mbf{I}_{2\times 2})$. All the remaining modes are generated by $\mbf{W}^{(w)}(\vf{\kappa}) = \mbf{a} + i\mbf{b}$ for $\mbf{a}, \mbf{b} \in \mathscr{N}(\mbf{0}, \frac{1}{2}\mbf{I}_{2\times 2})$ , and have the complex-conjugate symmetry: $\mbf{W}^{(w)}(\vf{\kappa}) = ( \mbf{W}^{(w)}(\vf{\kappa}'))^*$, where $\vf{\kappa}' = -\vf{\kappa} \mod M$.

\subsection{Block-diagonal preconditioning}
\label{sec:block-diagonal-precond}
This section presents the block-diagonal preconditioning technique introduced in \cite{RigidMultiblobs} for solving the ill-conditioned  linear system \eqref{discreteSaddlePoint} with GMRES, and for generating the near-field contribution of the random surface velocity with the Lanczos algorithm \cite{SquareRootKrylov,SquareRootPreconditioning}. The block-diagonal preconditioner for the linear system \eqref{discreteSaddlePoint} is obtained by neglecting all hydrodynamic interactions between different bodies, \ie,
\begin{equation}
\mbf{P} = \left[
\begin{array}{cc}
\widetilde{\mbf{M}} &  -\mbf{K} \\
-\mbf{K}^{\top} & \mbf{0}
\end{array}
\right],
\label{preconditionerGMRES}
\end{equation}
where $\widetilde{\mbf{M}}$ is a block-diagonal approximation of $\mbf{M}$ obtained by setting elements of $\mbf{M}$ corresponding to pairs of points on {\it distinct} bodies to zero,
\begin{equation}
\widetilde{\mbf{M}}_{\alpha \beta} = \delta_{\alpha \beta} {\mbf{M}}_{\alpha \beta},
\end{equation}
and the subscripts with Greek letters denote the sub-block containing interactions between $\Gamma_{\alpha}$ and $\Gamma_{\beta}$.
Applying the preconditioner to \Cref{discreteSaddlePoint} amounts to solving a linear system,
\begin{equation}
\left[
\begin{array}{cc}
\widetilde{\mbf{M}} &  -\mbf{K} \\
-\mbf{K}^{\top} & \mbf{0}
\end{array}
\right]
\left[
\begin{array}{c}
\vf{\mu} \\
\vf{U}
\end{array}
\right]
=
-\left[
\begin{array}{c}
\bfslip \\
\vf{F}
\end{array}
\right],
\label{preconditionsys}
\end{equation}
which requires the action of the approximate body mobility matrix (Schur complement),
\begin{equation}
\widetilde{\vf{N}} = \left( \mbf{K}^{\top} \widetilde{\mbf{M}}^{\dagger} \mbf{K} \right)^{-1}.
\end{equation} 
The approximate body mobility matrix $\widetilde{\vf{N}}$ is also block-diagonal, and can be efficiently applied for each body,
\begin{equation}
\widetilde{\vf{N}}_{\beta\beta} = \left( \mbf{K}_{\beta}^{\top}  \mbf{M}_{\beta\beta}^{\dagger} \mbf{K}_{\beta}  \right)^{-1}.
\label{approximate-body-mobility}
\end{equation}
The matrix block $\mbf{M}_{\beta \beta}$ in \Cref{approximate-body-mobility} is a small matrix with size $2N_p \times 2N_p$, which is precomputed for the reference configuration,
\begin{equation}
\mbf{M}_{\text{ref}} = \mbf{M}^{(a)}_{\text{ref}} + \mbf{M}^{(t)}_{\text{ref}} + \mbf{M}^{(w)}_{\text{ref}}.
\end{equation}
The action of $\mbf{M}_{\beta \beta}^{\dagger}$ can be efficiently applied by using $\mbf{M}_{\text{ref}}^{\dagger}$, which is also precomputed using a dense SVD decomposition or eigenvalue decomposition, and by using rotation matrices for different bodies because of the translational and rotational invariance of the free-space Stokeslet,
\begin{equation}
\mbf{M}_{\beta\beta}^{\dagger} \approx \mbf{R}_{\beta} \mbf{M}_{\text{ref}}^{\dagger} \mbf{R}_{\beta}^{\top}.
\label{M-block-inverse}
\end{equation}
We note that the two sides of \Cref{M-block-inverse} do not equal exactly, since the Alpert quadrature is not rotation-invariant for a general-shaped body. This is not an issue, since the error introduced by this artifact is within the error tolerance, and the preconditioner does not need to be exactly inverted to work effectively. 

For the near-field contribution of the random surface velocity, we use the preconditioned Lanczos algorithm  \cite{SquareRootPreconditioning} to generate 
\begin{equation}
\bfslip^{(r)} = \mbf{G}^{\dagger}\left(\mbf{G}\mbf{M}^{(r)} \mbf{G}^{\top} \right)^{1/2} \mbf{W}^{(r)},
\label{preconditioned-real-part}
\end{equation}
where $\mbf{G}$ is a block-diagonal preconditioner, whose diagonal blocks can be precomputed as a dense matrix for the reference configuration using the eigenvalue decomposition,
\begin{equation}
 \mbf{M}^{(r)}_{\text{ref}} = \mbf{V} \vf{\Sigma} \mbf{V}^{*},
\label{Gprecond-decomp}
\end{equation}
and we set
\begin{equation}
\mbf{G}_{\text{ref}} = ( \vf{\Sigma}^{\dagger})^{1/2} \mbf{V}^* \quad \text{and} \quad 
\mbf{G}_{\text{ref}}^{\dagger} =  \mbf{V} \vf{\Sigma}^{1/2},
\label{Gprecond-precompute}
\end{equation}
where the diagonal elements of $\vf{\Sigma}$ and $\vf{\Sigma}^{\dagger}$ corresponding to the spurious eigenvalues of $\vf{\Sigma}$ are set to be zero.
This gives the desired discrete fluctuation-dissipation balance, 
\begin{equation}
\left\langle \bfslip^{(r)}_{\text{ref}} ~ \bfslip^{(r)}_{\text{ref}}  \right\rangle = 
\mbf{G}_{\text{ref}}^{\dagger} \mbf{G}_{\text{ref}}\mbf{M}^{(r)} \mbf{G}_{\text{ref}}^{\top}  \left(\mbf{G}_{\text{ref}}^{\dagger}\right)^{\top} = \mbf{M}^{(r)}_{\text{ref}}.
\end{equation}
Generating $\bfslip^{(r)}$ as in \Cref{preconditioned-real-part} is efficient, since we can reuse the precomputed matrices in \Cref{Gprecond-precompute} throughout the simulation, and apply rotation matrices for each body,
\begin{equation}
\mbf{G}_{\beta \beta} = \mbf{R}_{\beta} \mbf{G}_{\text{ref}}  \mbf{R}_{\beta}^{\top}.
\end{equation}
We remark that the use of block-diagonal preconditioners does not increase the overall complexity of FBIM since they can be applied with $\mathcal{O}(N)$ work, and, furthermore, the cost is amortized over the length of the BD simulations.

%
%

%
%

\section{Numerical Results}
\label{sec:FBIM-results}
This section presents numerical results of applying the FBIM to a number of benchmark problems in two dimensions. We first address the effect of the Ewald splitting parameter $\xi$ on the accuracy of the first-kind mobility solver. Next we test the first-kind mobility solver by applying it to the steady Stokes flow through a square periodic array of disks, and compare the results to well-known analytical solutions.  In the third test, we consider suspensions of Brownian rigid disks, and assess  the effectiveness of FBIM by its accuracy and convergence, the robustness of iterative solvers, and its scalability to simulate suspensions of many-body particle systems. In the last set of numerical examples, we perform {\it Brownian Dynamics} (BD) simulations by combining the FBIM with previously-developed stochastic temporal integrators \cite{MultiscaleIntegrators,BrownianBlobs}. By simulating the free diffusion of a non-spherical (starfish-shaped) particle, we confirm that correctly handling the stochastic drift term in the temporal integrator is necessary in order to reproduce the equilibrium Gibbs-Boltzmann distribution. As a simple but nontrivial benchmark problem for many-body suspension, we investigate the dynamics of a pair of starfish-shaped particles connected by a harmonic spring, which includes interaction through the spring potential, hydrodynamic interaction between the particles and with their periodic images, as well as Brownian noises. 
 
\subsection{Choosing the Ewald splitting parameter}
\label{sec:PeriodicArrayDisks}
The first important aspect of FBIM that needs to be addressed is how the choice of $\xi$ affects the accuracy of the first-kind mobility solver. In the {\it Spectral Ewald} (SE) method \cite{SpectralEwald_Stokes} the choice of $\xi$ is only based on balancing the computational work between the real- and Fourier-space sums. In practice, we want to choose a sufficiently large $\xi$ (\ie, a short-ranged singular kernel $\rStokeslet_{\xi}$), so that the computational work in the real-space can be made cheap at the expense of the FFT in Fourier space. In this work, we demonstrate that the choice of $\xi$ is also limited by the accuracy of Alpert quadrature used to resolve the logarithmic singularity of $\pStokeslet$ (diagonal elements of $\mbf{M}$). 


\begin{figure}
 \centering
 \includegraphics[width=.5\linewidth]{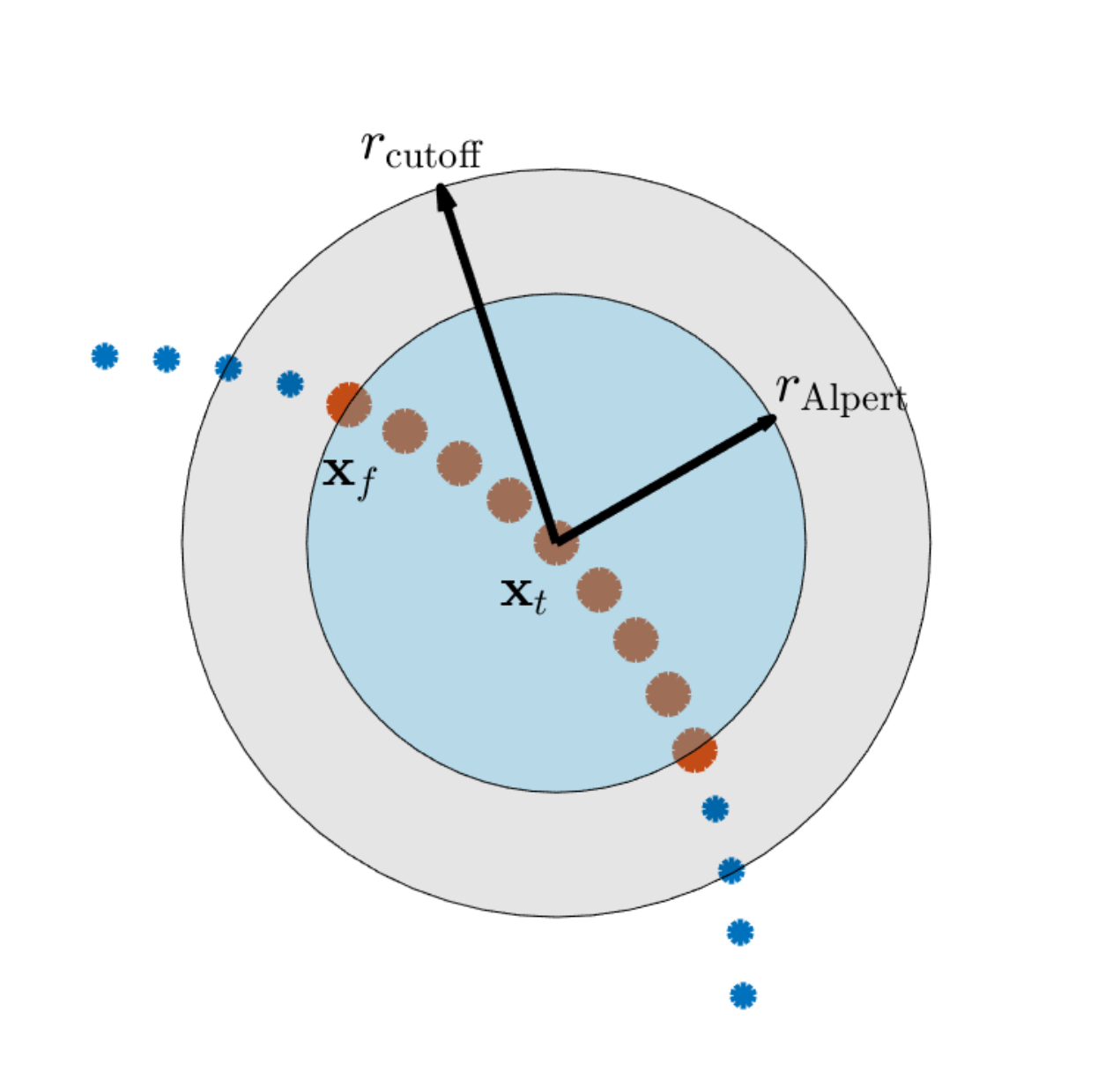} \hspace{-1em}
 \includegraphics[width=.5\linewidth]{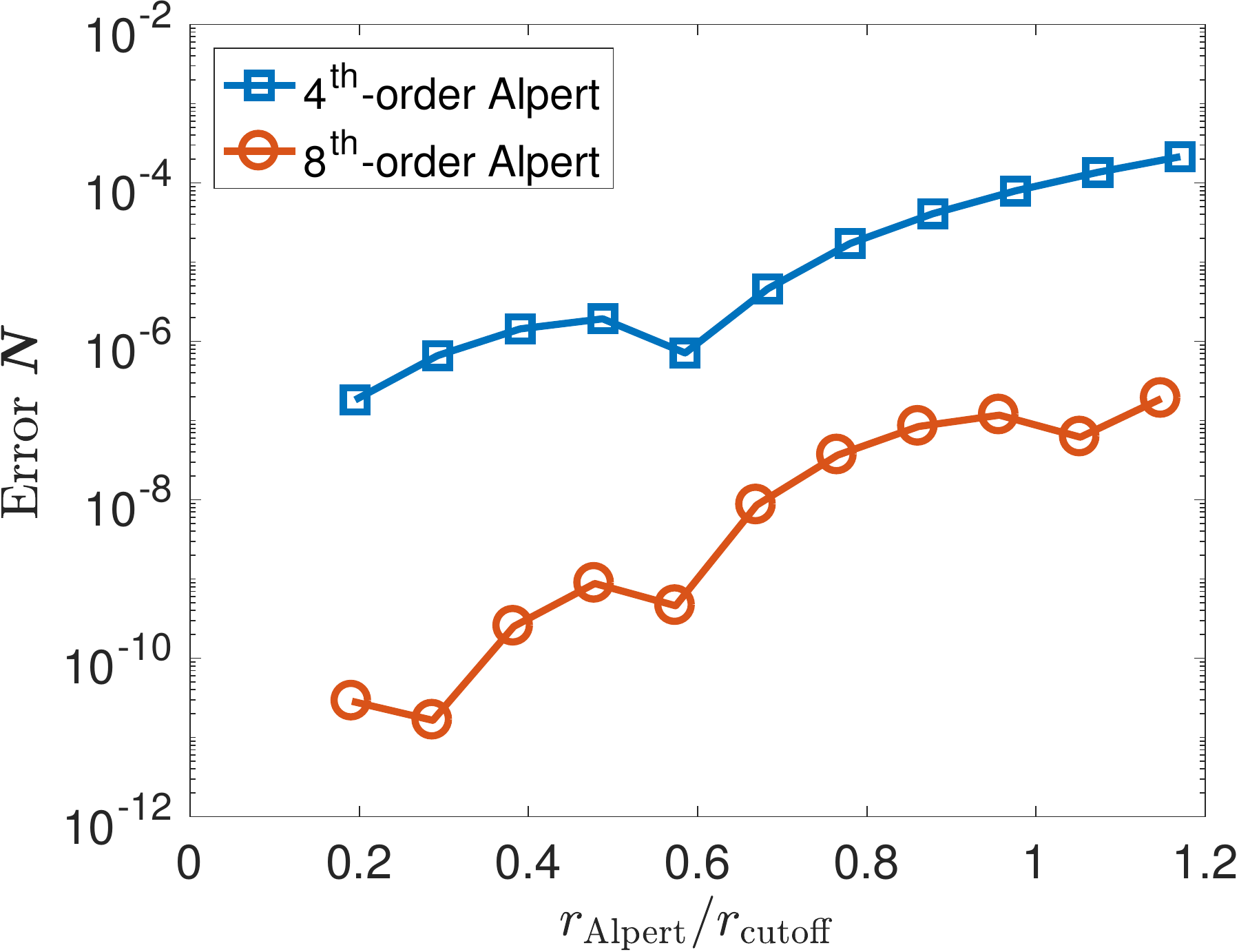} 
 \caption{({\it Left panel}) A cartoon illustration of radius of Alpert correction $\rAlpert$ and the cut-off radius of Ewald summation $\rcutoff$. Both blue and red nodes represent the quadrature nodes in the trapzoidal rule, and the red nodes are the special nodes with Alpert correction weights centered a target node $\mbf{x}_t$. The radius of Alpert correction is defined by $ \rAlpert = | \mbf{x}_t - \mbf{x}_f |$, where $\mbf{x}_f$ is the last node with Alpert correction weight for a chosen order of Alpert quadrature. ({\it Right panel}) Normalized error (in matrix 2-norm) of $\vf{N}$ with respect to $\op{N}$ (approximated with 12 digits of accuracy) versus the ratio $\rAlpert / \rcutoff$ for the $4^{th}$- and $8^{th}$-order Alpert quadratures. The mobility solver gradually loses accuracy because the singular kernel (with support $\rcutoff$) is not sufficiently resolved by the Alpert quadrature as the ratio increases.}
 \label{fig:AlpertXi}
\end{figure}

The Alpert quadrature can be viewed as assigning (interpolated) local correction weights to nearby quadrature nodes of a target point $\mbf{x}_{t}$, as illustrated in the left panel of \autoref{fig:AlpertXi}. For example, the number of (one-sided) quadrature nodes with nonzero correction weights is 4 and 8 for the $4^{th}$- and $8^{th}$-order Alpert quadrature, respectively. 
Since $\rStokeslet_{\xi}$ decays exponentially with length scale $\xi^{-1}$, the Alpert quadrature grid must resolve length scales smaller than $\xi^{-1}$ in order to capture the logarithmic singularity of $\rStokeslet_{\xi}$ at the origin. Thus for a fixed grid, as $\xi$ increases, we expect the accuracy of the Alpert quadrature to become progressively worse.

To demonstrate this issue, we study how the error of $\vf{N}$ changes with $\xi$ as follows. 
We fix the packing fraction $\phi = \pi/16 \approx 0.196$ and $N_p = 64$ but use different values of $\xi$ (\ie, different $\rcutoff$). The radius of Alpert correction $\rAlpert$ is defined by the distance between the target point $\mbf{x}_{t}$ and the farthest quadrature node $\mbf{x}_f$ with nonzero Alpert weight, \ie, $\rAlpert = |\mbf{x}_t - \mbf{x}_f|$. Note that $\rAlpert$ is fixed once the quadrature order is fixed. The Ewald sum is computed with accuracy $\tolEwald= 10^{-9}$ and the tolerance level of GMRES is $\tolIter=10^{-9}$. The body mobility matrix $\vf{N}$ is computed by solving the deterministic mobility problem with force and torque $\vf{F}$ set to be the columns of the $3\times 3$ identity matrix, which gives the corresponding columns of $\vf{N}$.
We can also compute the exact body mobility matrix $\op{N}$ to 12 digits of accuracy with $N_p = 256$ by using the second-kind boundary integral formulation \delete{(attached as supplementary material of this paper)}(see \ref{appendix:secondkind-formulation}). In two dimensions the resulting linear system of the second-kind formulation is well-conditioned and its solution is spectrally accurate. In the right panel of \autoref{fig:AlpertXi}, we show the normalized error of $\vf{N}$ (in matrix 2-norm) with respect to the 12-digit accurate approximation of $\op{N}$ for the $4^{th}$- and $8^{th}$-order Alpert quadrature, and the error increases  by at least two orders of magnitude for the range of $\rAlpert/ \rcutoff$ considered. We conclude that choosing $\xi$ with $\rAlpert / \rcutoff \lesssim 0.6$ maintains the accuracy of Alpert quadrature sufficiently well.

\subsection{Square lattice of disks}

Steady Stokes flow around a square periodic array of fixed disks in two dimensions is one of the classical problems in fluid mechanics, and its analytic solution is a thoroughly studied topic in the literature. Notably, Hasimoto \cite{Mobility2D_Hasimoto} obtained an analytical expression for the drag force $F$ on a dilute array of disks moving with velocity $U$ by solving the steady Stokes equations with Fourier series expansions. Later, Sangani and Acrivos \cite{PeriodicArray2D_Acrivos1982} extended Hasimoto's solution to the semi-dilute regime by including higher-order correction terms, and obtained the expansion
\begin{equation}
\frac{F}{\eta U} = \frac{4\pi}{ -\ln{\sqrt{\phi}} - 0.738 + \phi -0.887\phi^2 + 2.039 \phi^3 + O(\phi^4)},
\label{DiluteTheory}
\end{equation} 
where $a$ is the radius of disk, $l$ is the spacing of the square lattice of disks, and $\phi = \pi a^2/ l^2$ is the packing fraction. Another important regime for which there are theoretical solutions is the densely-packed regime. In this regime, the disks almost touch, so that there is flow through a narrow gap between two neighboring disks, and a lubrication correction is generally required to extend the solution to this regime. Sangani and Acrivos \cite{PeriodicArray2D_Acrivos1982b} obtained using lubrication theory the asymptotic formula
\begin{equation}
\frac{F}{\eta U} \approx \frac{9\pi}{2\sqrt{2}} \varepsilon^{-\frac{5}{2}}
\label{DenseTheory}
\end{equation}
where $\varepsilon = (l-2a)/l = 1-\sqrt{4\phi/\pi}$ is the relative gap between two neighboring disks. 

In the following test, we numerically solve the Stokes mobility problem \eqref{saddlesys} for the rigid body motion $\vf{U}$ with applied force and torque $\vf{F} = (1,0,0)$ to determine the relationship between the drag force and velocity. We set $a=1.0$ and $\tolEwald = \tolIter = 10^{-9}$, while varying the length of square lattice $l$ to  achieve different packing fractions. For this simple test problem, we do not seek to optimize code performance, and we set $\Nbox = 4^2$ ($\xi \approx 20.13$) for all the test cases. In \autoref{sec:SuspensionOfDisks} we will address the choice of optimal $\xi$ in the example of dense suspension of rigid particles. For a dilute or semi-dilute suspension ($\phi < 0.2$), the  boundary $\Gamma$ is discretized with $N_p = 64$ quadrature nodes, which is a sufficient resolution for the first-kind solver to give at least 6 digits of accuracy. For a higher packing fraction ($\phi > 0.2$), the number of quadrature nodes  is determined by the ratio of the smallest gap between any two neighboring disks $d_g = l-2a$ to the spacing between quadrature nodes $d_s = 2\pi a /N_p$.
In a moderately-resolved computation, we require that the quadrature node spacing  is comparable or smaller than the gap spacing, \ie, $d_s/d_g \lesssim 1$.
For example, for the highest packing fraction considered in this test $\phi = 0.76$ (note that for a close-packed square lattice $\phi_{\max} \approx 0.7854$), the number of quadrature nodes set by the ratio $d_s/d_g =1$ is $N_p \approx 190$.

\begin{figure}
 \centering
 \includegraphics[width=.48\linewidth]{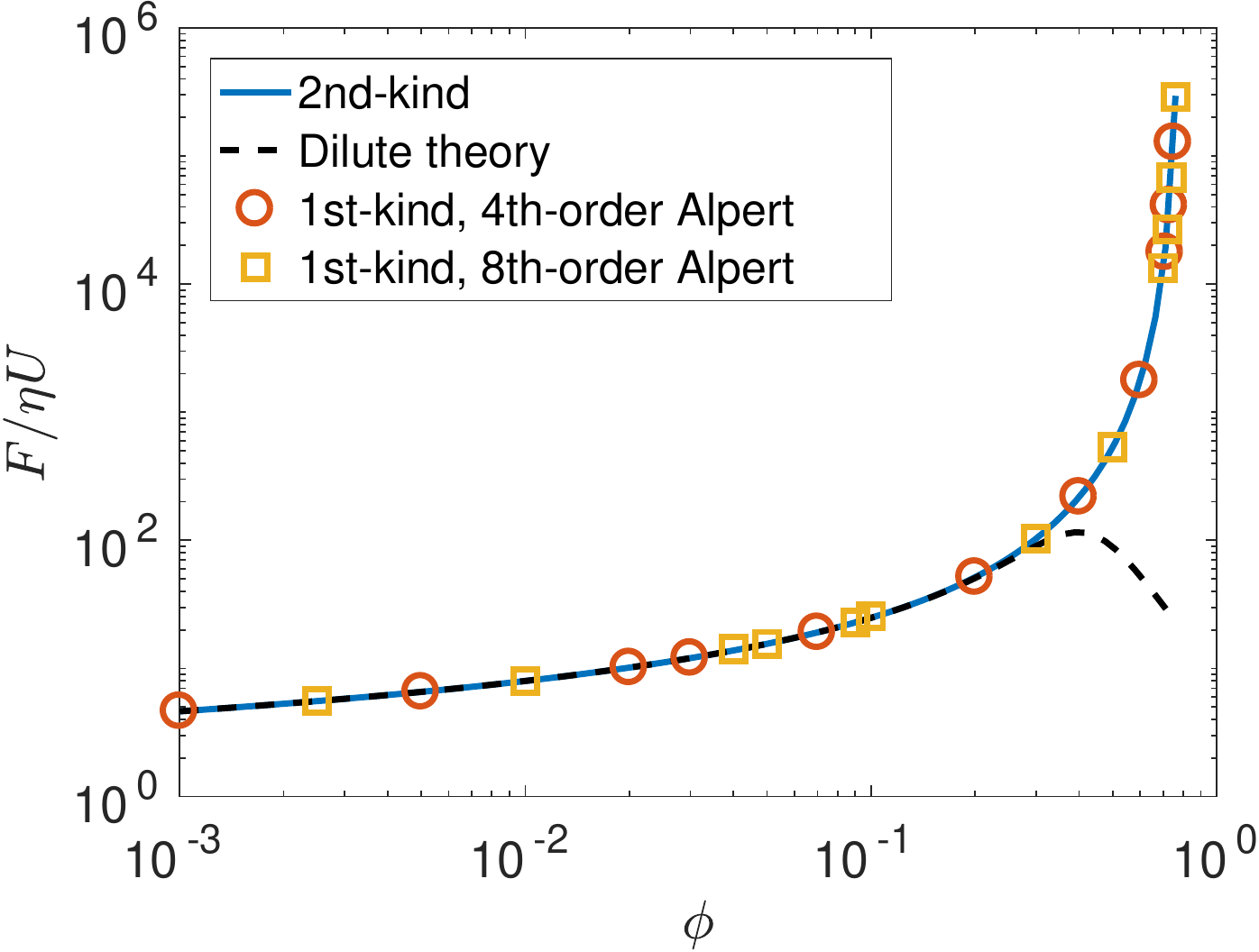} \hspace{1em}
 \includegraphics[width=.48\linewidth]{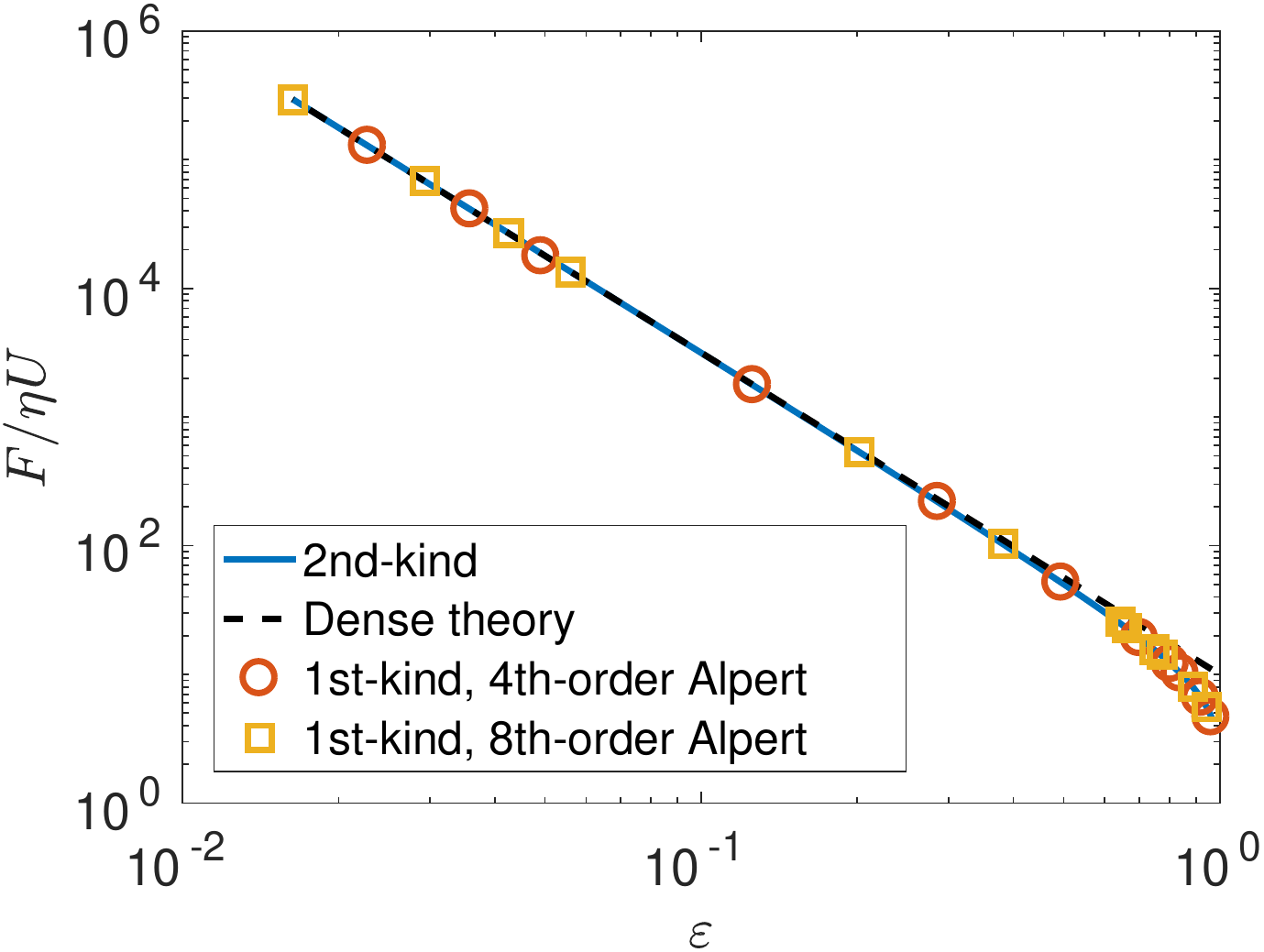} 
 \caption{The drag coefficients for a square periodic array of disks in steady Stokes flow, ({\it Left panel}) as a function of packing fraction $\phi$, and compared to the dilute theory \eqref{DiluteTheory}; ({\it Right panel}) as a function of the relative gap between disks $\epsilon$, and compared to the dense theory \eqref{DenseTheory}. In both panels, the numerical results obtained from the first-kind solver with Alpert quadrature match very well with the spectrally-accurate results from the second-kind solver. (For details of the second-kind boundary integral formulation of the Stokes BVP, the reader is referred to \ref{appendix:secondkind-formulation} of this paper.) }
 \label{fig:periodicArrayDrag}
\end{figure}

The linear system \eqref{discreteSaddlePoint} is solved by GMRES with block-diagonal preconditioning. Numerical results for the normalized drag force over a broad range of packing fractions are shown in \autoref{fig:periodicArrayDrag}. We have obtained very good agreement with the dilute theory \eqref{DiluteTheory} ($\phi < 0.2$). For the dense packing fractions, the dilute theory no longer provides a good description for the drag, but our solutions from the first-kind solver match the dense theory \eqref{DenseTheory} very well, as shown in the right panel of \autoref{fig:periodicArrayDrag}.
The solutions from the first-kind solver are also in excellent agreement with the highly-accurate solutions ( at least 9 digits of accuracy) from the second-kind solver for the range of packing fractions considered.

%

\subsection{Suspension of Brownian rigid disks}
\label{sec:SuspensionOfDisks}

\begin{figure}
\centering
\begin{minipage}{.5\textwidth}
\centering
 \includegraphics[width = .8\linewidth]{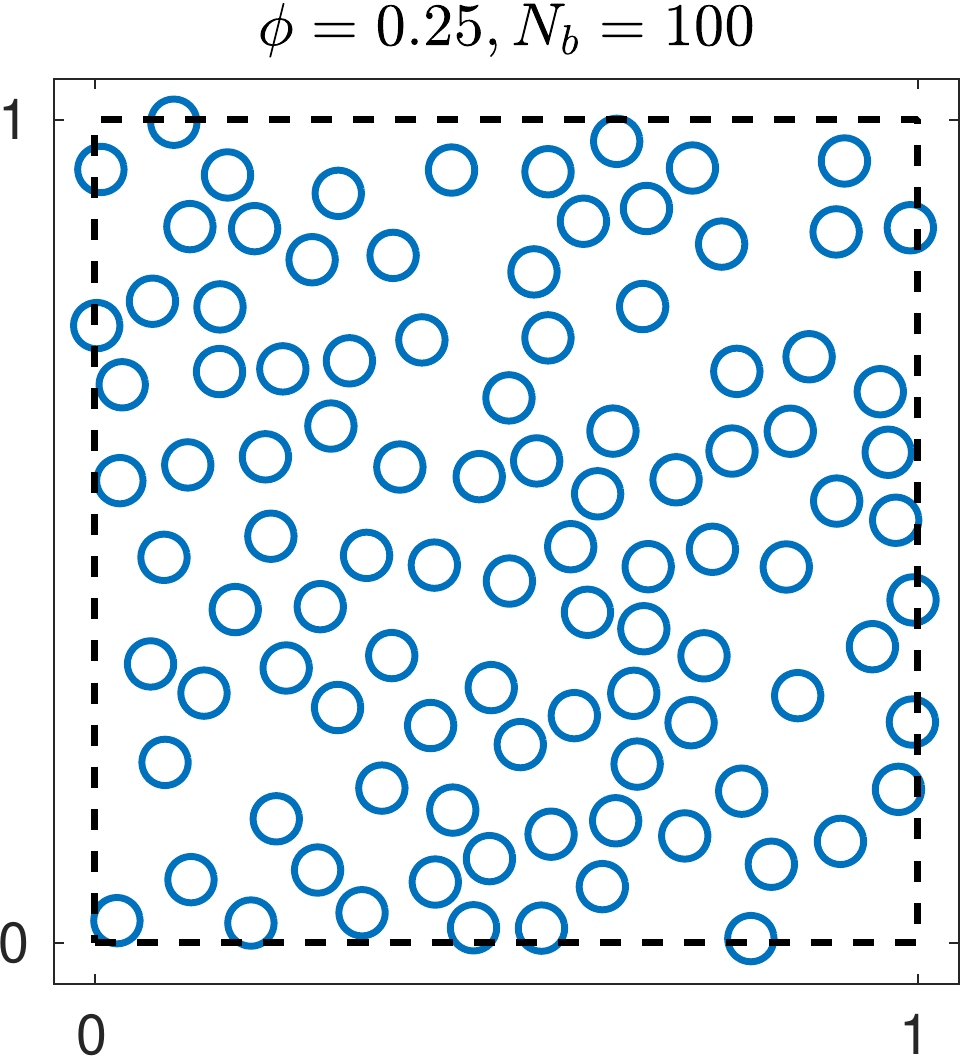} \\
 \includegraphics[width = 1.0\linewidth]{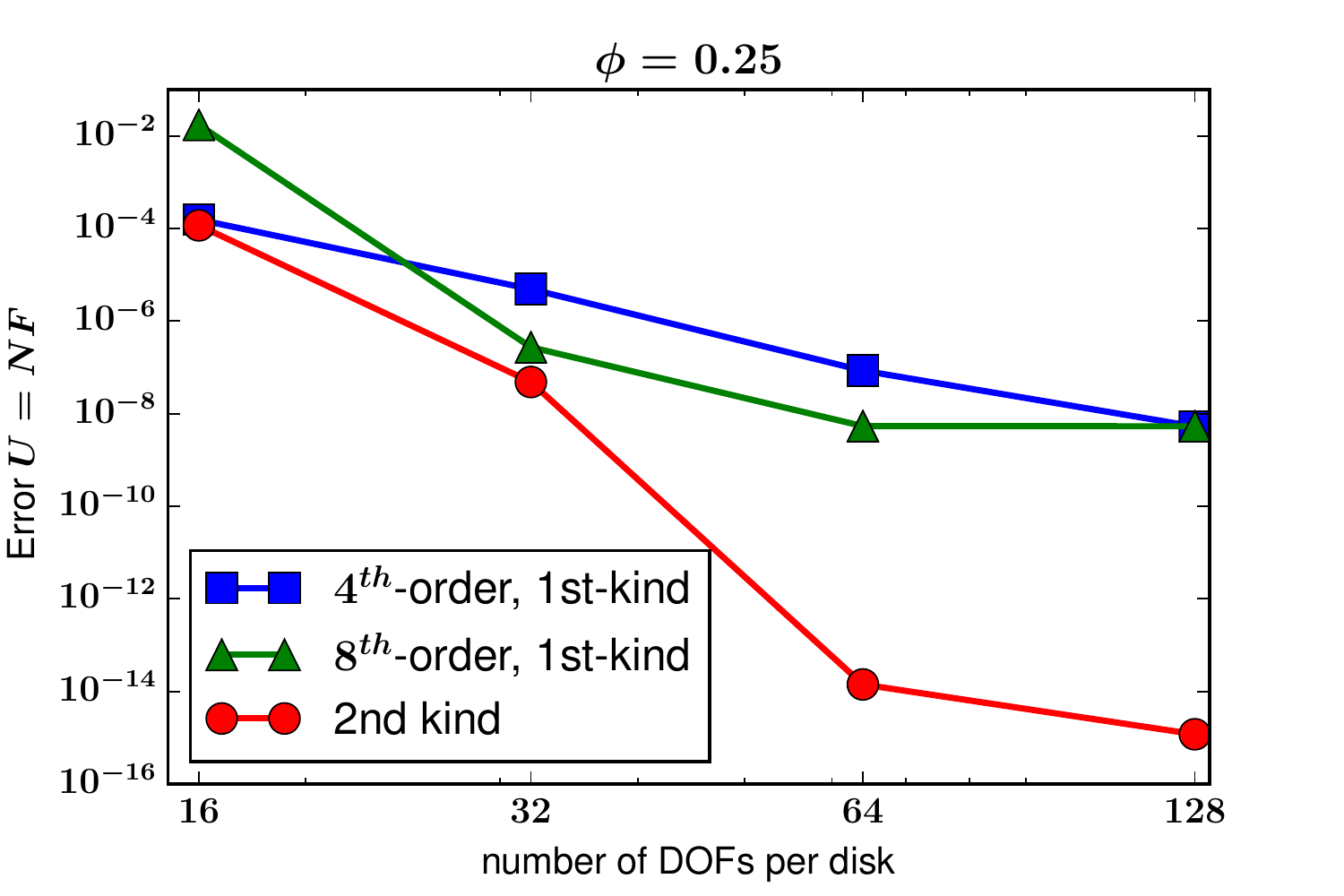}
\end{minipage}%
\begin{minipage}{.5\textwidth}
\centering
 \includegraphics[width = .8\linewidth]{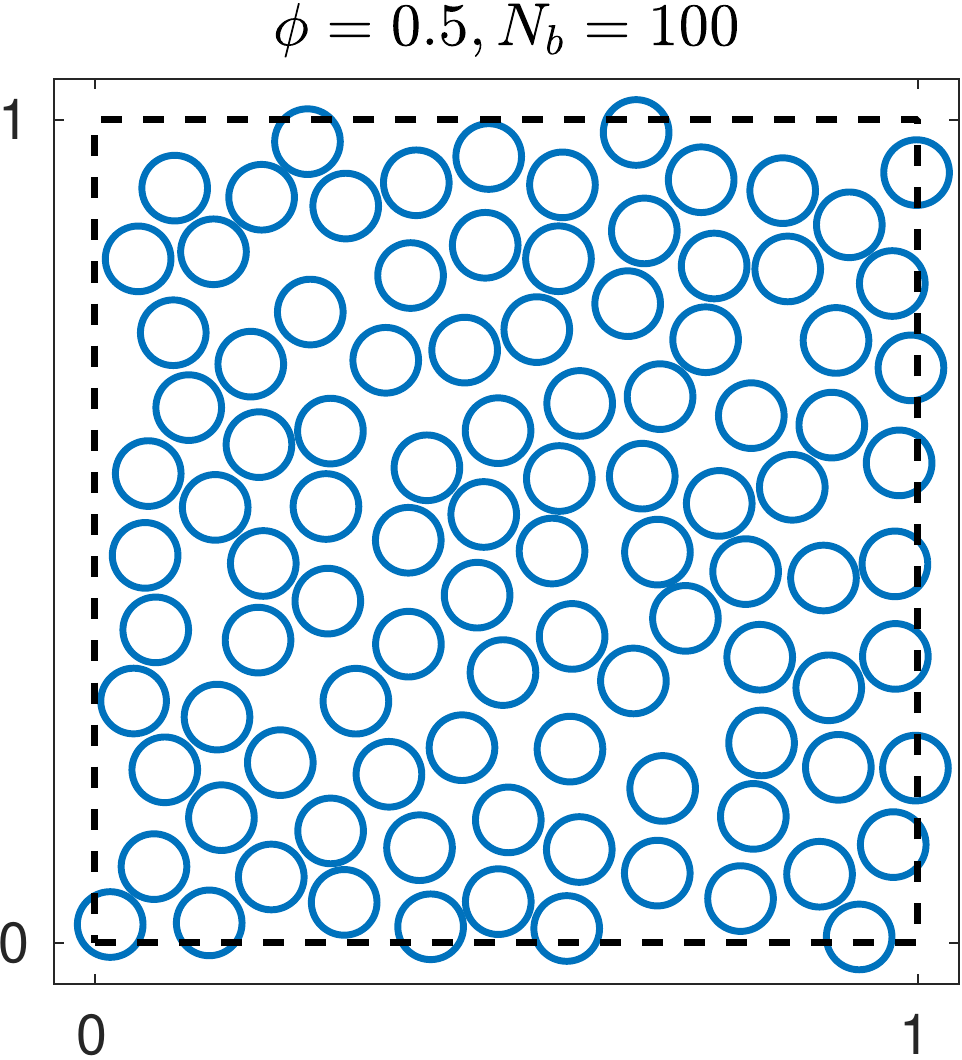} \\
 \includegraphics[width = 1.0\linewidth]{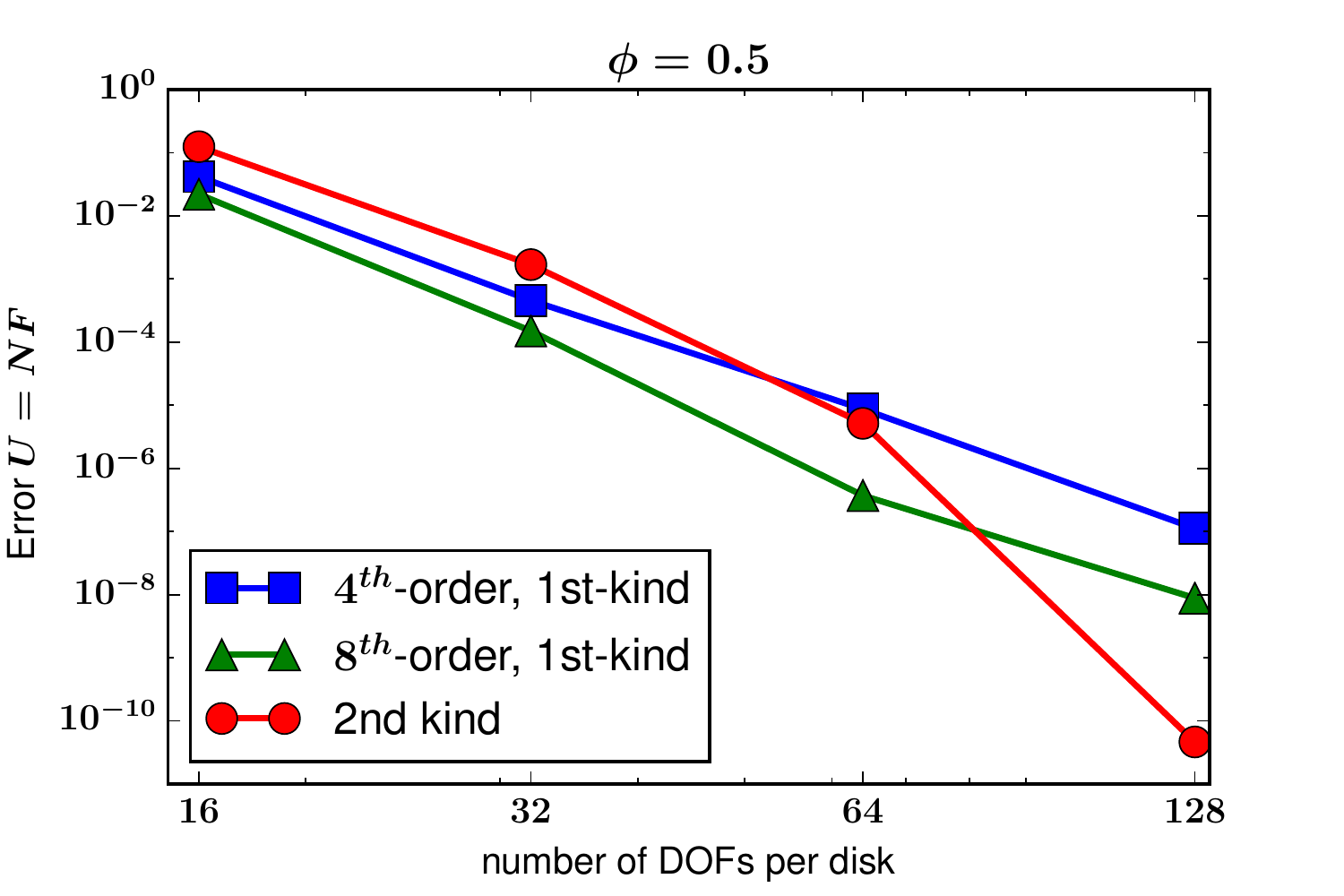}
\end{minipage}
 \caption{({\it Top panel}) Two random configurations of 100 rigid disks with packing fractions $\phi = 0.25$  (dilute) and $\phi = 0.5$ (dense) in a periodic unit cell. ({\it Bottom panel}) Normalized error of the mobility $\vf{U} = \vf{N}\vf{F}$ versus number of degrees of freedom (DOFs) per disk. For the dilute suspension, the second-kind solver converges spectrally fast and is generally more accurate than the first-kind solver. For the dense suspension, while the second-kind solver converges faster, the first-kind solver is actually more accurate at lower DOFs per disk. }
 \label{fig:accuracy-random-config}
\end{figure}

In this section our primary goal is to study the performance of FBIM applying to suspensions of Brownian rigid particles in two dimensions. 
We will assess the effectiveness of FBIM from three aspects: accuracy and convergence, robustness of iterative solvers, and its efficiency and scalability to many-body particle systems.
For simplicity, we focus on suspensions of disks only, however, there is no significant difficulty in applying the FBIM to more general bodies. 
In the top panel of \autoref{fig:accuracy-random-config}, we show two random configurations of $N=100$ disks with different packing fractions: 
a dilute suspension with $\phi = 0.25$ and a dense suspension with $\phi = 0.5$. These configurations are generated by using an event-driven molecular dynamics code. Since the random disks may nearly touch in the absence of (electrostatic) repulsive forces, especially when the packing fraction is high, we first generate a random configuration of disks with radius $a_0$ at a higher packing fraction $\phi_0$, and then adjust the actual radius of disks $a$ to achieve the desired packing fraction $\phi$. In this approach, the pairs of random disks are separated by a relative minimum distance $d_{\min}/a = 2\left(\sqrt{\phi_0/\phi}-1\right)$. For the dilute suspension ($\phi = 0.25$), we use $\phi_0 = 0.4$ so that $d_{\min}/a \approx 0.523 $, and for the dense suspension ($\phi = 0.5$), we use $\phi_0 = 0.6$ so that $d_{\min}/a \approx 0.191$.

First, we investigate the accuracy and convergence of the first-kind mobility solver by applying it to the random configurations of disks shown in  \autoref{fig:accuracy-random-config}, subject to random forces and torques $\vf{F}$ (without random surface velocity $\vslip$). In the first-kind mobility solver, we set $\tolEwald = \tolIter = 10^{-9}$, and choose the splitting parameter $\xi \approx 50$ (or $\Nbox=10^2$). Although this value of $\xi$ does not achieve the minimum CPU time (see \autoref{fig:profiling-scaling}), it ensures $\rAlpert / \rcutoff \lesssim 0.5$ for the $8^{th}$-order Alpert quadrature, so that the singularity is sufficiently resolved in all test cases for convergence study purpose. For this set of computations, we consider different numbers of quadrature nodes, or degrees of freedom (DOFs) per body, $N_p \in \left\{ 16,32, 64, 128\right\}$. The normalized error of $\vf{U} = \vf{N}\vf{F}$ is computed with respect to a 12-digit accurate solution computed from the second-kind solver with 256 DOFs per disk. For the dilute suspension (bottom left panel of \autoref{fig:accuracy-random-config}), the second-kind solver converges spectrally fast and is more accurate than the solutions from the first-kind solver. For the dense suspension (bottom right panel of \autoref{fig:accuracy-random-config}), while the second-kind solver converges faster and gives more accurate solutions for large number of DOFs per body, the first-kind solver is more accurate for low resolutions. This can be explained as follows. When two disks nearly touch, the singular kernel of the first-kind integral (Stokeslet) grows as $\log r$, but the singular kernel of the second-kind integral (stresslet) grows as $r^{-1}$  in two dimensions. This observation implies that the first-kind mobility solver is more practical than the second-kind solver for simulating suspensions of rigid particles with high packing fractions, since it can produce sufficiently accurate solutions with a smaller number of DOFs per body. 

\begin{figure}
\centering
 \includegraphics[width =  .8\linewidth]{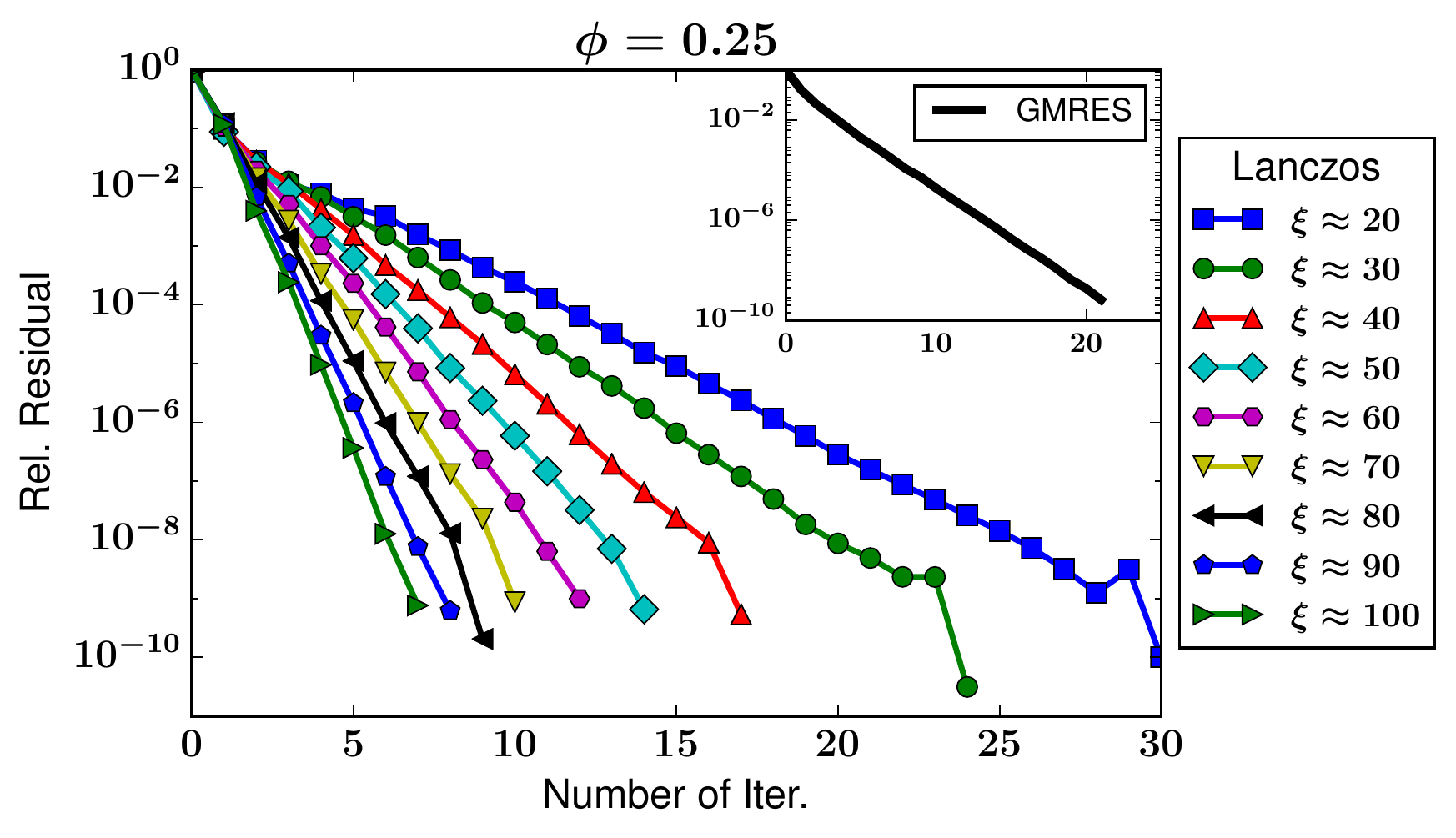} \\ 
  \includegraphics[width = .8\linewidth]{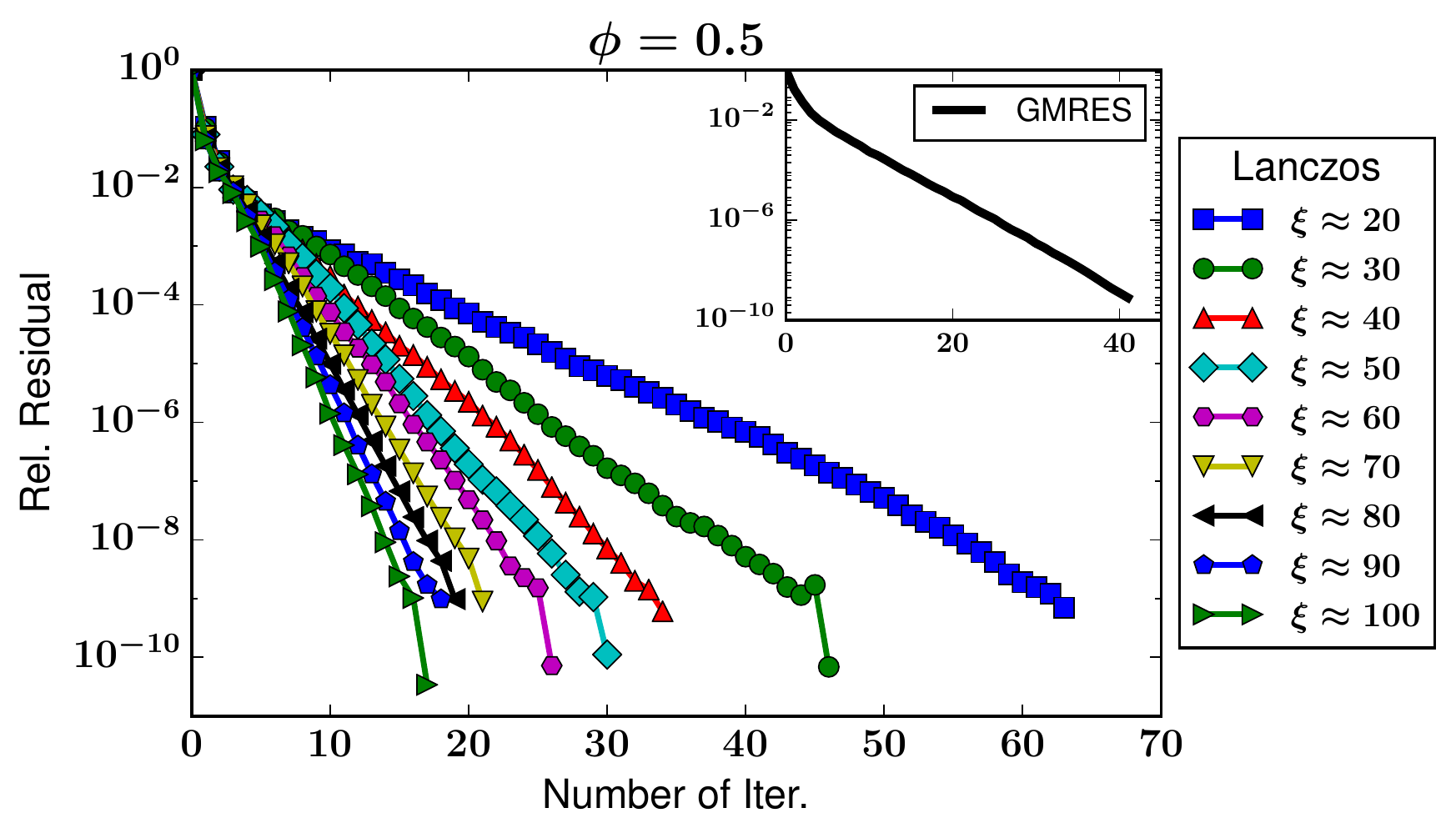} 
 \caption{Convergence of the Lanczos iteration for generating $(\mbf{M}^{(r)})^{1/2} \mbf{W}^{(r)}$ and the GMRES iteration (inset) for solving the saddle-point linear system \Cref{discreteSaddlePoint} with block-diagonal preconditioning for both iterative solvers. Generally, it requires more GMRES iterations for the dense packing ($\phi=0.5$), as shown in the inset. The number of GMRES iteration does not depend on $\xi$, since the choice of $\xi$ does not change $\mbf{M}$ in the linear system \eqref{discreteSaddlePoint}. The number of Lanczos iteration decreases with $\xi$, since the real-space kernel $\fStokeslet^{(r)}_{\xi}$ becomes more short-ranged as $\xi$ increases, and therefore, the block-diagonal preconditioner gets progressively better in approximating the inverse of $(\mbf{M}^{(r)})^{1/2}$.  }
 \label{fig:IterativeSolverConv}
\end{figure}

\begin{figure}
\centering
 \includegraphics[height = .34\linewidth]{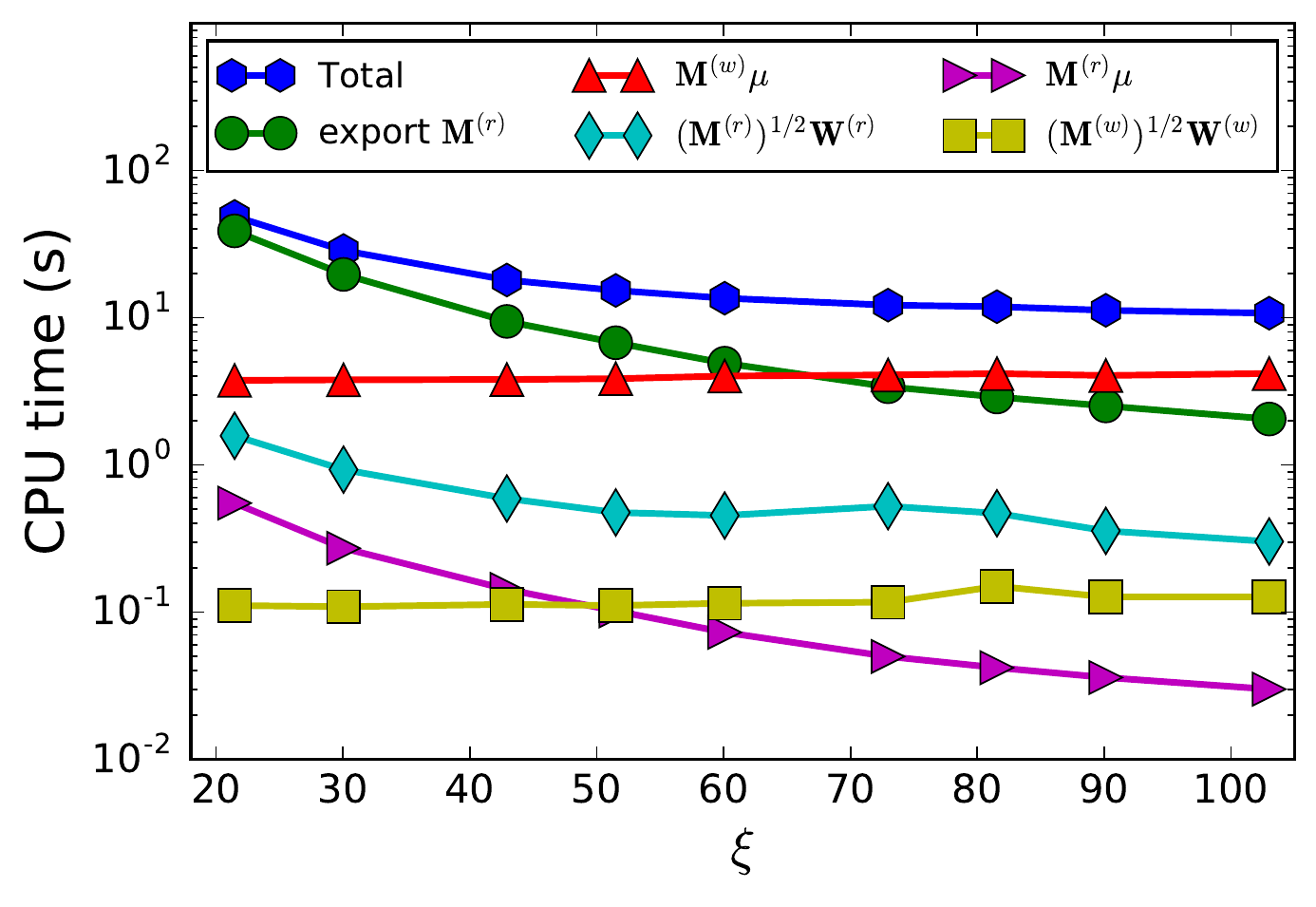} 
 \includegraphics[height = .34\linewidth]{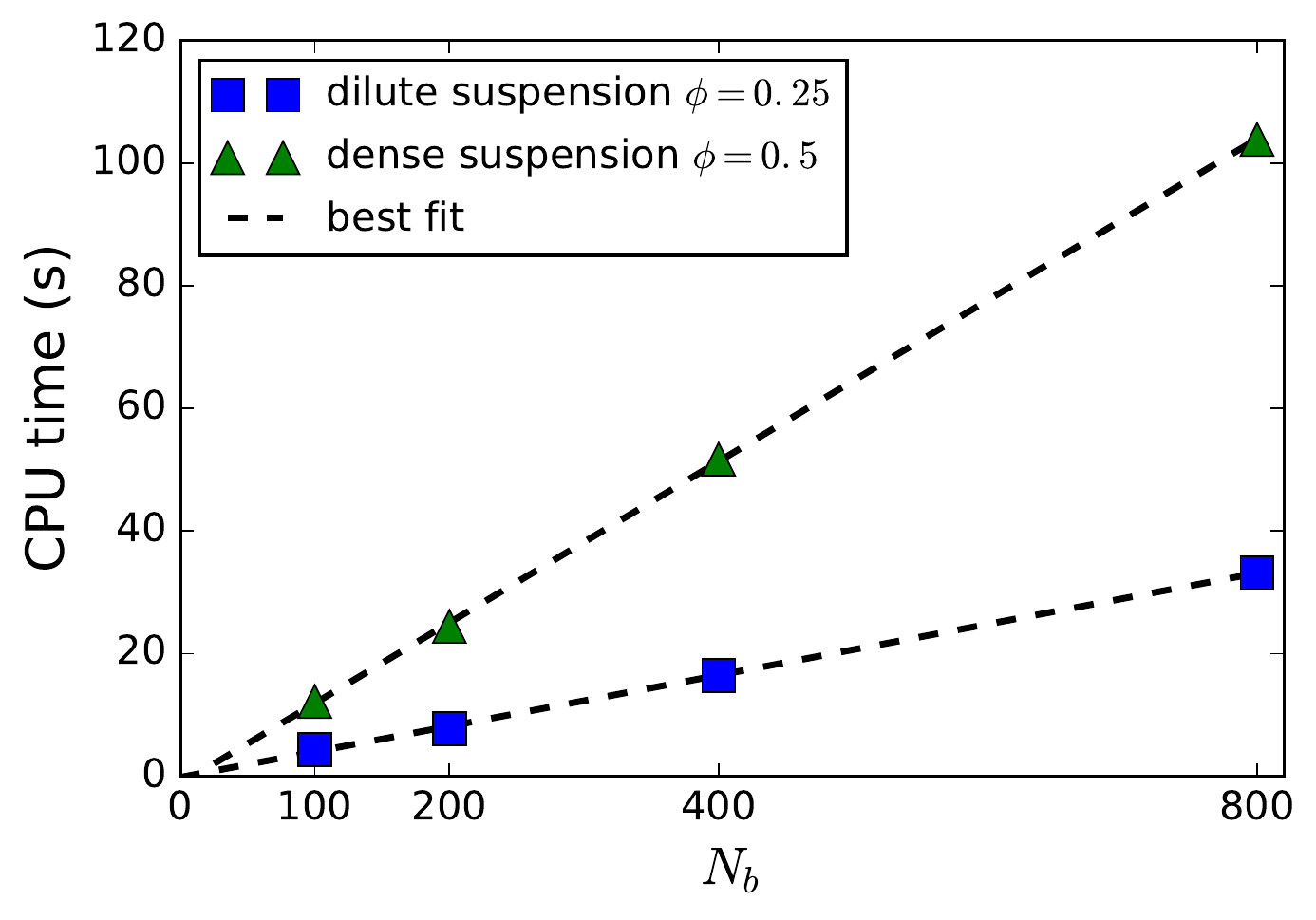} 
 \caption{({\it Left panel})  Execution time of FBIM and its algorithmic components versus $\xi$ for the dense packing $\phi=0.5$. In our implementation, the  computational work of FBIM is dominated by the subroutines for exporting $\mbf{M}^{(r)}$ sparsely and for evaluating the matrix-vector product $\mbf{M}^{(w)}\vf{\mu}$ (iteratively in GMRES). The total execution time of FBIM decreases with $\xi$, and the minimum CPU time is achieved near $\xi \approx 70$, where the CPU time for exporting $\mbf{M}^{(r)}$ and evaluating $\mbf{M}^{(w)} \vf{\mu}$ crosses. ({\it Right panel}) Linear scaling of FBIM with growing number of disks in the suspension, while the packing fraction is held fixed.}
 \label{fig:profiling-scaling}
\end{figure}

Next, we study the robustness of iterative methods for solving the saddle-point linear system \eqref{discreteSaddlePoint}. For the random configurations shown in \autoref{fig:accuracy-random-config}, in addition to random forces and torques, we also generate random surface velocity $\bfslip = \mbf{M}^{1/2}\mbf{W}$, and include it on the right-hand-side of \Cref{discreteSaddlePoint}. We show in the inset of \autoref{fig:IterativeSolverConv} the residual versus the number of GMRES iterations for solving  \Cref{discreteSaddlePoint} with block-diagonal preconditioning. It would take more than 3 times the number of GMRES iterations to converge to the same tolerance level without preconditioning. In general, when the packing fraction grows or when two disks get closer, the hydrodynamic interaction between the disks becomes stronger, and hence, the condition number of the linear system grows. As a result, it requires more GMRES iterations for the dense suspensions, as expected. The number of GMRES iteration is independent of $\xi$, since the choice of $\xi$ does not change $\mbf{M} = \mbf{M}^{(r)} + \mbf{M}^{(w)}$ in the linear system \eqref{discreteSaddlePoint}. We expect that GMRES will converge faster for three-dimensional problems because the Stokeslet decays faster ($r^{-1}$ in the far field) in three dimensions. 
We also show in \autoref{fig:IterativeSolverConv} the residual versus the number of Lanczos iterations for generating $(\mbf{M}^{(r)})^{1/2} \mbf{W}^{(r)}$ with block-diagonal preconditioning, for different values of $\xi$. The number of Lanczos iterations decreases with $\xi$ for both the dilute and dense suspensions. This can be explained by the observation that $\fStokeslet^{(r)}_{\xi}$ becomes more short-ranged as $\xi$ increases, and therefore, the block-diagonal preconditioner gets progressively better in approximating the inverse of $(\mbf{M}^{(r)})^{1/2}$.

We now present profiling and scaling analysis of FBIM.
First, we present the profiling results of FBIM and its algorithmic components in \autoref{fig:profiling-scaling}, and analyze the optimal choice of the splitting parameter $\xi$. We focus on the densely-packed configuration with 64 DOFs per disk ($d_s/d_g \approx 0.5$), and use the $4^{th}$-order Alpert quadrature. The accuracy of the first-kind mobility solver with these parameters is about $10^{-5}$ (see bottom right panel of \autoref{fig:accuracy-random-config}), so we set $\tolEwald = \tolIter = 10^{-6}$. We note that the execution time depends heavily on the choice of programming language, implementation and hardware. Our proof-of-concept serial implementation of FBIM in two dimensions is written in MATLAB with some subroutines accelerated by C with the aid of MEX files. As previously discussed in \autoref{sec:FBIM-method}, the main ingredients of FBIM are evaluating the matrix-vector products $\mbf{M}^{(r)} \vf{\mu}$ and $\mbf{M}^{(w)} \vf{\mu}$ in GMRES, and generating $(\mbf{M}^{(r)})^{1/2}\mbf{W}^{(r)}$ and $(\mbf{M}^{(w)})^{1/2}\mbf{W}^{(w)}$ using the Lanczos iteration.
In our implementation, we found it optimal to export and store $\mbf{M}^{(r)}$ sparsely for rapid matrix-vector multiplication in Lanczos and GMRES. In our profiling analysis, we profile the time to export $\mbf{M}^{(r)}$ sparsely, the cumulative time to evaluate $\mbf{M}^{(r)}\vf{\mu}$ and $\mbf{M}^{(w)}\vf{\mu}$ in GMRES separately, the time to generate $(\mbf{M}^{(r)})^{1/2}\mbf{W}^{(r)}$ and $(\mbf{M}^{(w)})^{1/2}\mbf{W}^{(w)}$ separately, and the total execution time of FBIM. We note that the total execution time of FBIM also includes the time of applying the preconditioners and other overhead time.  
The left panel of \autoref{fig:profiling-scaling} shows the profiling results of the densely-packed configuration for different values of $\xi$. 
First, we observe that our implementation of FBIM is dominated by two subroutines: exporting $\mbf{M}^{(r)}$ and evaluating $\mbf{M}^{(w)}\vf{\mu}$ (iteratively).  
The CPU time of the sparse matrix-vector product $\mbf{M}^{(r)}\vf{\mu}$ and generating $(\mbf{M}^{(w)})^{1/2}\mbf{W}^{(w)}$ (non-iteratively) is negligible.
The CPU time for generating $(\mbf{M}^{(r)})^{1/2} \mbf{W}^{(r)}$ is also small because of the rapid matrix-vector multiplication of the sparse matrix $\mbf{M}^{(r)}$. The CPU time for generating $(\mbf{M}^{(r)})^{1/2} \mbf{W}^{(r)}$ also includes the time for applying the preconditioner, which accounts for about 20\% - 30\% of the computational work.

Generally, the execution time of the real-space subroutines decreases with $\xi$, because the amount of work in the real space reduces as $\fStokeslet_{\xi}^{(r)}$ becomes more short-ranged with $\xi$. On the other hand, the execution time of Fourier-space subroutines remains almost constant for the range of $\xi$ considered in this test. This is because the cost of grid operations in Fourier-space sums (spreading/interpolation of a Gaussian to grid in NUFFT) dominates the cost of FFT in two dimensions. We expect the FFT cost would eventually become dominant in three dimensions. 
In our implementation, we observe that even the total CPU time of FBIM decreases with $\xi$, whereas Fiore {\it et al.} \cite{SpectralRPY} report an optimal $\xi$ for their GPU implementation of the PSE method in three dimensions. There is a wide range of $\xi$ that approximately minimizes the total CPU time (from $\xi \approx 70$ to $\xi \approx 100$), as shown in the left panel of \autoref{fig:profiling-scaling}. We recall, however, that the choice of $\xi$ in FBIM is also limited by the accuracy of Alpert quadrature. Using the criterion that $\rAlpert / \rcutoff \lesssim 0.6$ for the densely-packed configuration, we obtain that $\xi \lesssim 150$ for the $4^{th}$-order Alpert quadrature, and that $\xi \lesssim 75$ for the $8^{th}$-order Alpert quadrature. Similar optimal range of $\xi$ is also obtained in the profiling analysis for the dilute configuration.

Another important computational aspect that needs to be addressed is how the FBIM scales as the number of rigid particles grows while the packing fraction is held fixed. This aspect of FBIM is particularly important for applications involving a large number of particles. In the right panel of \autoref{fig:profiling-scaling}, we report the total execution time of FBIM with increasing number of rigid disks, while the packing fraction is held fixed at $\phi=0.25$ and $\phi=0.5$, and $\xi$ is also fixed ($\xi \approx 70$) for both configurations. We conclude that FBIM scales linearly in the number of rigid particles for both  dilute and dense suspensions.

%
%

\subsection{Brownian dynamics of rigid particles}
In the following tests, we combine the FBIM with stochastic temporal integrators to perform Brownian Dynamics (BD) for rigid particles. Applying the weakly first-order accurate Euler-Maruyama (EM) scheme to \Cref{BDeqns}, we obtain the BD algorithm, 
\begin{equation}
\begin{aligned}
\vf{Q}^{n+1} = \vf{Q}^n & + \Delta t  \vf{N}^n \vf{F}^n + \sqrt{2k_BT \Delta t} \, (\vf{N}^n)^{\frac{1}{2}} \mbf{W}^n \\
& + \Delta t \frac{k_BT}{\delta} \left[ \vf{N} \left( \vf{Q}^n + \frac{\delta}{2} \widetilde{\mbf{W}}^n \right)\widetilde{\mbf{W}}^n - \vf{N}\left( \vf{Q}^n - \frac{\delta}{2} \widetilde{\mbf{W}}^n \right) \widetilde{\mbf{W}}^n  \right],
\end{aligned}
\label{EM-RFD}
\end{equation} 
where $\Delta t$ is the time step size, the superscript denotes the time step level at which quantities are evaluated (e.g., $\vf{Q}^{n} = \vf{Q}(t=n\Delta t)$ and $\vf{N}^{n} = \vf{N}(\vf{Q}^n)$), $\delta$ is a small parameter, and $\mbf{W}^n$ and $\widetilde{\mbf{W}}^n$ are uncorrelated vectors of i.i.d standard Gaussian random variables. 

The last term in \Cref{EM-RFD} is a centered {\it random finite difference} (RFD) approximation to the stochastic drift term  that is equal in expectation to $(\Delta t \, \kbt ) (\partial_{\vf{Q}} \cdot \op{N})^n$ for sufficiently small $\delta$ \cite{BrownianBlobs,BrownianMultiBlobs,MagneticRollers}. 
The RFD term guarantees that the EM scheme is a consistent stochastic integrator of \Cref{BDeqns}, but is simpler and more efficient in practice than the Fixman midpoint scheme \cite{BD_Fixman}, in which the action of $\vf{N}^{-\frac{1}{2}}$ is required. The choice of $\delta$ is determined by a balance between truncation and roundoff error in the centered RFD, which gives $\delta / a \sim \epsilon^{1/3}$, where $\epsilon$ is the accuracy of the matrix-vector product $\vf{N} \vf{F}$ and $a$ is the characteristic length of the particle. We note that the EM scheme requires solving the saddle-point linear system \eqref{discreteSaddlePoint} three times: once for generating the velocity $ \vf{N}^n \vf{F}^n + \sqrt{2k_BT/ \Delta t} \, (\vf{N}^n)^{\frac{1}{2}} \mbf{W}^n$, and twice for the RFD approximation.

The EM scheme is not particularly accurate even for the deterministic motion. Another weakly first-order accurate temporal integrator that has been observed to give a better accuracy is the stochastic Adams-Bashforth (AB) scheme \cite{MagneticRollers}, 
\begin{equation}
\begin{aligned}
\vf{Q}^{n+1} = \vf{Q}^n & + \Delta t \left(  \frac{3}{2}\vf{N}^n \vf{F}^n - \frac{1}{2} \vf{N}^{n-1} \vf{F}^{n-1} \right) + \sqrt{2k_BT \Delta t} \, (\vf{N}^n)^{\frac{1}{2}} \mbf{W}^n \\
& + \Delta t \frac{k_BT}{\delta} \left[ \vf{N} \left( \vf{Q}^n + \frac{\delta}{2} \widetilde{\mbf{W}}^n \right)\widetilde{\mbf{W}}^n - \vf{N}\left( \vf{Q}^n - \frac{\delta}{2} \widetilde{\mbf{W}}^n \right) \widetilde{\mbf{W}}^n  \right],
\end{aligned}
\label{AB2-RFD}
\end{equation}
in which the deterministic mobility $ \op{N}\vf{F}$ in \Cref{BDeqns} is approximated by the second-order Adams-Bashforth approximation. We observe, however, that the AB scheme is more expensive to use in practice, because it requires four mobility problem solves instead of three in the EM scheme. More efficient schemes can be developed but are not the focus of our work \cite{BrownianMultiblobSuspensions}.

\begin{figure}
\centering
\subfloat[][]{\includegraphics[width=.35\linewidth]{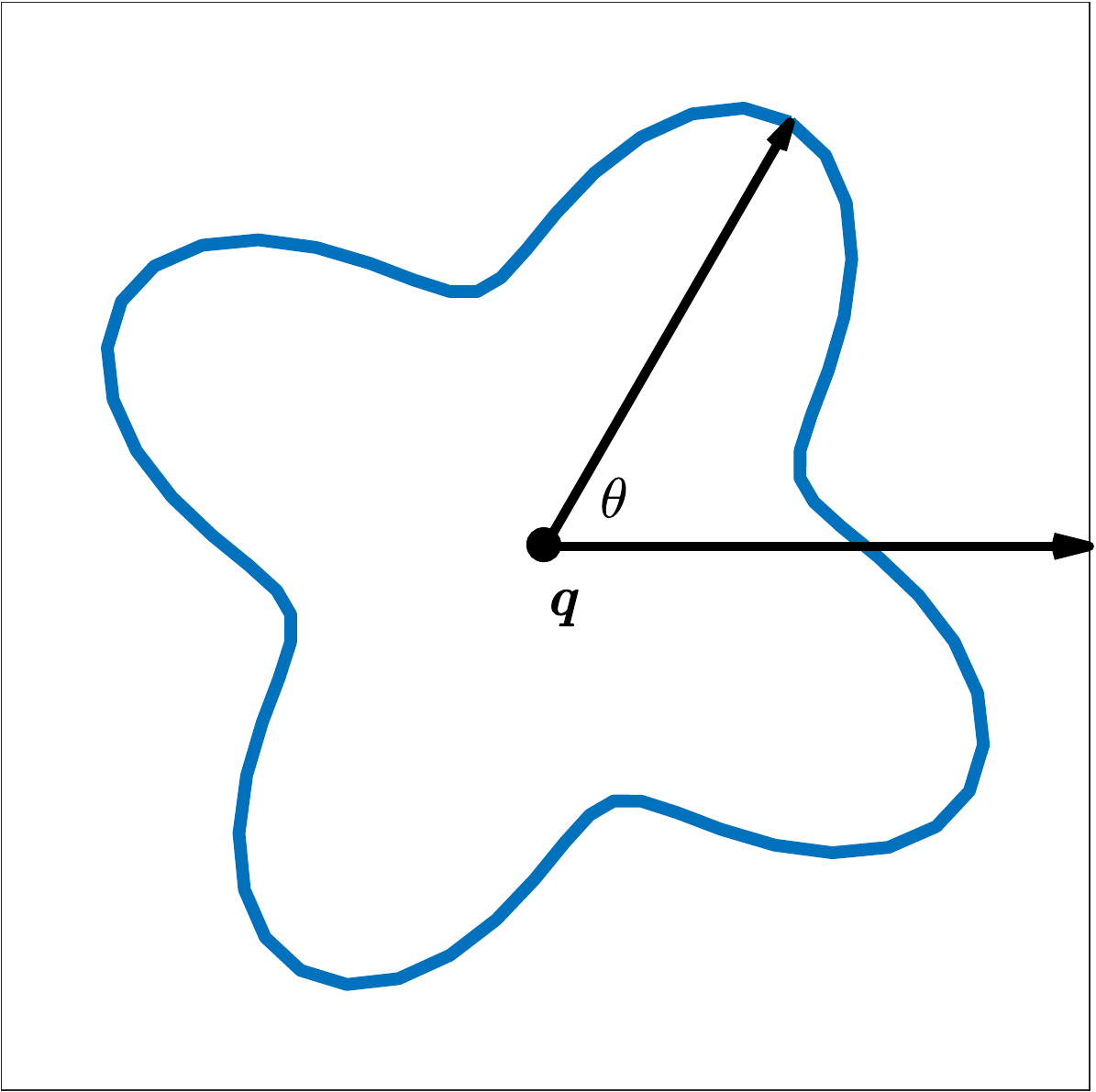}} \hspace{5em}
\subfloat[][]{\includegraphics[width=.35\linewidth]{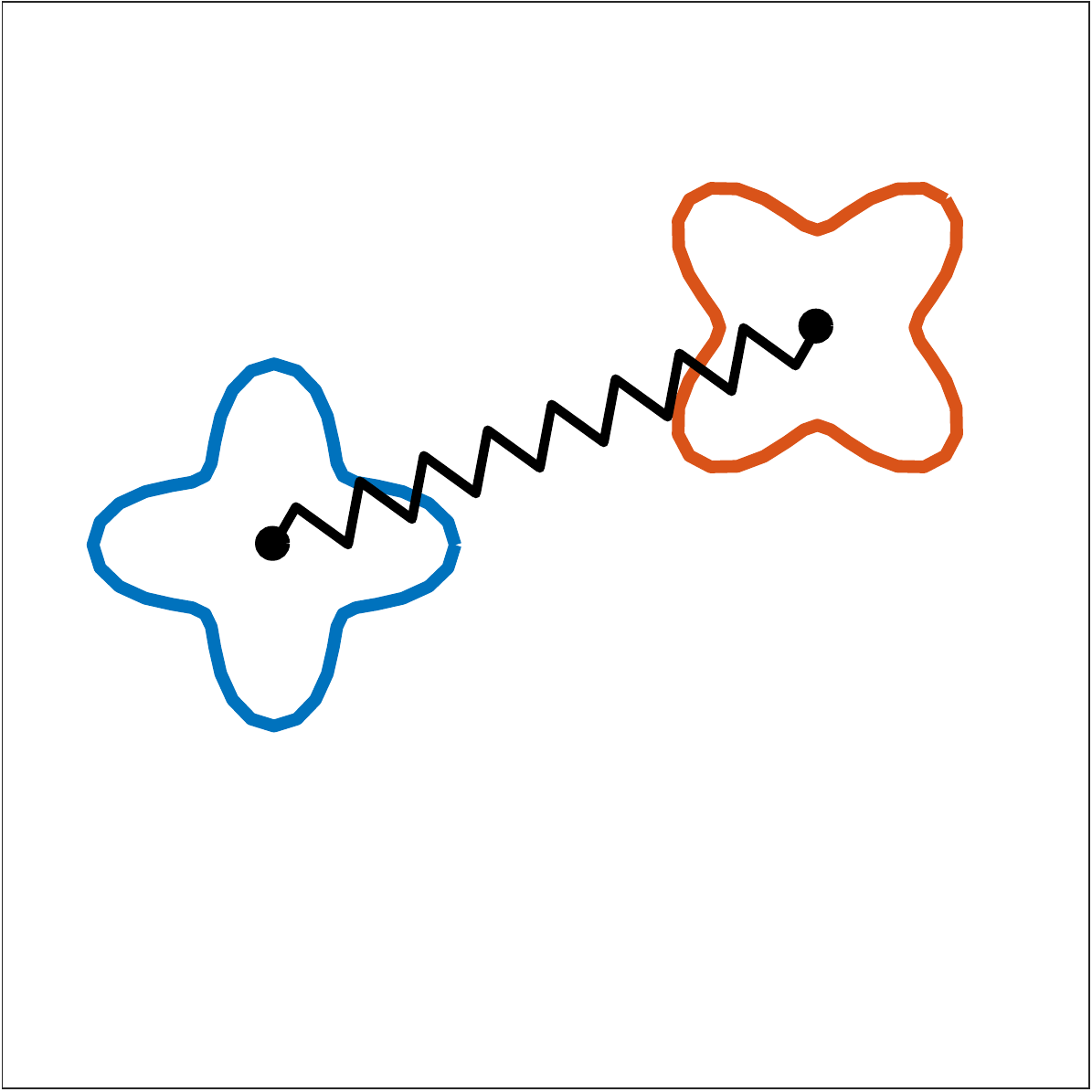}}
\caption{({\it Left panel}) A single starfish particle with four-fold rotational symmetry, described by the position of  tracking point $\vf{q}$ and its rotation $\theta$. ({\it Right panel}) Equilibrium position for a pair of starfish particles connected by a harmonic spring. The enclosing box shown in the figure is the periodic unit cell.}
\label{fig:starfish-config}
\end{figure}

\subsubsection{Free diffusion of a single starfish}
We consider a starfish-shaped particle freely diffusing (\ie, no applied force and torque) in a periodic unit lattice. 
A similar benchmark problem was studied by Delong {\it et al.} \cite{BrownianBlobs} using the FIB method, and by Delmotte {\it et al.} \cite{FluctuatingFCM_DC} using the FCM method for a single spherical particle in a periodic domain. In their case the body mobility matrix does not depend on the position of the particle, and therefore $\stochdrift{\vf{Q}}{\op{N}} = 0$. 
However, for a starfish $\op{N}$ depends on orientation due to interactions with periodic images and $\stochdrift{\vf{Q}}{\op{N}}$ is nonzero.


The starfish particle has four-fold rotational symmetry (left panel of \autoref{fig:starfish-config}), and is described by
\begin{equation}
\Gamma: \left(x(s), y(s) \right) = r_s(1+b\cos(4s)) \cdot ( \cos s, \sin s ), \quad s \in [0, 2\pi],
\end{equation}
where the characteristic length of the starfish particle  is its maximum radius $a = r_s(1+b)$.
In the continuum setting, the body mobility matrix of the four-fold starfish particle is a diagonal matrix and depends only on the rotation $\theta$, \ie, 
\begin{equation}
\op{N}(\theta) = 
\left[
\begin{array}{ccc}
\mathcal{N}_{xx}(\theta) & & \\
& \mathcal{N}_{yy}(\theta)  & \\
 & & \mathcal{N}_{\theta\theta}(\theta) \\
\end{array}
\right],
\end{equation}
where $\mathcal{N}_{xx}$, $\mathcal{N}_{yy}$ and  $\mathcal{N}_{\theta \theta}$ are the rotational and translational self-mobilities. Because of the symmetry of the starfish particle, we also have $\mathcal{N}_{xx} = \mathcal{N}_{yy}$.
The elements of $\op{N}(\theta)$ can be computed to 12 digits of accuracy using the second-kind solver by solving the deterministic mobility problem with $\vf{F}$ set to be the columns of the identity matrix for different $\theta$. 

Applying the EM scheme without RFD to the freely-diffusing starfish particle, we obtain the biased scheme,
\begin{equation}
\vf{Q}^{n+1} = \vf{Q}^n + \sqrt{2\kbt \Delta t} (\vf{N}^n)^{\frac{1}{2}}\mbf{W}^n.
\label{EM-noRFD}
\end{equation} 
In the limit $\Delta t \rightarrow 0$, the biased scheme \Cref{EM-noRFD} is consistent with the Ito SDE
\begin{equation}
\begin{aligned}
\frac{d\vf{Q}}{dt} 
&= \sqrt{2\kbt \Delta t} \left(\op{N}(\theta) \right)^{\frac{1}{2}}\op{W}(t) \\
&=  -\op{N}(\theta) (\partial_{\vf{Q}} \widetilde{U}) + \sqrt{2\kbt \Delta t} \left(\op{N}(\theta) \right)^{\frac{1}{2}}\op{W}(t) + (\kbt)  \stochdrift{\vf{Q}}{(\op{N}(\theta))}, 
\end{aligned}
\label{biased-SDE}
\end{equation}
where the bias potential is $\widetilde{U} \left( \vf{Q}=\{\vf{q},\theta \} \right) =  \kbt \log(\mathcal{N}_{\theta \theta} )$. By examining the corresponding Fokker-Plank equation, we can show that the biased SDE \eqref{biased-SDE} preserves the biased equilibrium distribution 
\begin{equation}
\widetilde{P}_{eq}(\vf{Q}) \equiv \widetilde{P}_{eq}(\theta) = {Z}^{-1} \exp \left( -\widetilde{U}(\theta)  / k_BT \right) =  \left( Z  \mathcal{N}_{\theta \theta}\right)^{-1}.
\label{biased-dist}
\end{equation}
We note that the correct Gibbs-Boltzmann distribution preserved by the unbiased scheme (EM with RFD) is a uniform distribution $P_{eq}(\vf{Q}) = \mathrm{constant}$.

\begin{figure}
\centering
\subfloat[][]{\includegraphics[width = .5\linewidth]{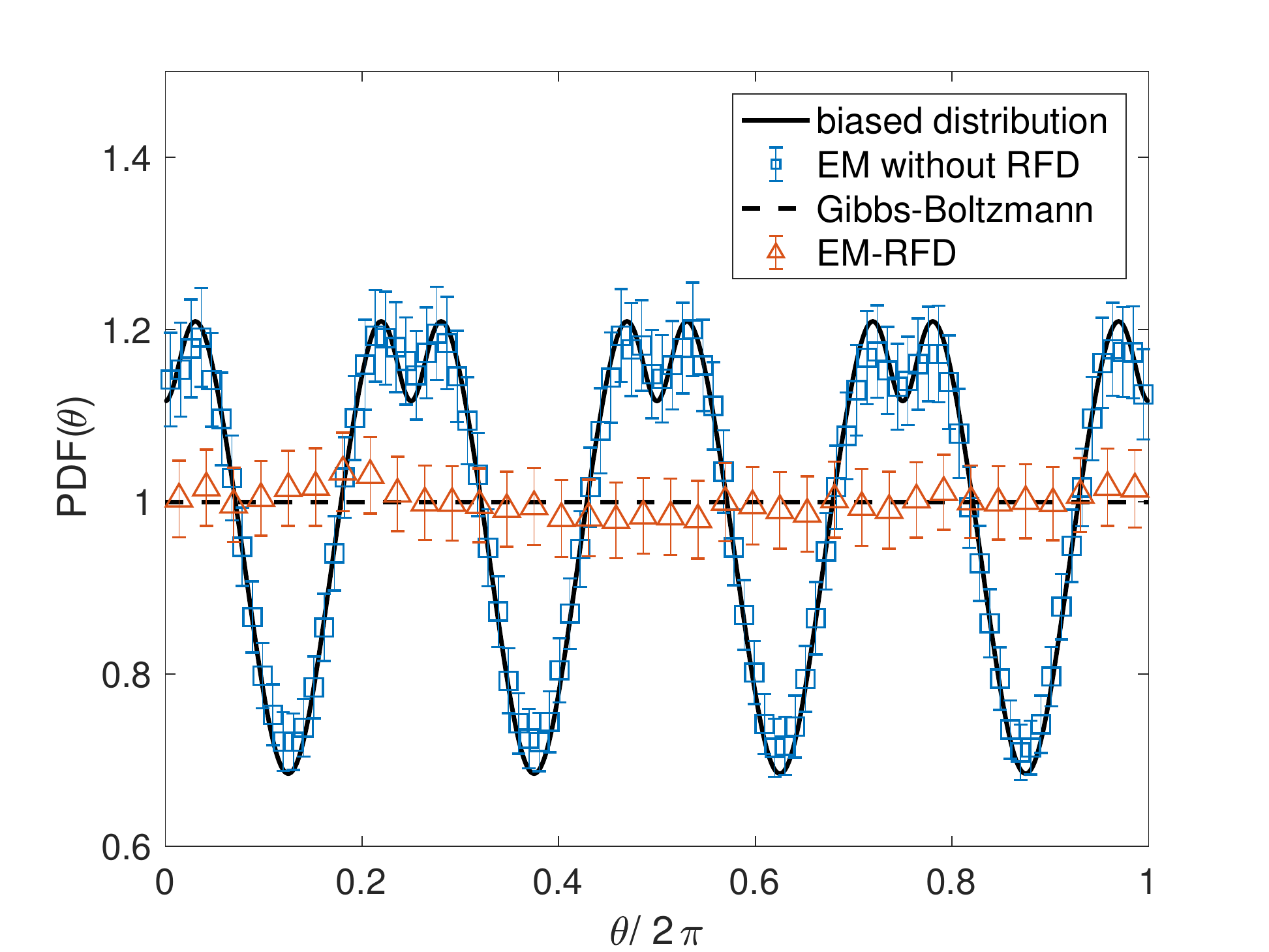}}
\subfloat[][]{\includegraphics[width = .5\linewidth]{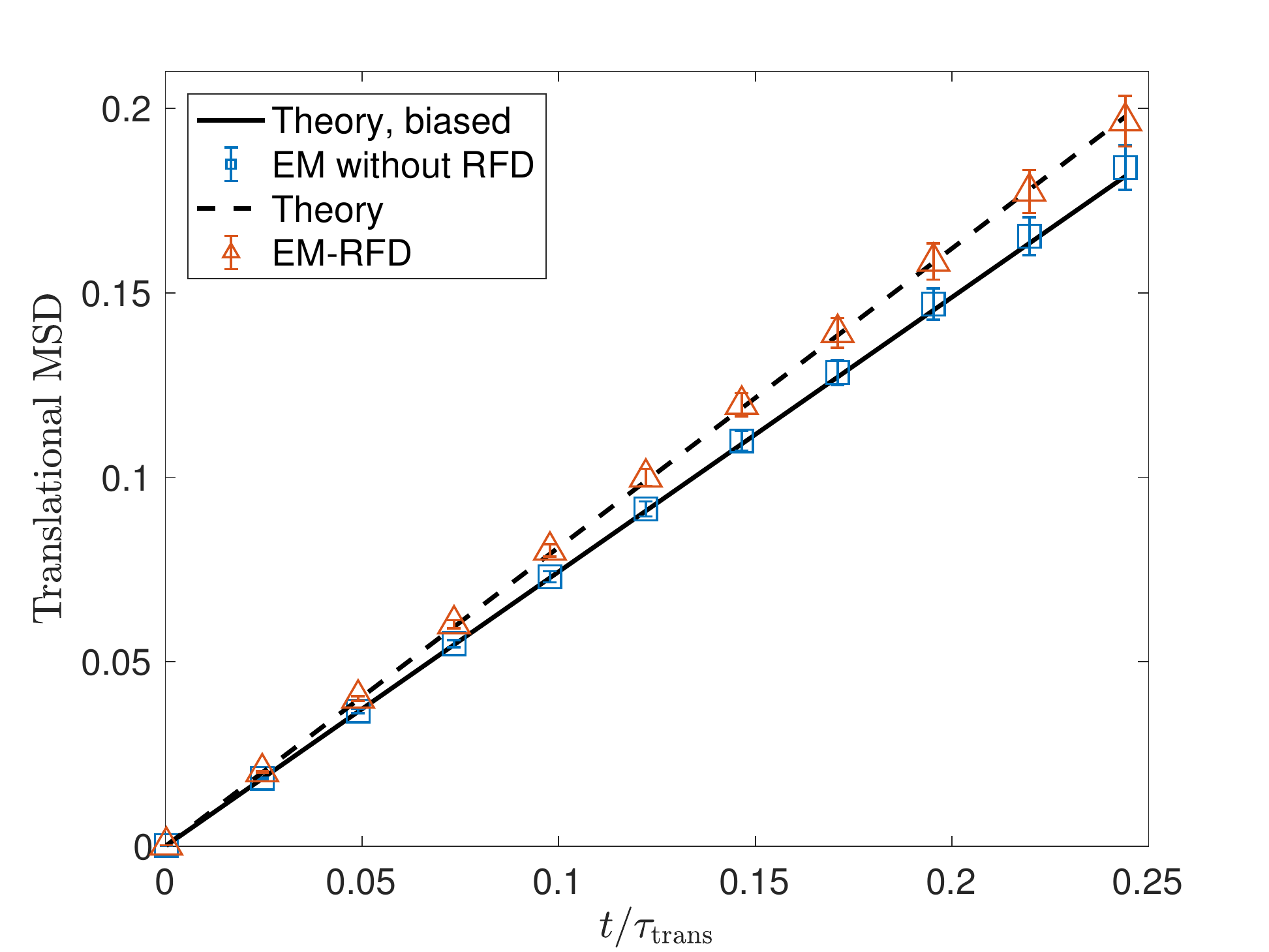}}
 \caption{({\it Left panel}) Equilibrium probability distribution of a single starfish particle freely diffusing in a periodic square lattice. The EM scheme without RFD produces a biased equilibrium distribution, and matches the theory \Cref{biased-dist} (solid). The EM scheme with RFD produces the correct Gibbs-Boltzmann distribution (dashed). ({\it Right panel}) Translational mean square displacement (MSD) for the freely-diffusing starfish particle using the biased and unbiased EM schemes.}
 \label{fig:EqDist-MSD}
\end{figure}

In our computation, we set $b=0.3$ and $a=1.3r_s= 0.45$, which gives a relatively high packing fraction $\phi \approx 0.393$ for the starfish particle, in order to amplify the difference between the biased distribution $\widetilde{P}_{eq}$ and the correct Gibbs-Boltzmann distribution $P_{eq}$ (see left panel of \autoref{fig:EqDist-MSD}). The starfish particle is discretized by $N_p = 64$ quadrature points using the $4^{th}$-order Alpert quadrature. In the first-kind mobility solver, we set the error tolerance level $\epsilon = 10^{-7}$. The RFD parameter $\delta$ is set by $\delta/ a \sim \epsilon^{1/3}$.  
The short-time translational $\chi_{\text{trans}}$ and rotational $\chi_{\text{rot}}$ diffusion coefficients  are determined by the Stokes-Einstein relations, 
\begin{equation}
\begin{aligned}
\chi_{\text{trans}} &= \kbt \langle \mathcal{N}_{xx}\rangle \approx 2.47 \times 10^{-3},   \\
\chi_{\text{rot}} &= \kbt \langle \mathcal{N}_{\theta \theta}\rangle \approx 2.081 \times 10^{-1},
\end{aligned}
\end{equation}
where the average $\langle \cdot \rangle$ is taken with respect to the equilibrium distribution. We use a small time step size $\Delta t = 0.02$ to minimize truncation errors, and each simulation is run for $T = 100$, where $\tau_{\text{rot}}$ is the rotational Brownian time scale.

Figure \ref{fig:EqDist-MSD} shows the estimated biased and unbiased equilibrium distributions (with error-bars of 2 standard deviations) obtained from 600 independent trajectories of the starfish particle freely diffusing in a periodic unit lattice. The numerical results using EM without RFD indeed match the biased equilibrium distribution $\widetilde{P}_{eq}(\theta)$, whereas EM with RFD preserves the correct Gibbs-Boltzmann distribution at equilibrium. In the right panel of \autoref{fig:EqDist-MSD}, we compare the translational mean square displacement (MSD)
for the biased and unbiased schemes, and also compare to the Einstein formula
\begin{equation}
\left\langle ||\vf{q}(t+s) - \vf{q}(t)||^2 \right\rangle =  4 \kbt \langle \mathcal{N}_{xx} \rangle \, t,
\end{equation}
where the average $\langle \cdot \rangle$ is taken with respect to the biased and unbiased equilibrium distribution, respectively. 
We conclude that consistent approximation to the stochastic drift term, such as the RFD approximation, is not only important for the equilibrium dynamics, but is also necessary  for producing the short-time dynamics (MSD) correctly. 

%

\subsubsection{A pair of interacting starfish}
In this test, we perform BD simulation of 
a pair of starfish particles connected by a Hookean spring with rest length $l_s$ in a periodic square lattice of length $l=2$ (right panel of \autoref{fig:starfish-config}). The two starfish particles interact through the spring  connecting  the two tracking points $\vf{q}_1, \vf{q}_2$ via a potential $U_{\text{spring}}$, and the rotation of each particle is also attached to a preferred angle through a harmonic potential $U_{\text{rot}}(\theta)$. The total potential $U(\vf{Q})$ is given by
\begin{equation}
\begin{aligned}
U(\vf{Q}) &= U(\vf{q}_1, \theta_1, \vf{q}_2, \theta_2) \\
 &= U_{\text{spring}}(\vf{q}_1, \vf{q}_2) + U_{\text{rot}}(\theta_1) + U_{\text{rot}}(\theta_2) \\
&= \frac{k_s}{2} ( |\vf{q}_1 - \vf{q}_2|- l_s)^2 + \frac{k_{\theta}}{2} \left(\theta_1 - \frac{\pi}{4}\right)^2+ \frac{k_{\theta}}{2} \left(\theta_2 - \frac{\pi}{2}\right)^2,
\end{aligned}
\label{total-potential}
\end{equation}
where $k_s$, $k_{\theta}$ are the stiffness coefficients. In the equilibrium, the Gibbs-Boltzmann preserved by \Cref{BDeqns} with $\vf{F} = -\partial_{\vf{Q}} U$ is
\begin{equation}
\begin{aligned}
P_{eq}(\vf{Q}) & \propto \exp \left( - U(\vf{Q})/ \kbt \right) \\
& \propto \exp \left( -\frac{U_{\text{spring}}(\vf{q}_1,\vf{q}_2)}{\kbt} \right) \cdot \mathscr{N}\left(\frac{\pi}{2}, \frac{\kbt}{k_{\theta}} \right) \cdot \mathscr{N}\left(\frac{\pi}{4}, \frac{\kbt}{k_{\theta}} \right),
\end{aligned}
\label{pair-starfish-dist}
\end{equation}
where $\mathscr{N}(\mu, \sigma^2)$ denotes the normal distribution with mean $\mu$ and variance $\sigma^2$.

In the computation, each particle has size $a/l = 1.3 \, r_s/l = 1/12$, and is discretized by $N_p = 48$ quadrature points using the $4^{th}$-order Alpert quadrature. 
The spring rest length is $l_s = 5a$ and the stiffness coefficient $k_s$ is set based on the criterion that 3 standard deviations of the distance $d = |\vf{q}_1 - \vf{q}_2 |$ is approximately $2a$, \ie, $3 \sqrt{\kbt/k_s} \approx 2a$, so that the particles very rarely overlap. The presence of harmonic springs defines the spring relaxation time scales,
\begin{equation}
\tau_{s} = \frac{1}{  k_s \langle \mathcal{N}_{xx} \rangle}  \quad \text{and} \quad
\tau_{\theta} = \frac{1}{  k_{\theta} \langle \mathcal{N}_{\theta \theta}  \rangle},
\end{equation}
where $\mathcal{N}_{xx}$ and $\mathcal{N}_{\theta \theta}$ are the translational and rotational self-mobilities of a single particle, which can be computed with high accuracy using the second-kind solver. 
The value of the parameters are $k_s = 81$, $k_{\theta} \approx 2.446$, and $\tau_s = \tau_{\theta} \approx 0.1189$.
Recall that the dominant computational work of FBIM is to compute the mobility $\vf{N}\vf{F}$ (see left panel of \autoref{fig:profiling-scaling}). To reduce the amount of computational work, we found that solving the mobility problems in RFD with a lower tolerance level $\epsilon=10^{-3}$ while maintaining a higher tolerance $\epsilon = 10^{-6}$ in the deterministic mobility does not introduce any observable statistical errors in the computation.

\begin{figure}
\centering
  \includegraphics[width = \linewidth]{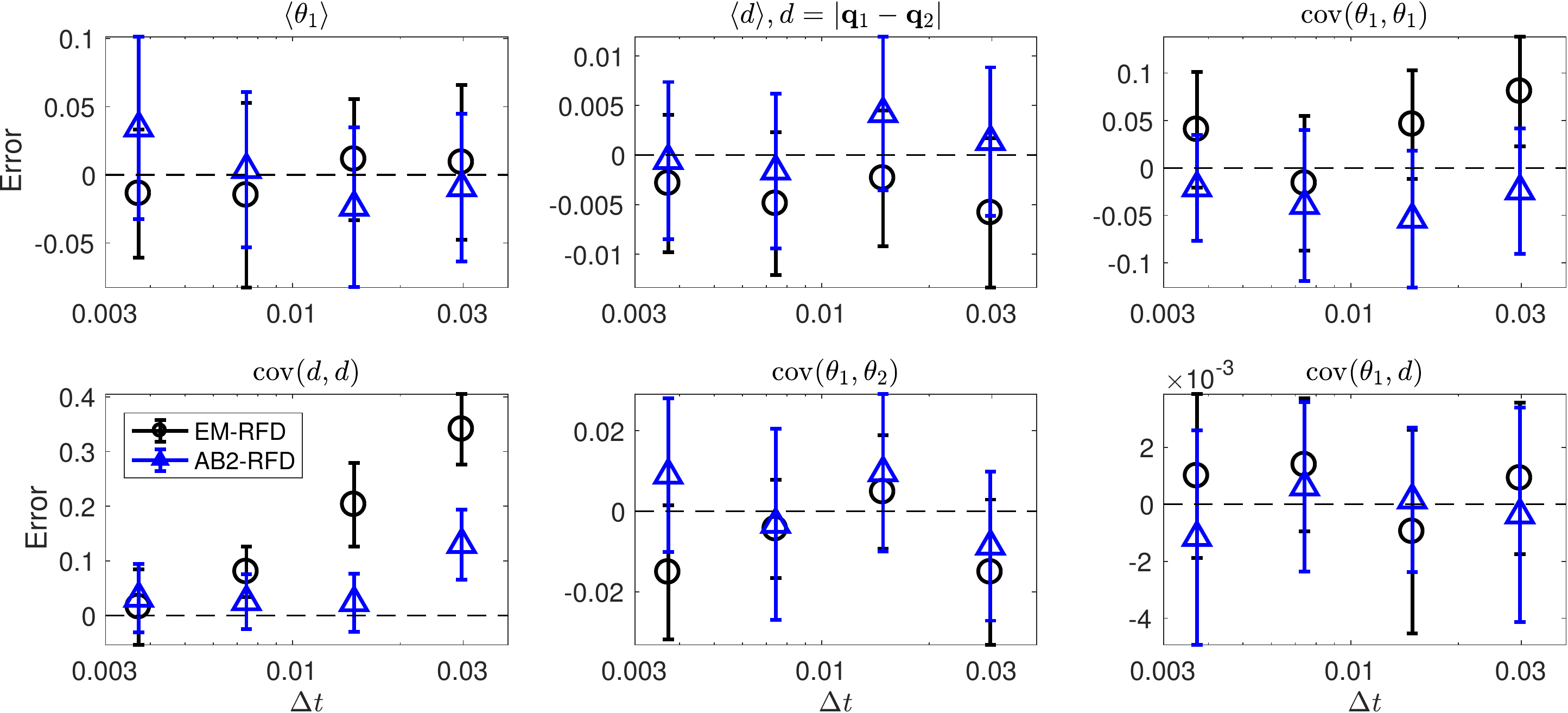} 
 \caption{Error bars (with two standard deviations) of the mean and covariance of distance between two tracking points $d = |\vf{q}_1 - \vf{q}_2|$ and rotations $\theta_1$, $\theta_2$ for a pair of starfish particles interacting with the potential \Cref{total-potential}. The numerical results are produced using the EM and AB schemes, with both including the RFD term. The AB scheme is more accurate than the EM scheme for $\mathrm{cov}(d,d)$. The cross covariances are statistically indistinguishable from zero, indicating the particle positions and their rotations  are uncorrelated, as expected from the equilibrium distribution \Cref{pair-starfish-dist}. }
 \label{fig:pairstarfish-error}
\end{figure}

In \autoref{fig:pairstarfish-error}, we compare errors of the mean and covariance of  $\theta_1$,  $\theta_2$ and the distance $d = |\vf{q}_1-\vf{q}_2|$ with respect to the equilibrium statistics calculated analytically from \Cref{pair-starfish-dist}. The numerical results are generated from 16 independent trajectories with length $T \approx 14.86$. We observe in   \autoref{fig:pairstarfish-error} that AB2 is more accurate than EM (see the panel that shows $\mathrm{cov}(d, d)$). We also observe that the cross covariances are statistically indistinguishable from zero for sufficiently small time step sizes, indicating that the distance between particles and their rotations are uncorrelated, as expected from the equilibrium distribution \Cref{pair-starfish-dist}.

\section{Conclusions}

In this paper we presented a fluctuating boundary integral method (FBIM) for simulating the overdamped Brownian Dynamics (BD) of rigid particles of complex shape in periodic domains.
To the best of our knowledge, it is the first boundary integral method that accounts for Brownian motion of nonspherical particles immersed in a viscous incompressible fluid.
Its main advantages are that particles of complex shape can be directly discretized with a surface mesh, and the deterministic mobility of the particles can be computed with high accuracy by using high-order singular quadrature techniques. Importantly, the Brownian displacements of the particles are computed along the way with only a marginal increase in the overall cost, and strictly satisfy discrete fluctuation-dissipation balance.
To accomplish this, instead of adding a stochastic stress tensor to the fluid equations as done in fluctuating hydrodynamics, we eliminated the fluid in the spirit of boundary integral representations. This lead to a Stochastic Stokes Boundary Value Problem (SSBVP) in which we prescribed a random surface velocity distribution that has covariance proportional to the (periodic) Stokeslet.
We found that using a first-kind boundary integral formulation is simplest since the first-kind integral operator inherits the SPD property from the Green's function for Stokes flow. The matrix discretizing the single-layer operator directly gives the covariance of the required fluctuating surface velocity, including a suitable handling of the singularity of the Oseen tensor.
While formal analysis suggests that second-kind boundary integral methods are to be preferred over first-kind methods, we found that first-kind methods can be more accurate for dense suspensions, and showed that a simple block-diagonal preconditioner can effectively handle the ill-conditioning of the first-kind formulation.
We confirmed through different benchmark problems that FBIM can efficiently compute both deterministic and Brownian motion of rigid particles in accordance with the order of accuracy of the singular quadrature scheme. Our preconditioned iterative solvers (GMRES and Lanczos) converged within a small number of iterations and only grew slowly with the packing fraction. We also confirmed that the computational cost of FBIM scales linearly with the number of particles, even for moderately dense packing fractions. Finally, we coupled FBIM with stochastic temporal integrators, and showed that it reproduced the correct equilibrium Gibbs-Boltzmann distribution of Brownian suspensions of particles of complex shape.

The FBIM presented in this work is a only a first step toward the overarching goal of performing accurate, efficient and robust BD simulations of a large collection of rigid particles of complex shape. Our proof-of-concept implementation of FBIM and numerical examples were presented in two spatial dimensions only. The continuum formulation of FBIM, which generates both deterministic and stochastic velocities in agreement with \Cref{DetPlusStochVelocities}, applies directly to three dimensions. In principle, our discrete formulation of FBIM can also be extended to three spatial dimensions. This extension requires, however, developing a suitable quadrature rule for the single-layer potential. While in two dimensions we were able to use the trapezoidal rule as the underlying quadrature rule and apply Alpert corrections to account for the singularity of the Oseen tensor, in three dimensions these pieces need to be developed anew. First, a suitable discretization of the particles' surfaces (e.g., using higher-order triangular elements) and a suitable non-singular quadrature rule need to be developed. For special particles shapes, notably spheres or spheroids, one can use specially chosen surface grids with the trapezoidal rule \cite{BoundaryIntegral_Periodic3D,BoundaryIntegral_SpheroidQBX}; however, for general particle shapes it is not straightforward to achieve spectral accuracy. One of the desired properties of the quadrature rule is to ensure that the far-field component of the discrete single-layer matrix $\mbf{M}$ is symmetric and positive definite. In particular, we expect that a good quadrature rule would yield the SPD matrix $\mbf{M}=\vf{\Psi}\mbf{G}\vf{\Psi}^T$, where $\mbf{G}_{ij} = \pStokeslet(\mbf{x}_i - \mbf{x}_j)$, and $\vf{\Psi}$ encodes the quadrature weights and mesh connectivity. Second, a singular quadrature near-field correction for the single-layer kernel needs to be constructed. A recently developed option is quadrature-by-expansion (QBX)  \cite{QBX_Epstein,QBX_Klockner}.
Recently, af Klinteberg and Tornberg \cite{BoundaryIntegral_SpheroidQBX} applied QBX to simulate a collection of non-Brownian spheroids immersed in a Stokes flow; however, the generalizations to more complex particle shapes requires an underlying, high-order smooth quadrature rule.

A particular challenge in developing singular quadrature schemes is preserving underlying symmetry properties of the single layer operator. Namely, the single layer integral operator is SPD and is rotationally invariant, meaning that if the body is rotated the result of applying the operator is also rotated in the same way. The Alpert quadrature correction used in this work gives a banded high-order log-singularity correction for the near-field component of $\mbf{M}$ that is neither symmetric, nor positive definite, nor rotationally invariant. We found all of these artifacts to be numerical errors below the order of accuracy of the quadrature scheme, as expected. However, it would be much better to have a singular quadrature scheme that does not have such unphysical artifacts {\em by construction}. The traditional focus in boundary integral methods has been on achieving higher-order accuracy.
What is more important for Brownian suspensions is preserving physical properties of the continuum operators in their discrete ``mimetic'' counterparts, so that even coarse resolutions give physically-consistent (even if not very accurate) discretizations.


{We note that a Galerkin boundary integral formulation combined with multipole expansions has been used by Singh and Adhikari and collaborators to model an active suspensions of spherical particles in an unbounded domain \cite{BoundaryIntegralGalerkin}. In some sense, this is an extension of the traditional Stokesian Dynamics method to include moments higher than the stresslet, as necessary to account for active flows. Recently the method has been extended to include Brownian contributions and to account for confinement near a no-slip boundary \cite{BoundaryIntegralWall_Adhikari}. While also based on a boundary integral formulation, this class of methods differs in significant ways from the one proposed here. Most importantly, our approach uses numerical quadrature instead of analytical integration, and therefore generalizes to arbitrary (smooth) particle shapes. Furthermore, we do not truncate a multipole hierarchy at a finite number of moments, and can therefore achieve controlled accuracy (\ie, a desired number of digits of accuracy). Nevertheless, generalizations of FBIM to three dimensions require a nontrivial amount of effort, and for spherical particles the Galerkin approach may be an effective alternative that yields sufficiently accurate answers in practice.}
\delete{Furthermore, we do not truncate a multipole hierarchy at a small number of moments (say,  after the dipole or quadrupole contribution), which can lead to uncontrolled loss of accuracy when the dilute approximation is no longer valid. In the FMM-accelerated boundary integral method, the number of terms used in multipole expansions is chosen to achieve a user-specified tolerance in the far field. High-order quadrature is used to achieve the same tolerance in the near field. That said, generalizations of FBIM to three dimensions require a nontrivial amount of effort, and for spherical particle the Galerkin approach mentioned above may be an effective alternative that yields sufficiently accurate answers in practice.}

Another key aspect of Brownian Dynamics simulations is the temporal integration schemes. There are three main challenges in this regard. The first is that in three dimensions orientation cannot be represented by a vector. This can most straightforwardly be handled by using normalized quaternions to represent particle orientation, as shown in \cite{BrownianMultiBlobs}. In the end, as long as a method to compute particle velocities in agreement with \Cref{DetPlusStochVelocities} is provided, handling the quaternion constraint simply amounts to using quaternion multiplication, rather than addition, to update orientations \cite{BrownianMultiBlobs,BrownianMultiblobSuspensions}. The second challenge is capturing the stochastic drift term proportional to the divergence of the body mobility matrix using linear-scaling iterative methods. In this work we used methods that perform random finite difference on the body mobility. This requires solving two mobility problems per time step just to capture the stochastic drift term, therefore at least doubling the cost of a time step. In future work we will describe novel temporal integrators that can be used with FBIM to give more accurate answers for larger time step sizes, and which use only a single mobility solve to capture the stochastic drift term \cite{BrownianMultiblobSuspensions}. A third challenge is to handle the fact that Brownian displacements can lead to overlaps between the particles, even for small time step sizes. Unlike the rigid multiblob method \cite{RigidMultiblobs}, which builds on a regularized first-kind boundary integral formulation, traditional boundary integral methods based on singular quadratures break down when particles overlap. Even if particles do not overlap, unless the surface discretization is refined adaptively (which would be too costly for denser suspensions), traditional boundary integral representations will give unphysical answers that can easily lead to breakdown especially in the presence of noise. A possible solution to this is to use a regularization of the formulation for particle gaps below the resolution of the method. This is done naturally in the rigid multiblob method \cite{BrownianMultiblobSuspensions}, and also more recently in a Galerkin multipole method \cite{BoundaryIntegralWall_Adhikari}, by using the Rotne-Prager-Yamakawa regularization of the Oseen tensor.

Another important direction of future work is to extend FBIM to model Brownian suspensions in other geometries, notably in unbounded domains. One option is to adapt our method to use a newly developed spectral Ewald summation for the free-space Stokeslet \cite{SpectralEwald_FreeSpace}. We note, however, that the computational cost of the FFTs in Ewald-type methods may become non-trivial for three dimensional problems. A challenge of great interest is incorporating the Fast Multipole Method (FMM) to generate the random surface velocity, thus developing a grid-free (near) linear-scaling method for Brownian suspension in an unbounded domain. For simple confined geometries, such as colloids sedimented in the vicinity of a single no-slip wall, one can use known analytical Green's functions for Stokes flow as done in \cite{BoundaryIntegralWall_Adhikari}. It turns out that because of the more rapid decay of the Green's function in the presence of a wall, iterative methods can efficiently generate the random surface velocity without requiring special handling of the far-field interactions \cite{MagneticRollers}. Nevertheless, achieving both linear scaling and controlled accuracy are challenges even in such simple geometries. For finite domains of more complicated geometry, one can discretize the domain boundary explicitly \cite{BoundaryIntegral_Wall}, and then employ the free-space Spectral Ewald method \cite{SpectralEwald_FreeSpace}.

\section*{Acknowledgements}

We gratefully thank Anna-Karin Tornberg for numerous informative discussions and for sharing with us notes on the Hasimoto splitting of the Stokeslet and Stresslet in two dimensions.
YB and AD were supported in part by the National Science Foundation under award DMS-1418706, and by the U.S. Department of Energy Office of Science, Office of Advanced Scientific Computing Research, Applied Mathematics program under Award Number DE-SC0008271.  EEK acknowledges support from EPSRC grant EP/P013651/1. 
LG and MR were supported in part by the
Office of the Assistant Secretary of Defense for Research and Engineering
and AFOSR under NSSEFF Program Award FA9550-10-1-0180 and in part
by the U.S. Department of Energy under contract DEFG0288ER25053.
MR was also supported in part by the Office of Naval Research under award N00014-14-1-0797/16-1-2123.

\appendix
\section{Body mobility matrix and Lorentz reciprocal theorem}
\label{appendix:LRT}

In this appendix, we derive the following expression for elements of the body mobility matrix $\op{N}$ in terms of the periodic Green's function $\pStokeslet(\vf{x}, \vf{y})$ of the Stokes equation,
\begin{equation}
(\op{N})_{ij} \equiv \mathcal{N}_{ij} = \int_{\Gamma} \int_{\Gamma} \vf{\lambda}^{(j)}(\vf{x}) \cdot \pStokeslet(\vf{x} -  \vf{y}) \cdot \vf{\lambda}^{(i)}(\vf{y}) \diff S_{\vf{y}} \diff S_{\vf{x}},\label{appendix-Nij}
\end{equation}
where the precise definition of $\vf{\lambda}^{(i)}, \vf{\lambda}^{(j)}$ appears later. 

For simplicity, we consider only a single rigid body $\Gamma$ immersed in a Stokes fluid with periodic boundary conditions. The generalization of \eqref{appendix-Nij} to account for many bodies is straightforward.  First, we recall that, the mobility problem that solves for the the  translational velocity $\vf{u}$ and the rotational velocity $\vf{\omega}$ of the body, in response to the force $\vf{f}$ and torque $\vf{\tau}$ exerted on the body  is described by steady Stokes equation with no-slip boundary condition, and force and torque balance conditions:
\begin{equation}
\begin{aligned}
-\div \vf{\sigma} &= \grad \pi - \eta \lapl \vf{v} = \vf{0}, \\
\div \vf{v} &= 0, \\
\vf{v}(\vf{x}) &= \vf{u} + \vf{\omega} \times (\vf{x}-\vf{q}), \quad \forall \vf{x} \in \Gamma, \\
\int_{\Gamma} \vf{\lambda}(\vf{x}) \diff S_{\vf{x}} &= \vf{f} \quad \text{and} \quad \int_{\Gamma} (\vf{x}-\vf{q}) \times \vf{\lambda}(\vf{x})\diff S_{\vf{x}} = \vf{\tau},
\end{aligned}
\label{appendix-mobility}
\end{equation}
where $\eta$ is the fluid viscosity,  $\vf{\sigma}$ the fluid stress tensor, and $\vf{v}$ is the fluid velocity, respectively.  Here $\vf{\lambda} = (\vf{\sigma} \cdot \vf{n})(\vf{x})$ is the surface traction of the body with $\vf{n}$ being the unit normal vector to the surface. The mobility problem \eqref{appendix-mobility}  can be viewed as a linear mapping 
\begin{equation}
\vf{U} = \op{N} \vf{F}, \label{appendix-RBM}
\end{equation}
where $\op N$ is the body mobility matrix that relates the rigid body motion $\vf{U}= \{ \vf{u},  \vf{\omega} \}$ to the applied force and torque $\vf{F}= \{ \vf{f},  \vf{\tau} \}$.

Let $\{ \vf{v}^{(i)}, \vf{\sigma}^{(i)}, \vf{u}^{(i)}, \vf{\omega}^{(i)} \}$ and $\{ \vf{v}^{(j)}, \vf{\sigma}^{(j)}, \vf{u}^{(j)}, \vf{\omega}^{(j)} \}$ denote the solutions to the mobility problem \eqref{appendix-mobility} with applied force and torque $\vf{F} = \{ \vf{f}^{(i)}, \vf{\tau}^{(i)} \} = \vf{e}^{(i)}$ and $\vf{F} = \{ \vf{f}^{(j)}, \vf{\tau}^{(j)} \} = \vf{e}^{(j)}$ respectively. Here $\vf{e}^{(i)}$ and $\vf{e}^{(j)}$ are the standard basis vectors (in two dimensions, $\vf{e}^{(i)} \in \mathbb{R}^3$).
It is not difficult to see that $\vf{U}^{(j)} = \{ \vf{u}^{(j)}, \vf{\omega}^{(j)} \}$ corresponds to the $j^{\text{th}}$ column of $\op{N}$.
It is understood that whenever the force or torque is made dimensionless in the canonical problem, the other quantities' dimensions are adjusted accordingly. Thus, if the force is made dimensionless, velocities have units of the force-velocity mobility, while the tractions have units of inverse area.  If the torque is made dimensionless, then the velocities have units of the torque-velocity mobility and the tractions have units of inverse length.

Invoking the Lorentz Reciprocal Theorem (LRT) \cite{BoundaryIntegral_Pozrikidis} and eliminating boundary terms arise from integration-by-parts using periodic BCs, we obtain
\begin{equation}
\int_{\Gamma} \vf{v}^{(i)} \cdot \vf{\lambda}^{(j)} \diff S = \int_{\Gamma} \vf{v}^{(j)} \cdot \vf{\lambda}^{(i)} \diff S. \label{appendix-LRT}
\end{equation}
After substituting the no-slip BC for $\vf{v}^{(j)}$ on the RHS of \eqref{appendix-LRT}, and make use of the force and torque balance condition for $\vf{\lambda}^{(i)}$, we obtain that
\begin{equation}
\begin{aligned}
\int_{\Gamma} \vf{v}^{(i)} \cdot \vf{\lambda}^{(j)} \,dS &= \int_{\Gamma} \left[ \vf{u}^{(j)} + \vf{\omega}^{(j)} \times (\vf{x}-\vf{q}) \right] \cdot \vf{\lambda}^{(i)} \diff S_{\vf{x}} \\
 &= \vf{u}^{(j)} \cdot \vf{f}^{(i)} + \vf{\omega}^{(j)} \cdot \vf{\tau}^{(i)}  \\ 
 &=  \vf{U}^{(j)} \cdot \vf{e}^{(i)} = \mathcal{N}_{ij}.
\end{aligned}
\label{appendix-LRT2}
\end{equation}
We recall that  $\vf{v}^{(i)}(\vf{x})$ for $\vf{x} \in \Gamma$ can be written as
\begin{equation}
\vf{v}^{(i)}(\vf{x}) = \int_{\Gamma} \pStokeslet(\vf{x} - \vf{y}) \cdot \vf{\lambda}^{(i)}(\vf{y})  \diff S_{\vf{y}}, \quad \vf{x} \in \Gamma, \label{appendix-firstkind}
\end{equation}
where $\pStokeslet$ is the Green's function of the Stokes equation with unit viscosity and periodic BCs (periodic Stokeslet), as dictated by the first-kind boundary integral formulation of the mobility problem \cite{BoundaryIntegral_Pozrikidis}.
Lastly, we substitute \eqref{appendix-firstkind} for $\vf{v}^{(i)}$ on the LHS of \eqref{appendix-LRT2} to conclude that
\begin{equation}
\mathcal{N}_{ij} = \int_{\Gamma} \int_{\Gamma} \vf{\lambda}^{(i)}(\vf{x}) \cdot \pStokeslet(\vf{x} - \vf{y}) \cdot \vf{\lambda}^{(j)}(\vf{y}) \diff S_{\vf{y}} \diff S_{\vf{x}}. 
\end{equation}


\section{A completed second-kind formulation in two dimensions}
\label{appendix:secondkind-formulation}
We present a completed second-kind boundary integral formulation of the deterministic Stokes boundary value problem (2.13), as used to compute highly-accurate reference results in the main body of the paper (see Sec.~4). We reproduce this formulation here for the benefit of the reader, since the full formulation is not contained in published work to our knowledge. Using a {\it double-layer} integral representation of the Stokes flow \cite{BoundaryIntegral_Pozrikidis}, we can write the exterior flow velocity $\vf{v}(\vf{x} \in E)$ in terms of an unknown double-layer density $\vf{\varphi}(\vf{x})$ as 
\begin{equation}
 \vf{v}(\vf{x} \in E) = (\op{D}_{\Gamma} \vf{\varphi})(\vf{x}) + (\op{G}[\vf{q}] \vf{f})(\vf{x}) + (\op{R}[\vf{q}]{\tau})(\vf{x}), \label{DLrepresentation}
 \end{equation}
where $\op{D}_{\Gamma}$ is the double-layer integral potential, $\op{G}[\vf{q}]$ and $\op{R}[\vf{q}]$ are the Stokeslet and rotlet at $\vf{q}$, respectively. In two dimensions, these operators are defined as
\begin{align}
(\op{D}_{\Gamma} \vf{\varphi})_j (\vf{x}) &= \frac{1}{4\pi \eta} \int_{\Gamma} T_{jlm}(\vf{x} - \vf{y} + \mbf{p}) n_m(\vf{y})  \varphi_l(\vf{y}) \, \diff S_{\vf{y}}, \\
(\op{G}[\vf{q}] \vf{f} )_j(\vf{x}) &=  \frac{1}{4\pi \eta} G_{jl}( \vf{x} - \vf{q} ) f_l, \\
(\op{R}[\vf{q}] {\tau})_j(\vf{x}) &=  \frac{1}{4\pi\eta} {R}_j( \vf{x} - \vf{q}) \tau \label{pRotlet},
\end{align}
where $G_{jl}$, $T_{jlm}$ and $R_j$ are the free-space Stokeslet, stresslet, and rotlet with unit viscosity, respectively. In two dimensions,
\begin{align}
 G_{jl} = -\delta_{jl} \log r + \frac{ x_j x_l}{r^2}, \quad
 T_{jlm} = -\frac{4 x_j x_l x_m}{r^4}, \quad
 R_{j}   = \frac{x_j^{\perp}}{r^2},
\end{align}
where $\vf{x}^{\perp} = (x_2, -x_1)$. It is well-known that the double-layer representation by itself cannot represent a flow that exerts force  and torque. Therefore, a completion flow is generally required to complete the representation. The choice of completion flow is not unique \cite{StokesSecondKind_Shelley}, and for simplicity, we use the completion flow proposed by Power and Miranda \cite{DoubleLayer_Power1987}, that is generated by a Stokeslet with applied force $\vf{f}$ at the tracking point $\vf{q}$ inside the particle, denoted as $\op{G}[\vf{q}]\vf{f}$, and a rotlet with torque $\tau$ (a scalar quantity in two dimensions), denoted as $\op{R}[\vf{q}] {\tau}$. We note that it is possible to avoid this kind of completion and still get a second-kind integral equation for the mobility problem \cite{Stokes2D_Manas,Stokes3D_FMM}.

Taking the limit as $\vf{x}$ approaches the boundary $\Gamma$ and using the jump condition for the double-layer potential \cite{BoundaryIntegral_Pozrikidis}, we obtain a second-kind boundary integral equation of the unknown translational and angular velocities $\{ \vf{u},\omega \}$ and the double-layer density $\vf{\varphi}$,
\begin{equation}
\vf{u} + \omega (\vf{x} - \vf{q})^{\perp} = \frac{1}{2} \vf{\varphi}(\vf{x}) + (\op{D}_{\Gamma} \vf{\varphi})(\vf{x}) +
 \op{G}[\vf{q}]\vf{f} + \op{R}[\vf{q}]\tau, \quad \forall \vf{x} \in \Gamma, \label{second-kind-formulation}
\end{equation}
and \Cref{second-kind-formulation} is closed by relating $\{ \vf{u},\omega \}$ to $\vf{\varphi}$ via\footnote{These are the two dimensional analogs of \cite[Eqs.~(4.9.22)-(4.9.23)]{BoundaryIntegral_Pozrikidis}. }
\begin{equation}
 \vf{u} = \frac{1}{|\Gamma|} \int_{\Gamma} \vf{\varphi}(\vf{x}) \, \diff S_{\vf{x}}, \quad
\omega = \frac{1}{W} \int_{\Gamma} \vf{\varphi}(\vf{x}) \cdot (\vf{x}-\vf{q})^{\perp} \, \diff S_{\vf{x}} \label{q2torque}, 
\end{equation}
 where $|\Gamma|$ is the length of $\Gamma$ and $W = \int_{\Gamma} |\vf{x}-\vf{q}|^2 \diff S_{\vf{x}}$. We remark that the random surface velocity $-\vslip (\vf{x})$ can also be included on the left-hand-side of \Cref{second-kind-formulation}. This provides us a mixed first- and second-kind formulation for solving the SSBVP (see Eq.~(2.10) of the paper). That is, we can generate the random surface velocity using the first-kind formulation described in the main body of the paper, and then, solve \Cref{second-kind-formulation,q2torque} for $\{ \vf{u},\omega\}$ to generate the action of $\op{N}$ and $\op{N}^{\frac{1}{2}}$ using a well-conditioned spectrally-accurate second-kind formulation. In this work we choose to use  the first-kind formulation for generating both the deterministic motions and the Brownian displacements in the main body of the paper, because the discrete fluctuation-dissipation balance is exactly preserved in the first-kind formulation, but not in the mixed formulation.

From a numerical standpoint, a desirable feature of the second-kind formulation in two dimensions is that the double-layer kernel (stresslet) has no singularities for points on the boundary, and discretizing the second-kind boundary integral equations using the regular trapezoidal rule leads to spectral accuracy. However, this holds neither in three dimensions nor for bodies that are near contact.

After discretizing \Cref{second-kind-formulation,q2torque}, we are required to evaluate the following periodic sums involving the Stokeslet, stresslet and rotlet in the resulting linear system,
\begin{align}
u^{G}_j(\mbf{x}_t) &= \frac{1}{4\pi \eta} \sum_{\mbf{p} \in \mathbb{Z}^2}  G_{jl}( \mbf{x}_t - \vf{q} + \mbf{p}L ) f_l, \label{pStokeslet-sum}\\
u^{R}_j(\mbf{x}_t) &= \frac{1}{4\pi\eta} \sum_{\mbf{p} \in  \mathbb{Z}^2}  {R}_j( \mbf{x}_t - \vf{q} + \mbf{p}L) \tau \label{pRotlet-sum},\\
u^{T}_j(\mbf{x}_t) &=  \frac{1}{4\pi \eta} \sum_{\mbf{p} \in\mathbb{Z}^2}  \sum_{s=1}^{N} T_{jlm}(\mbf{x}_t - \mbf{x}_s + \mbf{p}L) S_{lm}(\mbf{x}_s),  \label{pDoubleLayer-sum}
\end{align} 
where $\mbf{x}_t, \mbf{x}_s \in \Gamma$ denote the target and source points, respectively, and $S_{lm}(\mbf{x_s}) =  \varphi_{l}(\mbf{x}_s) n_{m}(\mbf{x}_s) \Delta s$.
We employ the Hasimoto function (see Eq.~(3.13) of the paper) to decompose the Stokeslet, rotlet and stresslet into near-field and far-field contributions, and thereafter, these sums can be efficiently computed using the Spectral Ewald method \cite{SpectralEwald_Stokes,BoundaryIntegral_Periodic3D}.
The decomposition of the periodic Stokeslet has been presented in the main body of the paper (see Eqs.~(3.14) and (3.15) of the paper).
The periodic sum of the rotlet can be decomposed as
\begin{equation}
 u^{R}_j(\mbf{x}_t) = \frac{1}{4\pi \eta} \sum_{\mbf{p} \in  \mathbb{Z}^2}  R^{(r)}_j(\mbf{x}_t-\vf{q} + \mbf{p}L ; \xi) \tau \nonumber + \frac{1}{\eta V} \sum_{\mbf{k} \neq \mbf{0}} {R}^{(w)}_j(\mbf{k} ; \xi)  \sin(\mbf{k}\cdot (\mbf{x}_t - \vf{q})) \tau ,
\end{equation}
where
\begin{subequations}
 \begin{align}
 R^{(r)}_j(\vf{x}; \xi) &= \frac{x^{\perp}_j}{r^2} (1-\xir) e^{-\xir}, \\
 {R}^{(w)}_j(\mbf{k}; \xi) &= \frac{2\pi k^{\perp}_j}{k^2}\left(1+\frac{k^2}{4\xi^2}\right)  e^{-\kxi}.
 \end{align}
 \end{subequations}
The complete decomposition of the periodic double-layer potential consists of four parts: the mean stresslet term $T^{(0)}_{jlm}$, which defines a zero mean flow in periodic domains \cite{BoundaryIntegral_Periodic3D}, the near-field $T^{(r)}_{jlm}$ and far-field  $T^{(w)}_{jlm}$ contributions,  and the term that arises from the diagonal elements of double-layer integral when $\mbf{x}_t = \mbf{x}_s$ (see \cite{StokesFastSolver_Biros}). In summary, we have
\begin{align}
  u^{T}_j(\mbf{x}_t) = \frac{1}{4\pi \eta} &\left(  \sum_{s=1}^{N} T^{(0)}_{jlm}(\mbf{x}_s) S_{lm}(\mbf{x}_s) \right. \nonumber  \\ 
  &~+ \sum_{\mbf{p} \in \mathbb{Z}^2} \sum_{n=1}^{N} T_{jlm}^{(r)}(\mbf{x}_t-\mbf{x}_s+\mbf{p}L; \xi)S_{lm}(\mbf{x}_s)  \nonumber \\
  &~+ \frac{1}{V} \sum_{\mbf{k} \neq \mbf{0}} {T}^{(w)}_{jlm}(\mbf{k} ; \xi) \sum_{n=1}^N S_{lm}(\mbf{x}_s)  \sin(\mbf{k} \cdot (\mbf{x}_t - \mbf{x}_s)), \nonumber  \\
  &+ \left. 2\pi \sum_{s=1}^N \kappa(\mbf{x}_s) \left[ (\vf{\hat{t}} \otimes \vf{\hat{t}})(\mbf{x}_s) \right]_{jl} \, \varphi_l(\mbf{x}_s)  \right),
\end{align}
 where $\kappa$ and $\vf{\hat{t}}$ are the curvature and unit tangent vector on $\Gamma$, and \cite{Stokes2D_VanDeVorst}
\begin{subequations}
\begin{align}
 T^{(0)}_{jlm}(\vf{x}) &= \frac{4\pi}{V}\delta_{lm} x_j, \\
 T^{(r)}_{jlm}(\vf{x} ; \xi) &= \left[ -\frac{4x_j x_l x_m}{r^4}(1+\xir) + 2\xi^2(\delta_{lm}x_j + \delta_{mj}x_l)\right] e^{-\xir},\\
 {T}^{(w)}_{jlm}(\mbf{k} ; \xi) &= -\frac{4\pi}{k^4} \left[ \left(1+\frac{k^2}{4\xi^2}\right) \left[k^2(k_j\delta_{lm}+k_l\delta_{mj}) - 2k_j k_l k_m \right] + k^2 k_m \delta_{jl} \right] e^{-\kxi}.
 \end{align}
 \end{subequations}
After we have efficient routines to compute the periodic sums described above, the discrete linear system obtained from \Cref{second-kind-formulation,q2torque} can be solved iteratively using GMRES without preconditioning, since the second-kind formulation is well-conditioned.








\end{document}